\documentclass{amsart}
\usepackage{amssymb, amsmath, amsfonts, amscd}

\input xypic

\theoremstyle{plain}
\newtheorem{theorem}{Theorem}[section]
\newtheorem{corollary}[theorem]{Corollary}
\newtheorem{lemma}[theorem]{Lemma}
\newtheorem{proposition}[theorem]{Proposition}

\newtheorem{example}[theorem]{Example}

\theoremstyle{definition}
\newtheorem{definition}[theorem]{Definition}

\theoremstyle{remark}

\numberwithin{equation}{theorem}

\newcommand{\M}{\mathcal{M}}

\newcommand{\E}{\mathcal{E}}

\renewcommand{\O}{\mathcal{O} }

\newcommand{\Hom}{\operatorname{Hom}} 
\newcommand{\T}{\operatorname{T}} 

\newcommand{\Der}{\operatorname{Der} }
\newcommand{\End}{\operatorname{End} }
\newcommand{\Spec}{\operatorname{Spec} }
\renewcommand{\P}{\operatorname{P} }
\renewcommand{\H}{\operatorname{H} }
\newcommand{\D}{\operatorname{D}^1_f(A) }
\newcommand{\Dl}{\operatorname{D} }
\newcommand{\DO}{\operatorname{D} }

\newcommand{\R}{\operatorname{R} }
\newcommand{\U}{\operatorname{U}} 
\newcommand{\Uo}{\operatorname{U}^{\otimes} }
\newcommand{\Ur}{\operatorname{U}^{\rho} }
\newcommand{\tUo}{\tilde{\operatorname{U}}^{\otimes} }
\newcommand{\tUr}{\tilde{\operatorname{U}}^{\rho} }
\newcommand{\K}{\operatorname{K} }
\renewcommand{\lg}{\mathfrak{g}}
 
\newcommand{\Diff}{\operatorname{Diff}}

\newcommand{\Z}{\operatorname{Z} }
\newcommand{\tL}{\tilde{L} }
\newcommand{\ta}{\tilde{\alpha} }
\newcommand{\tp}{\tilde{\pi}}
\newcommand{\sym}{\operatorname{Sym} }

\newcommand{\SETS}{\underline{\operatorname{Sets} }}
\newcommand{\LR}{\underline{\text{LR}(A/k)}}
\newcommand{\ALG}{\underline{\operatorname{Alg}}(k,A\otimes_k A)}
\newcommand{\ATiA}{\underline{\operatorname{D}^1_f(A)-\text{Lie}} }

\newcommand{\Conn}{\underline{\operatorname{Conn}} }
\newcommand{\Mod}{\underline{\operatorname{Mod} }}
\newcommand{\Ext}{\operatorname{Ext}}
\newcommand{\Tor}{\operatorname{Tor}}
\renewcommand{\SS}{\operatorname{SS} }

\newcommand{\HH}{\operatorname{HH}}
\newcommand{\C}{\mathbb{C} }

\newcommand{\fal}{\alpha^*(f) }

\newcommand{\falp}{\alpha'^*(f)  }

\begin{document}

\title{The enveloping algebra of a Lie algebra of differential operators}

\author{Helge \"{O}ystein Maakestad}

\email{\text{h\_maakestad@hotmail.com}}

%\thanks{"One Ring to rule them all, One Ring to find them, One Ring to bring them all and in the darkness bind them."  }


\date{March 2019}

\begin{abstract} The aim of this note is to introduce the notion of a $\Dl$-Lie algebra and to prove some general properties of the category of $\Dl$-Lie algebras, connections on $\Dl$-Lie algebras,
and universal enveloping algebras of $\Dl$-Lie algebras. We also define cohomology and homology of a connection on a $\Dl$-Lie algebra. One consequence of the construction is a functorial definition of Ext and Tor groups 
of arbitrary pairs $V,W$ of non-flat connections on an arbitrary $A/k$-Lie-Rinehart algebra $(L,\alpha)$. 
A $\Dl$-Lie algebra $\tL$ is a Lie-Rinehart algebra over $A/k$ equipped with an $A\otimes_k A$-module structure and a canonical central element $D\in Z(\tL)$ satisfying a compatibility property with the Lie-structure. Given a $\Dl$-Lie algebra $\tL$ and a connection $(\rho, E)$ we construct
the universal enveloping ring $\tUo(\tL,\rho)$ of $(\rho, E)$. The associative unital ring $\tUo(\tL,\rho)$ is a quotient of the associative ring  $\tUo(\End(\tL,E))$ corresponding to the non-abelian extension $\End(\tL,E)$ of the $\Dl$-Lie algebra $\tL$, and is a sub ring of $\Diff(E)$ - the ring of differential operators on $E$.  In the case when $A$ is Noetherian and $E$  and $\tL$ are finitely generated as left $A$-modules
it follows the ring $\tUo(\tL,\rho)$ is an almost commutative Noetherian ring. The ring $\tUo(\tL,\rho)$ is a quotient of the associative ring  $\Uo(\End(\tL,E))$ of the non-abelian extension $\End(\tL,E)$ and $\Uo(\End(\tL,E))$ is non-noetherian in general. If $E$ is a finitely generated $A$-module it follows the non-flat connection $(\rho, E)$ is a finitely generated $\tUo(\tL,\rho)$
module, hence we may speak of the characteristic variety $\SS(\rho,E)$ of $(\rho, E)$ in the sense of $D$-modules. We may define the notion of holonomicity for non-flat connections using the universal ring $\tUo(\tL,\rho)$. 
This was previously done for flat connections.
\end{abstract}

\maketitle

\tableofcontents

\section{Introduction}  Let $k$ be a  commutative unital ring and let $A$ be a commutative unital $k$-algebra. A $(k,A)$-Lie-Rinehart algebra $L$ is a $k$-Lie algebra and left $A$-module $L$ equipped with a map $a: L \rightarrow \Der_k(A)$ of $A$-modules and $k$-Lie algebras such that 

\[  [x,uy]=u[x,y]+ a(x)(u)y \]
for all $u \in A$ and $x,y \in L$. A connection on $L$ is a pair $(E, \nabla) $ where $E$ is a left $A$-module and $\nabla: L \rightarrow \End_k(E)$ is a $k$-linear map satisfying the derivation property.  There is a universal enveloping algebra $U(A,L)$ of $L$ with the property there is an "equivalence of categories" between the category of flat connections on $L$ and the category of left modules on $U(A,L)$.  The associative ring $U(A,L)$ has a canonical filtration and when $L$ is projective as left $A$-module there is a Poincare-Birkhoff-Witt isomorphism

\[ \rho: \sym_A^*(L) \cong Gr(U(A,L)) \]
where the right hand side is the associative graded ring of $U(A,L)$ with respect to the canonical filtration. In \cite {rinehart} Rinehart uses the isomorphism $\rho$ and the algebra $U(A,L)$ to construct cohomology and homology groups
of flat connections on $L$, and to prove basic properties of these cohomology groups. 

Given any 2-cocycle $c\in \Z^2(L,A)$ where $\Z^2(L,A):=ker(d^2)$ and $d^2$ is the differential in the Lie-Rinehart complex of $L$. A connection $(E,\nabla)$ has "curvature type $c$" if and only if the following holds:

\[ R_{\nabla}(x,y)(e)=c(x,y)e \]
for any $e\in E$ and $x,y \in L$. If $c\neq 0$ it follows $\nabla$ is a non-flat connection on $E$. There is a generalized universal enveloping algebra $U_c(A,L)$, with the property there is an "equivalence of categories" between the category of 
left $U_c(A,L)$-modules and the category of connections of curvature type $c$. The algebra $U_c(A,L)$ is a construction originating in a paper of Sridharan \cite{sridharan} and it has been studied by several people. When $c=0$ we get Rineharts universal enveloping algebra, hence $U_c(A,L)$ may be viewed as a family of associative rings parametrized by $\Z^2(L,A)$. The associative ring $U_c(A,L)$ has a canonical filtration, and  when $L$ is projective there is for any $c$ a PBW-isomorphism
\[ \rho_c: \sym_A^*(L) \cong Gr(U_c(A,L)) \]
generalizing the PBW-isomorphism of Rinehart (see \cite{maa1}).

 Let $(V,\nabla)$ and $(W,\nabla')$ be two arbitrary $L$-connections. Let $\Conn(L)$ denote the abelian category of $L$-connections and morphisms of connections. 
A fundamental problem in the study of $\Conn(L)$ is to give an explicit construction of the "Yoneda-Ext" group  
$\Ext^i((V,\nabla) ,(W,\nabla'))$ parametrizing equivalence classes of extensions of connections $(V,\nabla)$ and $(W,\nabla')$ of length $i$ in the category $\Conn(L)$:

\begin{align}
&\label{exti}  0 \rightarrow (W,\nabla') \rightarrow (V_1,\nabla_1) \rightarrow \cdots \rightarrow (V_i, \nabla_i) \rightarrow (V, \nabla) \rightarrow 0 .
\end{align}
A map between the connections $V$ and $W$ is an $A$-linear map $\phi: V \rightarrow W$ commuting with the action from $L$. One may check directly that the equivalence classes  of exact sequences of the type 
given in \ref{exti} form an abelian group. There is in general no way to given an abstract definition of Tor groups $\Tor_i((V,\nabla),(W,\nabla))$ for two arbitrary $L$-connections $V,W$ in $\Conn(L)$.
Since $\Conn(L)$ is a small abelian category, the Freyd-Mitchell Full Embedding Theorem (see \cite{freyd}) says there is an associative ring $R$ and an equivalence $\phi$ between $\Conn(L)$ and a subcategory of $\Mod(R)$. The equivalence
$\phi$ does not preserve injective and projective objects, hence we cannot use $\phi$ to define Ext and Tor groups of non-flat connections. 

In this paper we introduce the notion of a $\Dl$-Lie algebra $\tL$ and the universal
ring $\Uo(\tL)$ of  $\tL$. Using this construction we prove the following: Let $(L,a)$ be an arbitrary Lie-Rinehart algebra and let $c \in \Z^2(\Der_k(A),A)$ be any 2-cocycle. 
Let $L_c$ be the $\Dl$-Lie algebra associated to the pair $L$ and $c$. The construction of $L_c$ is functorial in $L$  (see Theorem \ref{main}). We prove there is an exact equivalence of categories

\begin{align}
&\label{equivconn}  \psi_c: \Conn(L) \cong \Mod(\Uo(L_c)) 
\end{align}
preserving injective and projective objects, where $\Uo(L_c)$ is the universal ring of the $\Dl$-Lie algebra $L_c$.
Hence we may use the associative unital ring $\Uo(L_c)$ to define the cohomology and homology of any pair of $L$-connections $(V,\nabla)$ and $(W,\nabla')$, flat or non-flat.  The functor $\psi_c$
realize the category of connections $\Conn(L)$ as a module category over an associative ring $\Uo(L_c)$ for any 2-cocycle $c\in \Z^2(L,A)$.
Using the equivalence $\psi_c$ we get for any integer $i \geq 0$ isomorphisms

\begin{align}
&\label{isoconn}     \Ext^i((V,\nabla),(W,\nabla')) \cong  \Ext^i_{\Uo(L_c)}(\psi_c(V,\nabla), \psi_c(W,\nabla'))) 
\end{align}
where the left side is the "Yoneda-Ext"-group. Hence the universal ring $\Uo(L_c)$ may be used to calculate the "true" $\Ext$-group of an arbitrary pair of $L$-connections $(V,\nabla)$ and $(W, \nabla')$, where $(L,a)$ is an arbitrary
Lie-Rinehart algebra. In the paper \cite{rinehart} this construction was done for flat $L$-connections. If $k$ is a field there is an isomorphism

\[    \Ext^i((V,\nabla),(W,\nabla'))  \cong \HH^i(\Uo(L_c), \Hom_k(V, W)) \]
where the right hand side is the Hochschild cohomology of the $\Uo(L_c)$-bimodule $\Hom_k(V,W)$. Hence the "Yoneda-Ext"-groups are calculated by the Hochschild complex.

%Note: The associative ring $R$ from the Freyd-Mitchell Full Embedding Theorem cannot be used 
%since the functor $\phi$ does not preserve injective and projective objects in general. 

The ring $\Uo(\tL)$ is almost commutative in general as is the case for $\U_c(A,L)$. When $A$ is Noetherian and $L$ a finite rank projective $A$-module it follows $\U_c(A,L)$ and $\Uo(\tL)$ are Noetherian in general. 
When the base ring $k$ is a field, it follows the groups $\Ext^i_{\Uo(\tL)}(V,W)$  and $\Tor_i^{\Uo(\tL)}(V,W)$ may be calculated using Hochschild cohomology 
and homology for modules on an associative unital ring. Since the groups 
\[ \Ext_{\Uo(\tL)}^i(V,W), \Tor_i^{\Uo(\tL)}(V,W) \]
are defined as Ext and Tor groups of modules over the associative unital ring $\Uo(\tL)$, they satisfy the usual functorial properties of such groups. It is not possible to prove such functorial properties working in 
the abelian category $\Conn(\tL)$ and with the construction in this paper all such properties are immediate. Hence the introduction of the notion $\Dl$-Lie algebra, the universal ring $\Uo(\tL)$ and the equivalence of categories
in \ref{equivconn} solves the problem of constructing cohomology and homology groups of arbitrary non-flat connections $V,W$ on an arbitrary Lie-Rinehart algebra. We get a definition of Ext and Tor groups of a pair of connections $(V,W)$ valid in complete generality. All functorial properties of the Ext and Tor groups follows from the classical book of Cartan and Eilenberg \cite{cartan}. It is impossible to give a general construction
of the Tor group $\Tor_i((V,\nabla),(W,\nabla))$ without the  equivalence of categories given in \ref{equivconn}.

 The Riemann-Hilbert correspondence in it's most naive form is a relation between the category of finite rank vector bundles with a flat connection on a simply connected complex projective manifold $X$ and finite dimensional complex representations of the topological fundamental group $\pi_1(X,x) $ of $X$. Given a flat connection $(E, \nabla)$, it follows the kernel $\ker(\nabla)$ of the connection $\nabla$ is a local system of finite dimensional complex vector spaces on $X$. The local system $ker(\nabla)$ gives rise to a finite dimensional complex representation $\rho: \pi_1(X,x) \rightarrow \operatorname{GL}(V)$. The Riemann-Hilbert correspondence says this is an "equivalence of categories".  A vector bundle with a flat connection is an algebraic object and a representation of $\pi_1(X,x)$ is a topological object. Hence this correspondence relates the topology of $X$ to the algebraic geometry of $X$ since $X$ is algebraic and since any finite rank vector bundle on $X$ is algebraic.  A vector bundle $E$ with a flat connection is canonically a left module on the sheaf of differential operators $D_X$ of $X$. There is a generalization of the sheaf $D_X$. A sheaf of generalized differential operators $A_X$ on $X$, is a sheaf of associative rings on $X$ with the property
that there is an open cover $U_i$ of $X$ and isomorphisms $D_{U_i} \cong (A_X)_{U_i}$ of sheaves of associative rings on $U_i$. Here $D_{U_i}$ and $(A_X)_{U_i}$ are the restrictions of the sheaves $D_X$ and $A_X$ to $U_i$. If $L$ is an invertible sheaf on $X$, it follows the sheaf $D_L$ of differential operators on $L$ is such a generalized sheaf of differential opertators on $X$. A generalized connection is a left module on $A_X$ and such connections have been studied by several people (see the papers of Bernstein, Beilinson and Simpson \cite{beilinson} and \cite{simpson}). If $U=\Spec(A) \subseteq X $ is an open affine subscheme of $X$ it follows $(A_X)_U \cong \U_c(A,L)$ for $L:=\Der(A)$ and some 2-cocycle $c$. Hence the sheaf $A_X$   is a global version of the associative ring $\U_c(A,L)$ studied in \cite{maa1}. The associative ring $\U_c(A,L)$ is a quotient of the ring $\Uo(L_c)$, hence the construction in this paper is related to the construction of Bernstein, Beilinson and Simpson.

%The aim of this note is to introduce and prove various general properties a generalization of a Lie-Rinehart algebra - a $\Dl$-Lie algebra.

A $\Dl$-Lie algebra $\tL$ is a Lie-Rinehart algebra over $A/k$ equipped with an $A\otimes_k A$-module structure that is compatible with the Lie-structure. There is a central element $D\in \tL$
satisfying a compatibility property with the Lie product. In the special case when $\tL$ as a left $A$-module is an abelian extension of $A$ by some 2-cocycle $f\in \Z^2(L,A)$ we may view $\tL$ as an Atiyah algebra
with an additional $A\otimes_k A$-structure. Hence $\tL$ with the underlying left $A$-module structure may be viewed as a simultaneous generalization of a Lie-Rinehart algebra and an Atiyah 
algebra with additional structure. We introduce the category $\ATiA$ of $\Dl$-Lie algebras, connections on $\Dl$-Lie algebras and prove various general properties of this construction. We also correct some mistakes in an earlier paper on this subject related to the  universal algebra $\U^{ua}(L(\fal))$  of a Lie-Rinehart algebra $(L,\alpha)$ (see \cite{maa1}, Appendix A).  

The main results in the paper are the following theorems: Given a $\Dl$-Lie algebra $(\tL,\ta, \tp,[,],D)$. Define the categories  of $\tL$-connections 

\[ \Mod(\tL, Id) \text{ and  }\Conn(\tL, Id),\]
 
as follows: An object $(E,\rho)$ in $\Mod(\tL, Id)$ is a left  $A$-module $E$ and an $A\otimes_k A$-linear map $\rho: \tL \rightarrow \End_k(A)$ with $\rho(D)=Id_E$. Morphisms in $\Mod(\tL, Id)$ are $A$-linear maps commuting with the action of $\tL$. An object in $\Conn(\tL, Id)$ is a pair $(E,\rho)$ where $E$ is a left $A$-module and and $A$-linear map $\rho: \tL \rightarrow \End_k(E)$ with
\[ \rho(x)(ae)=a\rho(x)(e)+ \tp(x)(a)\psi(e) \]
where $\psi\in \End_A(E)$ and and $\rho(D)=Id_E$. Morphisms in $\Conn(\tL, Id)$ are $A$-linear maps commuting with the action of $\tL$. We construct  two associative unital rings $\Uo(\tL)$ and $\Ur(\tL)$
with the following property (see Theorem \ref{univfunctor} and \ref{mainequiv}):
\begin{theorem} \label{main1} There are covariant functors 
\begin{align}
&\Uo: \ATiA \rightarrow \underline{Rings} \\
&\Ur: \ATiA \rightarrow \underline{Rings}.
\end{align}
with the following property: For any $\Dl$-Lie algebra $\tL$ there are exact equivalences of categories
\begin{align}
& F_1:\Mod(\tL, Id) \cong \Mod(\Uo(\tL)) \\
& F_2: \Conn(\tL, Id) \cong \Mod(\Ur(\tL))
\end{align}
with the property that $F_1$ and $F_2$ preserves injective and projective objects.
\end{theorem}

We use the associative rings $\Uo(\tL)$ and $\Ur(\tL)$ in Definition \ref{cohomology} to define the cohomology and homology of an arbitrary connection $(\rho,E)$. Previously
the notion of cohomology and homology was defined for flat connections. By Theorem \ref{main1} it follows the associative rings 
$\Uo(\tL)$ and $\Ur(\tL)$ may be viewed as \emph{universal enveloping algebras for non-flat connections}. The rings $\Uo(\tL)$ and $\Ur(\tL)$ are non-Noetherian in general. 
If $A$ is a Noetherian ring and the connection $(\rho,E)$ has the property that $E$ is a fintely generated $A$-module it follows from Proposition \ref{univnoetherian}
that the quotient ring $\Uo_E(\tL):=\Uo(\tL)/ann(\rho,E)$ is Noetherian.
A similar property holds for $\Ur(\tL)$. Hence even though the rings $\Uo(\tL)$ and $\Ur(\tL)$ are non-Noetherian in general, we may always pass to Noetherian quotients when studying connections $E$ that are
finitely generated as $A$-modules (see Example \ref{noetherianconn}).

We prove that the  rings $\Uo(\tL), \Ur(\tL)$ are solutions to universal problems in Proposition \ref{representconn}. Let $\ALG$ be the category with the following objects: Objects are
 associative unital $k$-algebras $R$ with $k$ in the centre of $R$ such that $R$ has a left $A\otimes_k A$-module structure. Maps in $\ALG$ are maps $f: R\rightarrow R'$ of unital $k$-algebras and left 
$A\otimes_k A$-modules. Let $J\subseteq \Uo(\tL)$ be a 2-sided ideal and let $\Uo_J(\tL):=\Uo(\tL)/J$.
Define the functor 
\[ \Conn_{\tL, J}: \ALG \rightarrow \SETS \]
by letting $\Conn_{\tL, J}(R)$ be the set of $A\otimes_k A$-linear maps $\rho: \tL \rightarrow R$ with $\rho(D)=1_R$ and $\Uo(\rho)(J)=0$. Here $\Uo(\rho): \Uo(\tL)\rightarrow \End_k(E)$ is the map induced by $\rho$.
There is a canonical structure of left $A\otimes_k A$-module on $\Uo(\tL)$
and a map of $A\otimes_k A$-modules $p^{\otimes}_J: \tL \rightarrow \Uo_J(\tL)$ with $p^{\otimes}(D)=1$.
It follows there is a functorial equality of sets
\[  \Hom_{k-alg}(\Uo_J(\tL), R) \cong \Conn_{\tL,J}(R) \]
hence the pair $(\Uo_J(\tL), p^{\otimes}_J)$ represents the functor $\Conn_{\tL,J}$. It follows the pair $(\Uo_J(\tL), p^{\otimes}_J)$ is unique up to unique isomorphism. A similar type result holds for $\Ur(\tL)$.
Hence it is justified to call the associative rings $\Uo(\tL)$ (and $\Uo_J(\tL)$) \emph{the universal ring of $\tL$}. In Example \ref{moduliconnection} we indicate how the functor $\Conn_{\tL,J}$ and the universal ring $\Uo(\tL)$  can be used to construct moduli spaces of arbitrary connections generalizing the classical case. The set
\[ \Conn_{\tL,J}(\End_k(E)) \]
is by definition the set of all connections $\rho:\tL \rightarrow \End_k(E)$ with $J$-curvature equal to zero. We say $\Conn_{\tL}:=\Conn_{\tL,(0)}$ where $(0)\subseteq \Uo(\tL)$ is the zero ideal, 
is the \emph{universal moduli functor for $\tL$-connections}  since
\[ \Conn_{\tL}(\End_k(E)) \]
by definition is the set of all $\tL$-connections $\rho: \tL\rightarrow \End_k(E)$.

In Example \ref{atiyah} we introduce the first order $\tL$-jet bundle $J^1_{\tL}(E)$ of a 
left $A$-module $E$, the $\Dl$-Atiyah sequence and the $\Dl$-Atiyah class 
\[ a_{\tL}(E) \in \Ext^1_{A\otimes_k A}(\tL\otimes_A E, E). \]
 The class $a_{\tL}(E)=0$ if and only if $E$ has an $(\tL, \psi)$-connection $\rho$.

To illustrate how the associative rings $\Uo(\tL)$ and $\Ur(\tL)$ can be used in the study of the classical curvature we construct in Example \ref{familyconn} the following: 
For any $A/k$-Lie-Rinehart algebra $(L,\alpha)$ and any 2-cocycle $f\in \Z^2(L,A)$ we construct a 2-sided ideal $I(f)\subseteq \Ur(L(0))$
where $L(0)$ is the abelian extension of $L$ with the zero cocycle. There is an equivalence of categories
\[ \Mod(\Ur(L(0))/I(f)) \cong \Mod(\U(A,L,f)) ,\]
where $\U(A,L,f)$ is the generalized universal enveloping algebra studied in \cite{maa1}. The associative ring $\U(A,L,f)$ has the property that left $\U(A,L,f)$-modules
correspond to $L$-connections of curvature type $f$. Hence any left $\Ur(L(0))$ module $(\tilde{\rho},E)$ annihilated by the ideal $I(f)$ corresponds to an $L$-connection $(\rho,E)$ with curvature type $f$. 
Hence we may use one fixed ring $\Ur(L(0))$ and the set of 2-sided ideals in $\Ur(L(0))$ to study the curvature $R_{\rho}$ of a connection $(\rho,E)$ on $L$. Hence the study of the set of
2-sided ideals in the rings $\Uo(\tL)$ and $\Ur(\tL)$ has applications in the study of the curvature of a classical connection. The algebra $\U(A,L,f)$ is a local version of a much studied object in 
the field $D$-modules.

The two associative unital rings $\Uo(\tL)$ and  $\Ur(\tL)$ are equipped with 2-sided ideals $I^{\otimes}\subseteq \Uo(\tL)$ and $I^{\rho}\subseteq \Ur(\tL)$ such that the following holds for the quotient
rings $\tUo(\tL):=\Uo(\tL)/I^{\otimes}$ and $\tUr(\tL):= \Ur(\tL)/I^{\rho}$ (see Theorem \ref{mainnoetherian}):

\begin{theorem} Let $(\tL, \ta, \tp,[,],D)$ be a $\Dl$-Lie algebra where $A$ is Noetherian and $\tL$ is finitely generated as left $A$-module. It follows the rings $\tUo(\tL)$ and $\tUr(\tL)$
are almost commutative unital Noetherian rings.
\end{theorem}

Hence we get many non-trivial examples of Noetherian quotients of the non-Noetherian rings $\Uo(\tL)$ and $\Ur(\tL)$.

Given an arbitrary $\Dl$-Lie algebra $(\tL, \ta,\tp, [,],D)$ and an arbitrary connection $(\rho,E)$ in $\Mod(\tL, Id)$ we may construct the non-abelian extension $\End(\tL,E)$ of $\tL$
by the $\tL$-connection $\End_A(E)$ as done in \cite{maa141}. We use this construction to construct the universal ring $\tUo(\tL, \rho)$ of the connection $(\rho,E)$. In Theorem \ref{mainuniversal}
we prove the following:

\begin{theorem} \label{mainuniversal} Let $(\tL, \ta,\tp,[,],D)$ be a $\Dl$-Lie algebra and let $(\rho, E)$ be an $\tL$-connection. There is a canonical map
\[ \rho^!:\End(\tL, E)\rightarrow \Diff^1(E) \]
and $\rho^!$ is a map of $B:=A\otimes_k A$-modules and $k$-Lie algebras. The map $\rho^!$ induce a map $\T(\rho^!):\tUo(\End(\tL,E))\rightarrow \Diff(E)$ of associative rings.
Let $\tUo(\tL,\rho):= Im(\T(\rho^!))$ be the image. We get an exact sequence of rings
\[ 0 \rightarrow ker(\T(\rho^!)) \rightarrow \tUo(\End(\tL, E)) \rightarrow \tUo(\tL, \rho) \rightarrow 0\]
where $\tUo(\End(\tL,E)):= \Uo(\End(\tL, E))/I$ where $I$ is the 2-sided ideal generated by the elements $u \otimes v-v\otimes u - [u,v]$ for $u,v\in \End(\tL,E)$. The rings
$\tUo(\End(\tL,E))$ and $\tUo(\tL,\rho)$ are almost commutative. If $A$ is noetherian and $\tL,E$ finitely generated as left $A$-modules it follows
$\tUo(\End(\tL,E))$ and $\tUo(\tL,\rho)$ are Noetherian rings.
\end{theorem}

Hence given an arbitrary connection $(\rho, E)$ on a $\Dl$-Lie algebra $(\tL,\ta,\tp,[,],D)$ we may construct the universal ring $\tUo(\tL,\rho)$ of $(\rho, E)$ and in the case when $A$ is Noetherian
and $\tL,E$ finitely generated $A$-modules it follows $\tUo(\tL,\rho)$ is an almost commutative Noetherian subring of $\Diff(E)$. The ring $\tUo(\tL, \rho)$ is defined for an arbitrary connection $\rho$ and one may use 
$\tUo(\tL,\rho)$ to define the characteristic variety $\SS(\rho,E)$ of $(E,\rho)$ and holonomiticy for non-flat connections. Previously notions such as holonomicity and characteristic variety have been studied for flat connections on holomorphic vector bundles on complex manifolds (see Example \ref{characteristic}).

\section{Functorial properties of $\Dl$-Lie algebras and connections}

In this section we introduce the notion of an $\Dl$-Lie algebra - a generalization of a Lie-Rinehart algebra. It is a Lie-Rinehart algebra equipped with the structure of 
an $A\otimes_k A$-module that is compatible with the Lie-structure.  Given any 2-cocycle $f\in Z^2(\Der_k(A),A)$ we construct in Theorem \ref{main} a functor

\[ F_f:  \LR \rightarrow \ATiA \]

from the category of $A/k$-Lie-Rinehart algebras to the category of $\Dl$-Lie algebras. We also consider connections on $\Dl$-Lie algebras and curvature of connections.

Let in the following $k \rightarrow A$ be an arbitrary map of unital  commutative rings. Let $\P^1:= \P^1_{A/k}$ be the module of principal parts. Its dual $\DO^1_0(A) :=\Hom_A(\P^1,A)=A\oplus \Der_k(A)$
has a canonical structure as a $k$-Lie algebra and $(A,A)$-module.
There is an inclusion $\DO^1_0(A) \subseteq \End_k(A)$ and we may define for any element $u=aI + x,v=bI+y$ with $I$ the identity endomorphism of $A$ and $x,y$ derivations and $a,b\in A$ the following:
\[ [u,v]:= u \circ v - v\circ u  = (x(b)-y(a))I+[x,y].\]
One checks this gives $\DO^1_0(A)$ the structure of a $k$-Lie algebra. It is the Lie algebra struture induced by the inclusion $\DO^1_0(A)\subseteq \End_k(A)$. Given any 2-cocycle $f\in Z^2(\Der_k(A),A)$ we may construct the following structure as $k$-Lie algebra on $A\oplus \Der_k(A)$:
\begin{align}
&\label{df} [aI+x,bI+y]=(x(b)-y(a)+f(x,y))I+[x,y]) \in A\oplus \Der_k(A).
\end{align}
The abelian group $A\oplus \Der_k(A)$ equipped with the $k$-Lie algebra structure $[,]$ in  \ref{df}  is denoted $\D$.  The abelian group $\D$ has two natural $A$-module structures:
\begin{align}
 cu&\label{m1}=c(aI+x)=(ca)I +cx \\
 uc&\label{m2}=(aI+x)c=(ac+x(c))I+cx.
\end{align}
The structures in \ref{m1} and \ref{m2} are induced by the left and right $A$-module structure on $\End_k(A)$. It follows $\D$ is an $A\otimes_k A$-module and a $k$-Lie algebra.
There is an endomorphism
\[ \pi: \D \rightarrow \D \]
defined by
\[ \pi(u)=\pi(aI+x)=x.\]
We get
\[ \pi([u,v])=\pi((x(b)-y(a)+f(x,y)I+[x,y])=[x,y]=[\pi(x),\pi(y)],\]
hence $\pi$ is a morphism of $k$-Lie algebras. We may think of $\pi$ as a map of $k$-Lie algebras and $A\otimes_k A$-modules
\[\pi:\D \rightarrow \Der_k(A) .\]
Here we give $\Der_k(A)$ the trivial right $A$-module structure. Let $\overline{z}:=(1,0)\in \D$. It follows $[\overline{z},u]=0$ for all elements $u\in \D$ hence $\overline{z}$ is a central element.

\begin{lemma} The following holds for every $u=aI+x,v=bI+y\in \D$ and $c\in A$: 
\begin{align}
&\label{lr} [u,cv]=c[u,v]+\pi(u)(c)v. \\
&\label{lr1}uc=cu+\pi(u)(c)\overline{z}
\end{align}
\end{lemma}
\begin{proof} we get
\[ [u,cv]=[aI+x, (cb)I+cy ] = (x(cb)-cy(a)+f(x,cy)I+[x,cy])=\]
\[ (cx(b)+x(c)b-cy(a)+cf(x,y))I+c[x,y]+x(c)y)=\]
\[  c((x(b)-y(a)+f(x,y))I+[x,y])+x(c)(bI+y)=c[u,v]+\pi(u)(c)v .\]
we get
\[ uc:=(aI+x)c:=(ac+x(c))I+cx=a(cI+x)+x(c)\overline{z}=cu+\pi(u)(c)\overline{z}\]
The Lemma follows.
\end{proof}

Let us sum this up in a Proposition:

\begin{proposition} Let $k\rightarrow A$ be a unital map of commutative rings and let $\D:=A\oplus \Der_k(A)$ with $f\in \Z^2(\Der_k(A),A)$ a 2-cocycle. Let $u:=aI + x, v:= bI+ y\in \D$.
Define the following left and right $A$-module structure on $\D$:
\[ c(aI+ x):= (ca)I+cx \]
and
\[ (aI+x)c:= (ac +x(c))I + c x\]
for $c\in A$. Define the following product $[,]$ on $\D$:
\[ [u,v]=[aI+ x,bI+ y]:=:=(x(b)-y(a)+f( x,y))I+ [x, y].\]
Define the map $\pi:\D \rightarrow \Der_k(A)$ by $\pi(aI+ x):=x$ and give $\Der_k(A)$ the trivial right $A$-module structure. It follows $\D$ is an $A\otimes_k A$-module and a $k$-Lie algebra. The map $\pi$ is a map of $A\otimes_k A$-modules and $k$-Lie algebras. The product $[,]$ satisfies
\[ [u,cv]=c[u,v]+\pi(u)(c)v \]
and
\[ uc=cu+\pi(u)(c)\overline{z}\]
for all $u,v \in \D$ and $c\in A$.
\end{proposition}
\begin{proof} The proof follows from the calculations above.
\end{proof}

Hence the underlying left $A$-module of the pair $(\D,\pi)$ is an ordinary Lie-Rinehart algebra. It is the abelian extension of $\Der_k(A)$ with the 2-cocycle $f\in \Z^2(\Der_k(A),A)$. 
It follows $\D$ is an Atiyah algebra in the sense of \cite{tortella1}.

Note: We see that it is impossible  to construct  a non-trivial right $A$-module structure on $\Der_k(A)$ induced by the inclusion $\Der_k(A) \subseteq \End_k(A)$. To get a non-trivial right $A$-module structure we must consider the abelian extension $\D$ for some 2-cocycle $f$.

%\begin{example} Special case: The dual of the module of principal parts. \end{example}

%If $f=0$ is the zero cocycle it follows $\operatorname{D^1(A,0)}:=\Hom_A(\P^1,A)$ is the dual of the module of principal parts $\P^1$ with its canonical $k$-Lie algebra structure and $A\otimes_k A$-module structure.
%The $A\otimes_ k A$-module structure and Lie algebra structure on $\operatorname{D^1(A,0)}$ comes from the structure on $\End_k(A)$. Hence $\D$ is a generalization of the dual of the module of principal parts $

%\Diff^1(A)$. One would like to give a definition of $\D$ in terms of $\P^1$: $\D=\Hom_A(\P^1_f,A)$ where $\P^1_f$ is some generalization of $\P^1$. One would also like to do a similar construction for higher order %differential operators $\Diff^l(A)$ and higher order principal parts $\P^l$.

%and to define "deformations of higher order connections" in terms "deformations of the principal parts".

We may define the following:

\begin{definition}\label{maindef}  A 5-tuple $(\tL,\ta,\tp, [,],D)$ where $\tL$ is an $A\otimes_k A$-module and $k$-Lie algebra  and 
\[ \ta: \tL \rightarrow \D \]
is a map of $A\otimes_k A$-modules and $k$-Lie algebras is a \emph{$\Dl$-Lie algebra} if the following holds: The element $D\in  Z(\tL)$ is a central element.
The map $\tp: \tL \rightarrow \Der_k(A)$ is a map of $A\otimes_k A$-modules and $k$-Lie algebras with $\tp(D)=0$ and $\pi \circ \ta=\tp$. Here $\Der_k(A)$ has the trivial right $A$-module structure.
For all $u,v\in \tL$ and $c\in A$ the following holds:
\begin{align}
&\label{al1} uc=cu +\tp(u)(c)D  \\
&\label{al3} [u,cv]=c[u,v]+\tp(u)(c)v 
\end{align}
A 4-tuple $(\tL,\tp,[,],D)$ (we remove the element $\ta$ from the definition of a $\Dl$-Lie algebra) satisfying the above criteria is a \emph{pre-$\Dl$-Lie algebra}.
Let $(\tL,\ta,\tp,[,],D)$ and $(\tilde{M}, \tilde{\beta}, \tilde{\gamma},[,],D')$ be $\Dl$-Lie algebras. A map $\psi: \tL \rightarrow \tilde{M}$ of $k$-Lie algebras and $A\otimes_k A$-modules is a  \emph{map of $\Dl$-Lie algebras} if 
$\tilde{\beta} \circ \psi = \ta$ and $\psi(D)=D'$.  Let $\ATiA$ denote the category of $\Dl$-Lie algebras and morphisms. Let $E$ be a left $A$-module. An $\tL$-\emph{connection} on $E$ is an 
$A\otimes_k A$-linear map
\[ \rho: \tL \rightarrow \End_k(E) .\]
The module $\End_k(E)$ has the $A\otimes_k A$-module structure defined by $a\otimes b \psi := a(\psi b)$ for $a\otimes b \in A\otimes_k A$ and $\psi \in \End_k(E)$.
Given two connections $(E,\rho_E)$ and $(F, \rho_F)$, a morphism $\phi: (E,\rho_E) \rightarrow (F, \rho_F)$ is a map of $A$-modules $\phi: E\rightarrow F$ 
such that for any element $u\in \tL$ we get a commutative diagram
\[
\diagram    E \dto^{\phi} \rto^{\rho_E(u)} &  E \dto^{\phi} \\
                  F \rto^{\rho_F(u) }                  &  F.
\enddiagram
\]
A $\Dl$-Lie algebra is also referred to as a \emph{Lie algebra of differential operators acting on $A/k$}.
Let $\ATiA$ denote the category of $\Dl$-Lie algebras and morphisms of $\Dl$-Lie algebras. Let $\Mod(\tL)$ denote the category of $\tL$-connections and morphisms of connections.
Let $\Mod(\tL,Id)$ denote the category of $\tL$-connections $(\rho,E)$ where $\rho(D)=Id_E$. We define similar notions for a pre-$\Dl$-Lie algebra.
\end{definition}

Note: By definition $(\D, id, \pi ,[,], z)$ where $\overline{z}:=(1,0)$ is a $\Dl$-Lie algebra.

Note: A $\Dl$-Lie algebra is an $A-A$-module in the sense of \cite{kontsevich} and non-commutative geometry and such objects are much studied in this field.

\begin{lemma} Let $(\tL, \ta, \tp,[,],D)$ be a $\Dl$-Lie algebra. The following formula holds for all $u,v\in \tL$ and $c\in A$:
\[  [u,vc]=  c[u,v]+\tp(u)(c)v+\tp(u)\tp(v)(c)D .\]
\end{lemma}
\begin{proof} We get
\[  [u,vc]=[u, cv+\tp(v)(c)D]=c[u,v]+\tp(u)(c)v+\tp(u)\tp(v)(c)D \]
and the Lemma is proved.
\end{proof}

\begin{lemma} Let $(L,\alpha)$ be an $A/k$-Lie-Rinehart algebra and let $f\in \Z^2(\Der_k(A),A)$ be a 2-cocycle. We get in a natural way a 2-cocycle $\fal\in \Z^2(L,A)$ defined by
$\fal(x,y):=f(\alpha(x),\alpha(y))$.
\end{lemma}
\begin{proof} Let $ u:=x\wedge y \wedge z\in \wedge^3\Der_k(A)$. We get
\[ d^2(\fal)(u)=d^2(f)(\alpha(x)\wedge \alpha(y) \wedge \alpha(z))=0 \]
since $d^2(f)=0$. The Lemma follows.
\end{proof}

Let $k \rightarrow A$ be a unital map of commutative rings and let $\alpha:L\rightarrow \Der_k(A)$ be a Lie-Rinehart algebra. Let $f\in \Z^2(\Der_k(A),A)$ be a 2-cocycle.
Let $L(\fal):=Az\oplus L$ with the following left and right $A$-module structure: Let $u:=az+x, v:=bz+y\in L(\fal)$.
\[ b(az+x):= (ba)z+bx \]
and
\[ (az+x)b:=(ab+\alpha(x)(b))z +ax.\]
Define the product $[u,v]$ as follows:
\[ [u,v]:=(\alpha(x)(b)-\alpha(y)(a)+\fal(x,y))z +[x,y] \in L(\fal).\]
Define the map $\alpha_f:L(\fal)\rightarrow \D$ by $\alpha_f(u):=aI+\alpha(x)$. Define $\pi_f:L(\fal)\rightarrow \Der_k(A)$ by $\pi_f(u)=\alpha(x)$. 
Assume $\phi: (L,\alpha)\rightarrow (L',\alpha')$ is a map of Lie-Rinehart algebras. Define the map $\phi_f:L(\fal) \rightarrow L'(\falp)$ by
\[ \phi_f(az+x):= az+\phi(x).\]
One checks that for any $u\in L(\fal)$ and $c\in A$ the following holds:
\[ uc=cu+\pi_f(u)(c)z \]
and that $[z,u]=0$ hence $z$ is a central element in $L(\fal)$.

\begin{lemma} \label{mainlemma} The abelian group $(L(\fal), [,])$ is an $A\otimes_k A$-module and
$k$-Lie algebra. The map $\alpha_f$ is a map of $A\otimes_k A$-modules and $k$-Lie algebras. The product $[,]$ satisfies 
\[ [u,cv]=c[u,v]+\pi_f(u)(c)v \]
and
\[ uc =cu+\pi_f(u)(c)z.\]
for all $u,v\in L(\fal)$ and $c\in A$. The map $\phi_f:L(\fal)\rightarrow L'(\falp)$ is a map of $A\otimes_k A$-modules and $k$-Lie algebras. There is an equality $\alpha'_f \circ \phi_f =\alpha_f$.
Hence the 5-tuple $(L(\fal), \alpha_f, \pi_f,[,], z)$ is a $\Dl$-Lie algebra.
\end{lemma}
\begin{proof} One checks $L(\fal)$ is an $A\otimes_k A$-module and $k$-Lie algebra and that $\alpha_f$ is a map of $A\otimes_k A$-modules and $k$-Lie algebras. One also checks that
\[ [u,cv]=c[u,v]+\pi_f(u)(c)v \]
for all $u,v\in L(\fal)$ and $c\in A$. Finally one checks that $\phi_f$ is a map of $A\otimes_k A$-modules and $k$-Lie algebras and that $\alpha'_f \circ \phi_f =\alpha_f$. The
Lemma follows.
\end{proof}

\begin{theorem} \label{main} Let $k \rightarrow A$ be an arbitrary map of unital commutative rings and let $f\in \Z^2(\Der_k(A),A)$ and let $\fal \in \Z^2(L,A)$ be the pull back of $f$ via $\alpha$. 
There is a covariant functor
\[ F : \LR \rightarrow \ATiA \]
defined by
\[ F(L,\alpha) := (L(\fal), \alpha_f, \pi_f,[,],z).\]
A map $\phi: (L,\alpha) \rightarrow (L',\alpha')$ of Lie-Rinehart algebras gives a map 
\[ \phi_f: L(\fal) \rightarrow L'(\falp) \]
of $\Dl$-Lie algebras. Hence $\phi_f(z)=z'$ with $z,z'$ the central elements of $L(\fal)$ and $L'(\falp)$.
The construction is functorial in the sense that if $\phi': (L', \alpha') \rightarrow (L", \alpha")$  is another map of Lie-Rinehart algebras it follows
\[ (\phi' \circ \phi)_f=\phi'_f \circ \phi_f .\] 
\end{theorem}
\begin{proof} The proof follows from Lemma \ref{mainlemma}.
\end{proof}

\begin{example} Atiyah algebras and $\Dl$-Lie algebras.\end{example}

Given a Lie-rinehart algebra $\alpha: L\rightarrow \Der_k(A)$ and a 2-cocycle $f\in \Z^2(\Der_k(A),A)$ it follows we get an exact sequence of $A\otimes_k A$-modules
\[ 0 \rightarrow Az \rightarrow L(\fal) \rightarrow L \rightarrow 0 \]

and the left $A$-module $L(\fal)$ is an \emph{Atiyah algebra} in the sense of \cite{ginzburg} if $L:=\Der_k(A)$. Hence $L(\fal)$ is an Atiyah algebra equipped with a canonical right
$A$-module structure and a marked central element $z\in L(\fal)$ satisfying.
\[  uc=cu +\pi_f(u)(c)z \]
for all $u\in L(\fal)$ and $c\in A$. A general $\Dl$-Lie algebra $\tL$ is not an extension of $\Der_k(A)$ by a rank one free $A$-module in general. Hence the underlying left $A$-module of 
$\tL$ is not an Atiyah algebra in general. 

Note: In \cite{tortella} Tortella constructs a canonical Hodge structure on the cohomology of the Atiyah algebra $At(\mathcal{L})$ of a holomorphic line bundle $\mathcal{L}$ on a complex projective manifold $X$. If 
$\mathcal{L}$ is a holomorphic line bundle on X and 
\[ 0\rightarrow \mathcal{L} \rightarrow At(\mathcal{L}) \rightarrow \Theta_X \rightarrow 0\]
an extension it follows $X$ has an open cover $\{U_i\}_{i \in I}$ where $\mathcal{L}$ trivialize. Hence we get exact sequences
\[ 0\rightarrow \O_{U_i} \rightarrow At(\O_{U_i}) \rightarrow \Theta_{U_i} \rightarrow 0 \]
where there is a 2-cocycle $f_i\in \Z^2(\Theta_{U_i},\O_{U_i})$ and an isomorphism $At(\O_{U_i})\cong \Theta_{U_i}(f_i)$ of sheaves of Lie-Rinehart algebras. Hence $At(\mathcal{L})$ is locally
the abelian extension of $\Theta_{U_i}$ by $\O_{U_i}$. Hence Tortella's Atiyah algebra is a global version of an abelian extension of Lie-Rinehart algebras valid for arbitrary complex manifolds.

\begin{definition} Let $(\tL,\ta,\tp,[,],D)$ be an $\Dl$-Lie algebra and let $(L,\alpha)$ be an $A/k$-Lie-Rinehart algebra. Let $E$ be an $A$-module and let $\psi\in \End_A(E)$. 
An $(\tL,\psi)-connection$ on $E$ is a map of left $A$-modules
\[ \nabla: \tL\rightarrow \End_k(E) \]
with the property that
\[ \nabla(u)(ae)=a\nabla(u)(e)+\tp(u)(a)\psi(e)\]
for all $u\in \tL, a\in A$ and $e\in E$.
An $(L,\psi)-connection$ on $E$ is a map of left $A$-modules
\[ \nabla: L\rightarrow \End_k(E) \]
with the property that
\[ \nabla(x)(ae)=a\nabla(x)(e)+\tp(x)(a)\psi(e)\]
for all $x\in L, a\in A$ and $e\in E$.
Let $\Conn(\tL, \End)$ be the category of $(\tL,\psi)$-connections and morphisms where we let $\psi\in \End_A(E)$ and $E$ vary.
Let $\Conn(\tL, Id)$ be the category of $(\tL,Id)$-connections $(E,\rho)$.
Let $\Conn(L,\End)$ denote the category of $(L,\psi)$-connections $\rho: L\rightarrow \End_k(E)$ where $\psi\in\End_A(E)$ may vary.
Let $\Conn(L,Id)$ denote the category of $(L,Id)$-connections $\rho: L\rightarrow \End_k(E)$.
\end{definition}

Note: It follows there are inclusions of categories 
\[  \Conn(\tL, Id) \subseteq \Conn(\tL) \subseteq \Conn(\tL,\End) \]

Recall that $\Conn(\tL)$ is the category of maps of $A\otimes_k A$-modules
\[ \rho: \tL \rightarrow \End_k(E) \]
and morphisms.

Note: The morphisms in $\Conn(\tL, \End)$ $\phi:(\rho, E)\rightarrow (\rho',F)$ are maps of left $A$-modules $\phi: E\rightarrow F$ such that
\[ \rho'(u) \circ \phi = \phi \circ \rho(u) \]
for all $u\in \tL$.
The morphisms $\phi:(\rho,E,\psi)\rightarrow (\rho',F, \psi')$ in $\Conn(L,\End)$ are by definition maps of $A$-modules
\[ \phi: E\rightarrow F \]
with the property that
\[  \psi' \circ \phi =\phi \circ \psi\text{ and } \rho'(x)\circ \phi = \phi \circ \rho(x) \]
for all $x\in L$.

\begin{lemma} A $\psi$-connection $\nabla:\tL\rightarrow \End_k(E)$ gives a map
\[ \nabla: \tL\rightarrow \Diff^1_k(E) \]
where $\Diff^1_k(E):=\Hom_A(\P^1\otimes_A E,E)$ is the module of first order differential operators of $E$.
\end{lemma}
\begin{proof} We must show that for any $b \in A$ and $x\in L$ it follows $[\nabla(x), bId_E]\in \End_A(E)$. We get
\[ [\nabla(x), bId_E] = \tp(x)(b)\psi \in \End_A(E) \]
and the Lemma follows.
\end{proof}

\begin{lemma} \label{corr} Let $(\tL, \ta,\tp,[,],D)$ be a $\Dl$-Lie algebra and let $E$ be a left $A$-module.
Any $\tL$-connection 
\[ \rho: \tL \rightarrow \End_k(E) \]
is an $(\tL,\psi)$-connection with $\psi:=\rho(D)\in \End_A(E)$.
It follows we get a map
\[ \rho:\tL \rightarrow \Diff^1(E) \]
of $A\otimes_k A$-modules.
\end{lemma}
\begin{proof} Assume $\rho: \tL \rightarrow \End_k(E)$ is an $A\otimes_k A$-linear map. We get for any $x\in \tL, a\in A$ and $e\in E$ the following calculation:
\[ \rho(x)(ae)=\rho(xa)(e)=\rho(ax+\tp(x)(a)D)(e)=a\rho(x)(e)+\tp(x)(a)\rho(D)(e).\]
We get
\[ \rho(D)(ae)=\rho(Da)(e)=\rho(aD+\tp(D)(a)D)(e)=\rho(aD)(e)=a\rho(D)(e) \]
since $\tp(D)=0$. It follows $\psi:=\rho(D)\in \End_A(E)$. It follows $\rho$ is a $\rho(D)$-connection with $\rho(D)\in \End_A(E)$. 
\end{proof}

\begin{lemma} \label{corr} Let $f\in \Z^2(\Der_k(A),A)$ and let $F(L,\alpha)=(L(\fal), \alpha_f, \pi_f,[,],z)$ be the $\Dl$-Lie algebra from Theorem \ref{main}. 
Let $E$ be a left $A$-module.
There is a one-to-one correspondence between $A\otimes_k A$-linear maps
\[ \rho: L(\fal) \rightarrow \End_k(E) \]
and $\psi$-connections $\nabla: L\rightarrow \End_k(E)$ with $\psi \in \End_A(E)$.
\end{lemma}
\begin{proof} Assume $\rho: L(\fal) \rightarrow \End_k(E)$ is an $A\otimes_k A$-linear map. Let $i: L\rightarrow L(\fal)$ be the canonical inclusion map
and define $\nabla:= \rho \circ i$. Since $\rho$ is $A\otimes_k A$-linear it follows $\nabla$ is left $A$-linear. Let $x\in L, a\in A$ and $e\in E$. We get
\[ \nabla(x)(ae)=\rho(i(x))(ae)=\rho(xa)(e)=\rho(\alpha(x)(a)z+ax)(e)=  \]
\[ a\rho(i(x))(e)+\alpha(x)(a)\rho(z)(e).\]
Put $\rho(z):=\psi$. It follows $\psi\in \End_A(E)$ and we get
\[ \nabla(x)(ae)=a\nabla(x)(e)+\alpha(x)(a)\psi(e) \]
hence $\nabla$ is a $\psi$-connection. Assume $\nabla: L\rightarrow \End_k(E)$ is a $\psi$-connection where $\psi\in \End_A(E)$. Define
\[ \rho: L(\fal) \rightarrow \End_k(E) \]
by
\[ \rho(u)=\rho(az+x):= a\psi +\nabla(x).\]
It follows $\rho$ is a left $A$-linear map. It is right $A$-linear for the following reason:
\[  \rho(uc):= \rho((az+x)c)=\rho((ac+\alpha(x)(c))z+ cx)=\]
\[ (ac+\alpha(x)(c))\psi+ \nabla(cx)=\]
\[ ac\psi+ c\nabla(x)+\alpha(x)(c)\psi =\]
\[ ac\psi+\nabla(x)c = a\psi c +\nabla(x)c=(a\psi +\nabla(x))c=\rho(u)c.\]
Hence
\[ \rho(uc)=\rho(u)c.\]

It follows $\rho$ is a map of $A\otimes_k A$-modules. This gives a one-to-one correpondence as claimed and the Lemma follows.
\end{proof}

Note: The definition of a connection $\rho: L(\fal)\rightarrow \End_k(E)$ depends on the $A\otimes_k A$-module structure on $L(\fal)$ and not on the Lie-algebra structure. Since
for any 2-cocycles $f,g\in \Z^2(\Der_k(A),A)$ it follows $L(\fal)\cong L(g^{\alpha})$ as $A\otimes_k A$-modules it follows there is a one-to-one correspondence between
$L(\fal)$-connections and $L(g^{\alpha})$-connections. In fact there is an equivalence of categories
\[ \Conn(L(\fal), \End) \cong \Conn(L(g^{\alpha}), \End) \]
for any pair of 2-cocycles $f,g\in \Z^2(\Der_k(A),A)$.

\begin{example} $\Dl$-Lie algebras, connections and $(A,A)$-vector bundles. \end{example}

Note: An ordinary $L$-connection $\nabla:L\rightarrow \End_k(E)$ corresponds by Lemma \ref{corr} to an $A\otimes_k A$-linear map

\[ \rho: L(\fal) \rightarrow \Diff^1_k(E) \]

with $\rho(z)=Id_E$. Usually a connection is a $k$-linear map
\[ \nabla: E \rightarrow \Omega^1_{A/k}\otimes_A E \]
satisfying $\nabla(ae)=a\nabla(e)+d(a)\otimes e$. Hence $ker(\nabla)$ and $Im(\nabla)$ are $k$-vector spaces. The vector spaces $ker(\nabla), Im(\nabla)$  are infinite dimensional in general.

If $A$ is a finitely generated and regular ring over  a field of characteristic zero, 
$E$ is finitely generated and projective as $A$-module and $L(\fal)=\D$ it follows $\Diff^1_k(E)$ and $L(\fal)$ are locally trivial $A$-modules of finite rank as left and right $A$-modules separately.
Hence the map $\rho$ is a more "geometric" object: One of the reasons to define a connection as a  map $\rho$ of $A\otimes_k A$-modules is because
we want to study the kernel $ker(\rho)$ and image $Im(\rho)$ and these modules are "geometric" objects since they are vector bundles from the left and right in many cases. 
With $\nabla$ the kernel $ker(\nabla)$ and image $Im(\nabla)$ are infinite dimensional  $k$-vector spaces and not  $A$-modules, and such objects are "more difficult" to handle. 

\begin{example} Left and right $A$-module structures on modules of principal parts.\end{example}

Note: An $A\otimes_k A$-module $W$ that is finitely generated and projective as left and right $A$-module separately is called an \emph{$(A,A)$-vector bundle}. 
The module of principal parts $\P^l(E)$ is an $(A,A)$-vector bundle in many cases. There are examples where the left structure on $\P^l(E)$ is different from the right structure (see \cite{maa15}).
Similar results hold for the module of differential operators $\Diff^l(E,E)$. From \cite{maa15} we get the following example. Let $C:=\mathbb{P}^1$ be the projective line over a field of characteristic zero and let 
$\O(n)$ be the invertible sheaf with $n\in \mathbb{Z}$ an integer. The module of $l$'th order differential operators $\Diff^l(\O(n))$ from $\O(n)$ to $\O(n)$ has a left and right structure as $\O_C$-module and we get the following classification:

\begin{theorem} Let $k\geq1$ and $n\in \mathbb{Z}$ be integers. The following holds:
\[ \Diff^l(\O(n))^{right}\cong \O_{C}\oplus \O(l+1)^k \text{ for all $l\geq 1$ and $n\in \mathbb{Z}$.}\]
\[ \Diff^l(\O(n))^{left}\cong  \O(l)^{l+1} \text{ for all $1 \leq l\leq n$.}\]
\[ \Diff^l(\O(n))^{left}\cong \O_{C}^{n+1}\oplus \O(l+1)^{n-l} \text{ for all $n < l$ and $ l \geq 1$}\]
\end{theorem}
\begin{proof} The proof follows from \cite{maa15} since the sheaf of differential operators $\Diff^l(\O(n))$ is the dual of the sheaf of principal parts.
\end{proof}

Hence $\Diff^l(\O(n))^{left} \neq \Diff^l(\O(n))^{right}$ as $\O_C$-module in general.

We denote the left and right $A$-module structure on $\P^l(E)$ as $\P^l(E)^{left}$ and $\P^l(E)^{right}$. The class
\[   \gamma^l(E):=[\P^l(E)^{left}]- [\P^l(E)^{right}] \in \K_0(A) \]
is zero in most cases. This is because $\P^l(E)$ is an extension $\P^{l-1}(E)$ with  $\sym^l(\Omega^1)$ and $\sym^l(\Omega^1)^{left}\cong \sym^l(\Omega^1)^{right}$ as $A$-modules for all $l\geq 1$.

\begin{lemma} \label{fund-diff}Let $k\rightarrow A$ be a commutative ring that is a finitely generated and regular $k$-algebra where $k$ is a field of characteristic zero and let $E$ be a finitely generated and projective $A$-module. There is for every $l\geq 1$
an exact sequence of left and right $A$-modules
\[ 0 \rightarrow \Diff^{l-1}(E)\rightarrow \Diff^l(E) \rightarrow \sym^l(\Der_k(A))\otimes \End_A(E) \rightarrow 0.\]
We let $\Diff^0(E):=\End_A(E)$.
\end{lemma}
\begin{proof} There is an exact sequence of left and right $A$-modules 
\[ 0\rightarrow \sym^l(\Omega^1_{A/K})\otimes_A E \rightarrow \P^l(E) \rightarrow \P^{l-1}(E) \rightarrow 0 \]
where $\P^l(E)$ is the $l$'th module of  principal parts. Since $A$ is regular it follows $\P^l(E)$ is a projective $A$-module of finite rank as left and right $A$-module. When we apply the functor $\Hom_A(-,E)$
to the sequence we get the claimed sequence and the Lemma follows.
\end{proof}

\begin{lemma} \label{class}Let $A$ be a commutative ring satisfying the hypothesis from Lemma \ref{fund-diff}.  Let $\Diff^l(E)^{left}$ and $\Diff^l(E)^{right}$ denote the left and right $A$-module structure on
$\Diff^l(E)$. The following holds in the Grothendieck group $\K_0(A)$ of $A$:
\[ [\Diff^l(E)^{left}]=[\End_A(E)^{left}](1+[\Der_k(A)^{left}]+[\sym^2(\Der_k(A)^{left}]+\cdots + \]
\[ [\sym^l(\Der_k(A)^{left}]).\]
A similar formula holds when we consider the right $A$-module structure.
\end{lemma}
\begin{proof} The Lemma follows from Lemma \ref{fund-diff} and an induction on $l$.
\end{proof}

\begin{theorem} Let $A$ be a commutative ring satisfying the hypothesis in Lemma \ref{fund-diff}. Let $E$ be a finitely generated and projective $A$-module and let 
\[ \eta^l(E):=[\Diff^l(E)^{left}]-[\Diff^l(E)^{right}]\in \K_0(A). \] 
It follows $\eta^l(E)=0$.
\end{theorem}
\begin{proof} The Theorem follows from Lemma \ref{class} since there is for every $l\geq 1$ an isomorphism $\sym^l(\Der_k(A))^{left}\cong \sym^l(\Der_k(A))^{right}$ of $A$-modules. 
\end{proof}

Hence the Grothendieck group $\K_0(A)$ does not detect that $\P^l(E)^{left} \neq \P^l(E)^{right}$ and $\Diff^l(E)^{left} \neq \Diff^l(E)^{right}$ in general (see \cite{maa14} for a more detailed discussion). 
The aim of this study is to construct "generalized jet bundles"  where the classes $\gamma^l(E)$ and $\eta^l(E)$ are non trivial in $\K_0(A)$ and to apply this in the study of Chern classes and Hodge theory. One want to construct non-trivial classes in $\H^{2*}(L,A)$ coming from $\K_0(L)$.

\begin{lemma} \label{curvrho} Let $k \rightarrow A$ be a map of unital commutative rings and let $f\in \Z^2(\Der_k(A),A)$ be a 2-cocycle. Let $(L,\alpha)$ be an $(A/k)$-Lie-Rinehart algebra and let $\fal \in \Z^2(L,A)$ be the pull back of $f$.
Assume $\rho: L(\fal) \rightarrow \End_k(E)$ is an $A\otimes_k A$-linear map, with $E$ a left $A$-module and let $i:L\rightarrow L(\fal)$ be the canonical injective map. Let $u=az+x, v=bz+y \in L(\fal)$ and define $R_\rho(u,v):=[\rho(u),\rho(v)]- \rho([u,v]) $. The following holds:
\[ R_\rho(u,v)(e)=R_{\rho \circ i}(x,y)(e) -\fal(u,v)e \]
where 
\[ R_{\rho \circ i}(x,y):= [\rho(i(x)),\rho(i(y))]- \rho(i([x,y])) .\]
\end{lemma}
\begin{proof} Let $\rho: L(\fal) \rightarrow \End_k(E)$ be an $A\otimes_k A$-linear map. We get
\[ R_{\rho}(u,v)(e)=[\rho(u),\rho(v)](e)-\rho([u,v])(e) .\]
We get
\[ \rho([u,v])(e)=\rho((\alpha(x)(b)-\alpha(y)(a)+f(\alpha(x),\alpha(y)))z +[x,y])(e)=\]
\[(\alpha(x)(b)-\alpha(y)(a)+f(\alpha(x),\alpha(y)))e+\rho(i([x,y]))(e).\]
Similarly we get
\[ [\rho(u),\rho(v)](e)=[\rho(i(x)),\rho(i(y))](e).\]
We get
\[ R_{\rho}(u,v)(e)=[\rho(i(x)),\rho(i(y))](e)-\rho([u,v])(e)-f(\alpha(x),\alpha(y))e= \]
\[ R_{\rho \circ i}(x,y)(e)-f(\alpha(x),\alpha(y))e=\]
\[  R_{\rho \circ i}(x,y)(e)- \fal(u,v)e .\]
The Lemma is proved.
\end{proof}

\begin{example} \label{families} Families of connections. \end{example}

Given two 2-cocycles $f,g\in \Z^2(\Der_k(A),A)$. It follows from Lemma \ref{corr} we get for any $\psi$-connection $\nabla:L\rightarrow \End_k(A)$
two connections
\[ \rho_{\fal}:L(\fal)\rightarrow \End_k(E)\]
and 
\[ \rho_{g^{\alpha}}:L(g^{\alpha})\rightarrow \End_k(E).\]
If $g=f+d^1\phi$ for an element $\phi \in \operatorname{C}^1(L,A)$, there is an isomorphism of Lie-Rinehart algebras $L(\fal)\cong L(g^{\alpha})$.
Hence the construction gives for a fixed $A$-module $E$, a family of connections $\rho_{\fal}:L(\fal)\rightarrow \End_k(E)$ parametrized by the set of 2-cocycles $\Z^1(\Der_k(A),A)$ 
and from Lemma \ref{curvrho} we see the curvature $R_{\rho_{\fal}}$ varies with the 2-cocycle $f$. One may ask if it is possible to use the family $(E,\rho_{\fal})$ to 
study the original connection $\nabla$ and its characteristic classes. If $g=f+d^1\phi$ it follows there is a canonical isomorphism
\[ L(\fal)\cong L(g^{\alpha}) \]
of Lie-Rinehart algebras, hence we may for any cohomology class $c:=\overline{\fal}\in \H^2(L,A)$ define $L(c):=L(\fal)$. We get for any $\psi$-connection $\nabla$ a family of connections
\[ \rho_{c}:L(c) \rightarrow \End_k(E) \]
parametrized by the cohomology group $\H^2(L,A)$. If $B:=\sym_k^*(\H^2(L,A)^*)$ We get in a canonical way a connection
\[  q^*(\rho_{c}):q^*L(c)\rightarrow \End_B(B\otimes_k E) \]
On $B\otimes_k A$. The element $c$ gives in a canonical way a $k$-rational point $x(c)\in \Spec(B)(k)$. Hence we may view $\rho_{c}$ as the restriction of $q^*(\rho_{c})$ to the fiber
$p^{-1}(x(c))$ where $x(c)\in \Spec(B)$. Here the maps $p,q$ are the canonical maps $q:A\rightarrow B\otimes_k A$ and $p:B\rightarrow B\otimes_k A$. We may ask if there is
a globally defined $\operatorname{D}(B\otimes_k A, g)$-Lie algebra 
\[ \beta: \overline{L}\rightarrow \operatorname{D}^1(B\otimes_k A,g) \]
for some 2-cocycle $g\in \Z^2(\Der_k(B\otimes_k A), B\otimes_k A)$ with $\overline{L}_{p^{-1}(x(c))}=L(c)$.

\begin{example} The curvature of a $\psi$-connection with $\psi=Id_E$. \end{example}

The map $R_{\rho \circ i}(x,y)$ from Lemma \ref{curvrho} is not in $\End_A(E)$ in general. We always have $R_{\rho \circ i}(x,y)\in \End_k(E)$.
In general it follows we get a map
\[ \nabla :=\rho \circ i: L\rightarrow \End_k(E) \]
and $R_{\nabla}(x,y)=R_{\rho \circ i}(x,y) \in \End_A(E)$ if $\rho(z)=\psi =Id_E$. In this case $\nabla$ is an ordinary connection.

\section{Functorial properties of universal rings of $\Dl$-Lie algebras}

In general if $L$ is a $k$-Lie algebra over a commutative ring $k$, we may construct the universal enveloping algebra $\U(L)$ of $L$. Let $\underline{\operatorname{Alg}}(k)$  be the category of 
associative unital $k$-algebras $R$ with $k$ contained in the center $Z(R)$ of $R$, and maps of unital associative $k$-algebras.
There is a functor
\[  F: \underline{\operatorname{Alg}}(k) \rightarrow \SETS \]
defined by
\[ F(R):=\{\text{the set of maps of $k$-Lie algebras $\rho: L \rightarrow R$} \} .\]
There is a canonical map of $k$-Lie algebras $i: L\rightarrow \U(L)$ and an isomorphism of functors
\begin{align}
&\label{isofunct} \Hom_{k-alg}(\U(L), R) \cong F(R) ,
\end{align}
hence the pair $(\U(L), i)$ represents the functor  $F(-)$. It follows the pair $(\U(L), i)$ is uniquely determined by this property. 
It follows there is an equality of sets
\[ \Hom_{k-alg}(\U(L), \End_k(E)) \cong F(\End_k(E)) \]
for any $k$-module $E$. Hence the set of $L$-modules $(E,\rho)$ are in one to one correspondence with left $\U(L)$-modules. We get an equivalence of categories 
\begin{align}
&\label{equivmodL} \Mod(L) \cong \Mod(\U(L)) 
\end{align}
between the category of $L$-modules and the category of $\U(L)$-modules preserving projective and  injective objects. 
We may use the equivalence \ref{equivmodL} and $\Ext$ and $\Tor$-groups to define the cohomology and homology of any $L$-module $V$:
\[ \H^i(L,V):=\Ext^i_{\U(L)}(k,V) \text{ and }\H_i(L,V):=\Tor_i^{U(L)}(k,V).\]
It follows functorial properties of the groups $\H^i(L,V)$ and $\H_i(L,V)$ follow from well known functorial properties of $\Ext$ and $\Tor$-groups for associative rings as proved in Cartan and 
Eilenberg's classical book \cite{cartan}. The cohomology and homology groups $\H^i(L,V)$ and $\H_i(L,V)$ may in the case when $k$ is a field or $L$ a projective 
$k$-module be calculated using Hochschild cohomology and homology.

The aim of this section is to construct an isomorphism and equivalence similar to \ref{isofunct} and \ref{equivmodL} and to define $\Ext$ and $\Tor$-groups
for the universal enveloping algebra $\U(A,L,f)$ of an $A/k$-Lie-Rinehart algebra $(L,\alpha)$ and a 2-cocycle $f\in \Z^2(L,A)$. We also do a similar construction for any $\Dl$-Lie algebra $(\tL, \ta, \tp,[,],D)$ and construct the universal ring $\Uo(\tL)$ of $\tL$.  The associative unital ring $\Uo(\tL)$ is applied to the study of the category $\Mod(\tL, Id)$ of $\tL$-connections $\rho:\tL\rightarrow \End_k(E)$ with $\rho(D)=Id_E$. We construct a functor $\psi$ realizing  the category $\Mod(\tL, Id)$ as a module category of left modules over  $\Uo(\tL)$, and use $\psi$ to construct $\Ext$ and $\Tor$-groups of arbitrary connections in $\Mod(\tL, Id)$.
We do a similar construction for the category $\Conn(\tL, Id)$.
More precisely we construct for any $\Dl$-Lie algebra $\tL$ two associative rings $\Uo(\tL)$ and $\Ur(\tL)$ containing the base ring $k$ in the center, with the property that there are exact equivalences of categories
\begin{align}
 F_1:\Mod(\tL, Id) \cong \Mod(\Uo(\tL)) \\
F_2: \Conn(\tL, Id) \cong \Mod(\Ur(\tL))
\end{align}
such that $F_1$ and $F_2$ preserve injective and projective objects. The categories $\Mod(\tL, Id)$ and $\Conn(\tL, Id)$ are categories of non-flat connections and the rings $\Uo(\tL)$ and $\Ur(\tL)$ are
non-Noetherian in general. Hence the rings $\Uo(\tL)$ and $\Ur(\tL)$ may be viewed as universal enveloping algebras for non-flat connections. The rings $\Uo(\tL)$ and $\Ur(\tL)$ contain 2-sided ideals $I^{\otimes}$ and $I^{\rho}$ with the property that the quotient rings $\tUo(\tL)$ and $\tUr(\tL)$ are almost commutative rings. If $A$ is Noetherian and $\tL$ is a finitely generated left $A$-module it follows $\tUo(\tL)$ and $\tUr(\tL)$ are
Noetherian (see Theorem \ref{mainequiv} and Theorem \ref{mainnoetherian}). We use $\Uo(\tL)$ and $\Ur(\tL)$ to construct $\Ext$ and $\Tor$ groups for non-flat connections. 
This construction was previously done for flat connections in \cite{rinehart}.

Let in this section $k$ be an arbitrary commutative unital ring and let $A$ be an arbitrary commutative unital $k$-algebra.
Let $(\tL, \ta,\tp,[,],D)$ be a $\Dl$-Lie algebra. A connection on $\tL$ is by Definition \ref{maindef} an $A\otimes_k A$-linear map
\[ \rho: \tL \rightarrow \End_k(E) \]
where $E$ is a left $A$-module. 

\begin{definition}
Let $\T^j_k(\tL):=\tL\otimes_k \otimes \cdots \otimes_k \tL= \tL^{\otimes^j_k}$ be the tensor product of $\tL$ with itself $j$ times over the ring $k$. By definition $\tL^{\otimes_k ^0}:=k$.
Let $\T^*_k(\tL):=\oplus_{i \geq 0}\T^i_k(\tL)$ be the \emph{tensor algebra} of $\tL$ over $k$. Let $\T^*_k(\tL)^j:=\oplus_{i \geq j} \T^i_k(\tL)$ for an integer $j \geq 0$.
Let $\T^*_k(\tL)_j:=\oplus_{i=0}^j \T^i_k(\tL)$.
\end{definition}

Note: It follows $\T^*_k(\tL)^j \subseteq \T^*_k(\tL)$ is a 2-sided ideal for every integer $j \geq 0$. There is moreover a filtration
\[ \T^*_k(\tL)_0  \subseteq \T^*_k(\tL)_1 \cdots \subseteq \T^*_k(\tL)_j \subseteq \cdots \T^*_k(\tL)  \]
that is compatible with the multiplication on $\T^*_k(\tL) $. This means for any elements $u\in \T^*_k(\tL)_i, v\in \T^*_k(\tL)_j$ it follows $u\otimes v\in \T^*_k(\tL)_{i+j}$.

Let $R$ be an associative unital $k$-algebra with $k$ in the center of $R$. 

\begin{definition} Let $E$ be an abelian group. A \emph{left $R$-module structure on $E$} is a map
\[ \sigma: R\times E \rightarrow E \]
where we write $\sigma(a,e):=ae$ for $a\in R, e\in E$. The action $\sigma$ should verify
$(a+b)e=ae+be, a(e+e')=ae+ae', (ab)e=a(be)$ and $1e=e$ for all $a,b\in R, e,e'\in E$ and $1\in R$ the multiplicative unity.
\end{definition}

\begin{lemma} \label{leftmodule} Let $E$ be an abelian group and let $R$ be an associative unital ring with $k$ in its center. 
There is a one-to-one correspondence between the set of all left $R$-module structures on $E$ and the set of all pairs $(\rho, \sigma_k)$ 
where $\sigma_k$ is a left $k$-module structure on $E$ and $\rho: R\rightarrow \End_k(E)$ is a map of associative $k$-algebras.
\end{lemma}
\begin{proof} If $\sigma: R\times E\rightarrow E$ is a left $R$-module structure on $E$ it follows $E$ is a left $k$-module and the map
$\rho: R\rightarrow \End_k(E)$ defined by $\rho(a)(e):=\sigma(a,e)$ is a map of $k$-algebras. Conversely, if $E$ is a $k$-module and $\rho: R\rightarrow \End_k(E)$ 
a map of $k$-algebras it follows in particular that $E$ is an abelian group and $\sigma: R\times E\rightarrow E$ defined by $\sigma(a,e):=\rho(a)(e)$ defines a left
$R$-module structure on $E$. This proves the Lemma.
\end{proof}

\begin{lemma} \label{tensoralgebra}Let $R$ be an associative $k$-algebra where $k$ is in the center of $R$. For each $k$-linear map $\rho: \tL\rightarrow R $ there is a canonical map
\[ \rho^i:\T^i_k(\tL)\rightarrow R \]
defined by
\[ \rho^i(u_1\otimes \cdots \otimes u_i):=\rho(u_1)\cdots \rho(u_i).\]
The abelian group $\T^*_k(\tL)$ is an associative $k$-algebra with $k$ in its center. There is a functorial equality of sets
\[ \Hom_k(\tL, R) \cong \Hom_{k-alg}(\T^*_k(\tL) ,R).\]
If $E$ is an abelian group it follows $E$ is a left $\T^*_k(\tL)$-module if and only if $E$ is a $k$-module and there is a map of $k$-algebras
\[ \phi: \T^*_k(\tL) \rightarrow \End_k(E).\]
\end{lemma}
\begin{proof} The proof is straight forward and is left to the reader.
\end{proof}

Since there is a functorial equality of sets
\[  \Hom_k(\tL, R) \cong \Hom_{k-alg}(\T^*_k(\tL) ,R)\]
it follows an abelian group $E$ is a left $\T^*_k(\tL)$-module if and only if $E$ is a left $k$-module and there is a map of $k$-modules $\phi': \tL \rightarrow \End_k(E)$.

Let $(\tL, \ta, \tp,[,], D)$ be a $\Dl$-Lie algebra and let $(E,\rho)$ be an $\tL$-connection. It follows $\rho:\tL \rightarrow \End_k(E)$ is a map of $A\otimes_k A$-modules.

Recall that  $\Mod(\tL)$ (resp. $\Mod(L(\fal)$) denote the categories of $\tL$-connections and morphisms (resp. $L(\fal)$-connections and morphisms). Recall that $\Mod(\tL, Id)$ denotes the category of $\tL$-connections
$(E,\rho)$ with $\rho(\iota)=Id_E$. Let for an associative ring $R$, $\Mod(R)$ denote the 
category of left $R$-modules and maps of $R$-modules.

\begin{definition} Let $(\tL, \ta, \tp,[,],\iota)$ be a $\Dl$-Lie algebra. A sequence of maps of $\tL$-connections
\[  \cdots \rightarrow (E_{i+1},\rho_{i+1}) \rightarrow^{\phi_{i+1}} (E_i,\rho_i) \rightarrow^{\phi_i} (E_{i-1}, \rho_{i-1}) \rightarrow \cdots \]
is \emph{exact at $(E_i,\rho_i)$} if and only if $Im(\phi_{i+1}=ker(\phi_i)$.
\end{definition}

Note: In \cite{maa1}, Definition A.7 we defined the \emph{universal algebra} $\U^{ua}(L)$ of a Lie-Rinehart algebra $(L,\alpha)$  where $L(\fal)$
was the abelian extension of $L$ with the 2-cocycle $\fal \in \Z^2(L,A)$. Definition A.7 in \cite{maa1} is not correct and in the following give a correct construction of the universal ring for any $\Dl$-Lie algebra.

\begin{definition}\label{ideals}  Let $(\tL, \ta, \tp,[,],D)$ be a $\Dl$-Lie algebra. Define the following 2-sided ideals in $\T^*_k(\tL)$:
\[ J_1:=\{ D-1\text{ and }au-(aD)\otimes u\text{ and }vb-v\otimes (bD):\text{ such that  $a,b\in A$ and $u,v\in \tL$.}   \} \]
and
\[ J_2:=\{ D-1\text{ and }au -(aD)\otimes u\text{ and } v\otimes(bD)-(bD)\otimes v-\tp(v)(b)D:  \] 
\[ \text{ such that $a,b\in A$ and $u,v\in \tL$} \}  .\]
\end{definition}

\begin{definition} \label{univring} Let $\Uo(\tL):=\T^*_k(\tL)/J_1$ and $\Ur(\tL):=\T^*_k(\tL)/J_2$. When we speak of the \emph{universal ring} of the $\Dl$-Lie algebra $\tL$ we refer to
$\Uo(\tL)$ or $\Ur(\tL)$. Let $p^{\otimes}: \tL \rightarrow \Uo(\tL)$ and $p^{\rho}:\tL \rightarrow \Ur(\tL)$ be the canonical maps. Let $J\subseteq \Uo(\tL)$ be a 2-sided ideal and define $\Uo_J(\tL):=\Uo(\tL)/J$.
There is a canonical map $p_J^{\otimes}: \tL \rightarrow \Uo_J(\tL)$. Define similarly $p_J^{\rho}: \tL \rightarrow \Ur_J(\tL)$ for any 2-sided ideal $J\subseteq \Ur(\tL)$.
\end{definition}

\begin{lemma} \label{assunital} Let $(\tL, \ta, \tp,[,],D)$ be a $\D$-Lie algebra.
The abelian groups $\Uo(\tL), \Uo_J(\tL), \Ur(\tL)$ and $\Ur_J(\tL)$ are associative unital rings with the element $\overline{D}:=1$ as unit.
\end{lemma}
\begin{proof} The associative ring $\T^*_k(\tL)$ has multiplicative unit $1\in k$. Passing to the quotient $\Uo(\tL)$ it follows $\overline{D}$ is the multiplicative unit.
The rest follows similarly and the Lemma follows.
\end{proof}

\begin{lemma} Let $(\tL,\ta, \tp,[,],D)$ be a $\Dl$-Lie algebra. There is a canonical $A\otimes_k A$-module structure on $\Uo(\tL)$ defined as follows: Let $a\in A, z\in \Uo(\tL)$. define
$az:=p^{\otimes}(aD)z, za:=zp^{\otimes}(aD)$. There is a similar structure on $\Uo_J(\tL), \Ur(\tL)$ and $\Ur_J(\tL)$. 
The maps $p^{\otimes}p_J^{\otimes}, p^{\rho}$ and $p_J^{\rho}$ are $A\otimes_k A$-linear maps mapping $D$ to the multiplicative identity. 
\end{lemma} 
\begin{proof} Let $a,b\in A, z,w\in \Uo(\tL)$. We get
\[ (a+b)z:=p((a+b)D)z=p(aD)z+p(bD)z=az+bz.\] 
Similarly it follows $1z=z1=z, a(z+w)=az+aw, (z+w)a=za+wa.$
We get
\[ (ab)z=p((ab)D)z=\overline{(ab)D}z=\overline{(aD)\otimes (bD)}z=p(aD)p(bD)z =a(bz).\]
Similarly $z(ab)=(za)b$ and it follows $\Uo(\tL)$ is a left $A\otimes_k A$-module with $p^{\otimes}$ an $A\otimes_k A$-linear map. The rest follows similarly and the Lemma
follows.
\end{proof}

\begin{definition} Let $\ALG$ be the category with objects associative unital $k$-algebras $R$ with $k$ in the center of $R$, where $R$ is equipped with the structure of a left $A\otimes_k A$-module.
A map of objects $R,R'\in \ALG$ is a map $f:R\rightarrow R'$ of unital $k$-algebras such that $f$ is also $A\otimes_k A$-linear.
Given a $\Dl$-Lie algebra $(\tL, \ta, \tp,[,],D)$. And let $\rho: \tL\rightarrow R$ be a $k$-linear map. Let $\T(\rho):\T^*_k(\tL) \rightarrow R$ be the induced map of $k$-algebras.
Let $J\subseteq \Uo(\tL)$ be a 2-sided ideal and let $J_{\T}\subseteq \T^*_k(\tL)$ be the inverse image in $\T^*_k(\tL)$.
Define the functor 
\[ \Conn_{\tL,J}: \ALG \rightarrow \underline{Sets} \]
by defining $\Conn_{\tL,J}(R)$ to be the set of $A\otimes_k A$-linear maps $\rho: \tL \rightarrow R$ such that $\rho(D)=1$ is the multiplicative identity in $R$, and $\T(\rho)(J_{\T})=0$. 
\end{definition}

Note: It is clear $\Conn_{\tL,J}$ is a covariant functor. By definition it follows the set
\[ \Conn_{\tL, J}(\End_k(E)) \]
is the set of $J$-flat connections
\[ \rho: \tL \rightarrow \End_k(E) .\]

\begin{lemma}\label{represent}  Let $(\tL, \ta,\tp,[,],D)$ be a $\Dl$-Lie algebra and let $J\subseteq \Uo(\tL)$ be a 2-sided ideal.
Let $R$ be an associative $k$ algebra with $k$ in the centre of $R$. There is a one-to-one correspondence between the set of maps of $k$-algebras $\tilde{\rho}:\Uo_J(\tL) \rightarrow R$ 
and the set of pairs $(\rho, \sigma)$ where $\sigma: A\otimes_k A \times R \rightarrow R$ is a left $A\otimes_k A$-module structure on $R$ and $\rho: \tL \rightarrow R$ is an $A\otimes_k A$-linear map
with $\rho(D)=1_R$ and $\T(\rho)(J_{\T})=0$.
\end{lemma}
\begin{proof} Let $a\in A, z\in R$ and let $\tilde{\rho}:\Uo_J(\tL)\rightarrow R$ be a map of unital $k$-algebras. Let $p:=\tilde{p}\circ p^{\otimes}_J:\tL \rightarrow R$.
Define $az:=p(aD)z$ and $za:=zp(aD)$. It follows $R$ is a left $A\otimes_k A$-module and the map $p$ is by definition $A\otimes_k A$-linear with $p(D)=1_R$ and $\T(p)(J_{\T})=0$. 
Let $\sigma$ be the $A\otimes_k A$-module structure defined and map $\tilde{\rho}$ to the pair $(p, \sigma)$. This correspondence is one-to-one and the Lemma is proved.
\end{proof}

\begin{proposition} \label{isofunctors} The notation is as in Lemma \ref{represent}. \label{representconn} There is an isomorphism of functors 
\[ \Conn_{\tL,J}(-) \cong \Hom_{k-alg}(\Uo_J(\tL), -), \]
hence the pair $(\Uo_J(\tL), p^{\otimes}_J)$ represents the functor $\Conn_{\tL, J}$.
It follows the pair $(\Uo_J(\tL), p^{\otimes}_J)$ is unique up to unique isomorphism.
\end{proposition}
\begin{proof} Assume $\rho: \tL \rightarrow R$ is an $A\otimes_k A$-linear map with $\rho(D)=1_R$ and $\T(\rho)(J_{\T})=0$. It follows $\rho$ induce a map of $k$-algebras
 $\tilde{\rho}: \Uo_J(\tL) \rightarrow R$ of associative $k$-algebras. By Lemma \ref{represent} this gives an isomorphism of functors
and the Proposition is proved.
\end{proof}

\begin{definition} let $(\tL,\ta, \tp,[,],D)$ be a $\Dl$-Lie algebra and let $J\subseteq \Uo(\tL)$ be a 2-sided ideal. The functor $\Conn_{\tL, J}$ is the 
\emph{moduli functor for $\tL$-connections with $J$-curvature zero}. The functor $\Conn_{\tL}:=\Conn_{\tL, (0)}$ where $(0)$ is the zero ideal, is the \emph{universal moduli functor for $\tL$-connections}.
\end{definition}

Note: It makes sense to call $\Conn_{\tL}$ universal since it parametrize all connections $\rho: \tL \rightarrow \End_k(E)$ with no condition on the $J$-curvature $R^J_{\rho}$.

\begin{example} \label{moduliconnection} Moduli spaces of connections.\end{example}

Proposition \ref{isofunctors} gives an alternative approach to the study of moduli spaces of connections using modules on associative rings. 
See \cite{simpson} for an approach using the quot scheme and GIT theory. In \cite{simpson} the author does the following: If $X$ is a smooth projective complex variety and $\Dl_{\lambda}$
a generalized sheaf of rings of differential operators, Simpson parametrize left $\Dl_{\lambda}$-modules $M$ where the underlying coherent $\O_X$-module $M$ is p-semistable and has a given Hilbert polynomial $P(x)$. Simpson uses this construction to study moduli spaces of representations of the toppological fundamental group of $X$ via the Riemann-Hilbert correspondence. Simpson's ring $\Dl_{\lambda}$ is a global version of the generalized enveloping algebra  $\U(A,L,f)$ studied in the  paper \cite{maa1}. Since the ring $\U(A,L,f)$ is a quotient of $\Uo(L(\alpha^*f))$ we may ask if it is possible to do a similar construction with $\Uo(\tL)$. We 
may globalize $\Uo(\tL)$ to get a sheaf of associative unital rings on a projective scheme $Y$ and ask if it is possible to construct parameter spaces $\M(\Uo_J(\tL), P(x))$ for left $\Uo_J(\tL)$ modules $E$, 
where $E$ is a coherent left  $\O_Y$-module with fixed Hilbert polynomal $P(x)$.  If such a construction is possible, the space $\M(\Uo(\tL), P(x))$ would in some sense be a \emph{universal moduli space} for connections, since it contains all left $\Uo(\tL)$-modules $E$ that are coherent as $\O_Y$-modules with a fixed Hilbert polynomial $P(x)$. There are a lot of technical details that has to be checked: One has to define $\Dl$-Lie algebras in a relative setting for morphisms of schemes. The construction of $\Uo(\tL)$ is functorial in $\tL$, hence this may give a construction of general moduli spaces for connections. 

 If $Y$ is a complex projective manifold it follows in particular that $Y$ is a smooth projective algebraic variety. Any holomorphic finite rank complex vector bundle $E$ is algebraic. A flat connection 
$\nabla: E\rightarrow \Omega^1_{Y}\otimes E$ corresponds to a left $\Dl_Y$-module $E$ where $\Dl_Y$ is the sheaf of rings of polynomial differential operators on $Y$. Associated to $(E,\nabla)$ we may
construct the characteristic variety $SS(E,\nabla)$ and we use $SS(E,\nabla)$ to define holonomicity. Hence the category of holonomic $\Dl_Y$-modules is a sub category of the category of flat connections.
The Riemann-Hilbert correspondence in its most elementary form says that there is an equivalence of categories between the category of holonomic $\Dl_Y$-modules of finite rank as $\O_Y$-module and the
category of finite dimensional complex representations of the topological fundamental group $\pi_1(Y)$ of $Y$. Hence one way to construct explicit non-trivial examples of flat algebraic connections, is to construct
a non-trivial finite dimensional complex representation $\rho$ of $\pi_1(Y)$ and to use the Riemann-Hilbert correspondence to pass from $\rho$ to a flat algebraic connection $(E_{\rho}, \nabla_{\rho})$ on $Y$.
It is easier to check if a representation of $\pi_1(Y)$ is non-trivial than to check if the corresponding flat connection is non-trivial.

\begin{example}\label{atiyah} Atiyah classes and Atiyah sequences for $\Dl$-Lie algebras. \end{example}

For any $\Dl$-Lie algebra $\tL$ and any left $A$-module $E$ there is an exact sequence of $A\otimes_k A $-modules (a generalized Atiyah sequence)
\[ 0 \rightarrow E \rightarrow J^1_{\tL}(E) \rightarrow \tL\otimes_A E \rightarrow 0 \]
which  is right split by an $A\otimes_k A$-linear map $s:\tL\otimes_A E\rightarrow J^1_{\tL}(E)$ if and only if $E$ has an $(\tL,\psi)$-connection $\rho$ with $(\rho,E)\in \Conn(\tL, \End)$. 

\begin{lemma} Define $J^1_{\tL}(E):=E\oplus \tL\otimes_A E$ with the obvious left $A$-module structure and the following right $A$-module structure: Let $(x, u\otimes y)\in J^1_{\tL}(E), a\in A$ and define
$(x,u\otimes y)a:=(xa+\tp(u)(a)y, u\otimes (ya))$. It follows $J^1_{\tL}(E)$ is an $A\otimes_k A$-module and the canonical sequence
\begin{align}
&\label{atiyahsequence} 0 \rightarrow E \rightarrow J^1_{\tL}(E) \rightarrow \tL \otimes_A E \rightarrow 0 
\end{align}
is an exact sequence of $A\otimes_k A$-modules, where we have given $E$ the trivial right $A$-module structure. 
The sequence \ref{atiyahsequence} is right split by an $A\otimes_k A$-linear map $s$  if and only if there is an $(\tL, \psi)$-connection $\rho: \tL \rightarrow \End_k(E)$
with $\rho$ a left $A$-linear map and $\rho(u)(ae)=a\rho(u)(e)+\tp(u)(a)\psi(e)$ for some $\psi\in \End_A(E)$. Hence the pair $(\rho,E)$ is an object in $\Conn(\tL, \End)$.
\end{lemma}
\begin{proof} The proof is an exercise.
\end{proof}

\begin{definition} Let $(\tL, \ta, \tp,[,],D)$ be a $\Dl$-Lie algebra and let $E$ be a left $A$-module. The class $a_{\tL}(E) \in \Ext^1_{A\otimes_k A}(\tL\otimes_A E,E)$ is the \emph{$\Dl$-Atiyah class} of $E$
The sequence \ref{atiyahsequence} is the \emph{$\Dl$-Atiyah sequence} of $E$.  The $A\otimes_k A$-module $J^1_{\tL}(E)$ is the \emph{first order $\tL$-jet bundle} of $E$.
\end{definition}

The $\Dl$-Atiyah class  $a_{\tL}(E)\in \Ext^1_{A\otimes_k A}(\tL\otimes_A E,E)$ has the propety that  $a_{\tL}(E)=0$ if and only if $E$ has an $(\tL, \psi)$-connection $\rho$ and this construction globalize. Hence for any coherent $\O_Y$-module  $\mathcal {E}$ with Hilbert polynomial $P(x)$ there is a class
$a_{\tL}(\mathcal{E})$ which measures when $\mathcal{E}$ has an $(\tL,\psi)$-connection. Hence if the set $\M((\tL, \psi), P(x))$ parametrize the set of coherent $\O_Y$-modules
$E$ with Hilbert polynomial $P(x)$ that has an $(\tL, \psi)$-connection, we may describe $\M((\tL,\psi), P(x)) \subseteq \operatorname{Quot}(Y,P(x))$ as a subset using the $\Dl$-Atiyah class $a_{\tL}(E)$.

Note: If $Y$ is smooth of finite type over a field $k$ of characteristic zero and $E$ is a coherent $\O_Y$-module with a connection $\nabla: E \rightarrow \Omega^1_Y \otimes E$, it follows $E$ is locally free.
This does not hold in general for a $\Dl$-Lie algebra $\tL$ and a connection $\rho: \tL \rightarrow \End_k(E)$.

%Hence "pointwise" there is an "inclusion" of spaces $\M(\Uo_J(\tL),P(x)) \subseteq \M(\Uo(\tL),P(x))$ for any 2-sided ideal $J$.
%One has to give a rigorous definition of the functor $\Conn_{\tL}$ for projective schemes and to check if the moduli space $\M(\Uo(\tL), P(x))$ is a closed subscheme of the corresponding Hilbert scheme  
%$\operatorname{Hilb}(Y, P(x))$. Moreover we must check if there are closed immersions of schemes $\M(\Uo_J(\tL), P(x)) \subseteq \M(\Uo(\tL), P(x))$.

Note: There are relationships betwen the cohomology of Hilbert schemes of points and infinite dimensional Lie algebras (see \cite{qin1}, \cite{qin2}).

\begin{example} The ring of differential operators.\end{example}

Let $k$ be a field of characteristic zero and let $A$ be a regular $k$-algebra of finite type. It follows $\Diff(A)\cong \U(A,L)$ where $L:=\Der_k(A)$. Let $f:=0$ be the zero 2-cocycle for $\Der_k(A)$
with values in $A$. Let $\tL$ be the $\Dl$-Lie algebra associated to the pair $(L,f)$. We may define $\Diff(A)$ as follows.
\[  \Diff(A):= \T^*_k(\tL)/J_1 \]
where $J_1$ is the 2-sided ideal generated by $D-1, au-(aD)\otimes u, ua-u\otimes (aD), u\otimes v-v\otimes u -[u,v]$ for $a\in A, u,v\in \tL$. Hence $\Diff(A)$ equals $\Uo(\tL)/J$ where
$J$ is the 2-sided ideal generated by $u\otimes v-v\otimes u-[u,v]$ for $u,v\in \Uo(\tL)$. There is a canonical map $p:\tL \rightarrow \Diff(A)$ which is $A\otimes_k A$-linear. Define the functor
\[ \Conn_{\tL, flat}: \ALG \rightarrow \SETS \]
by  letting $\Conn_{\tL,flat}(R)$ to be the set of $A\otimes_k A$-linear maps $\rho: \tL \rightarrow R$ with $\rho(D)=1_R$ and $\rho$ a map of $k$-Lie algebras. It follows there is an equality of sets
\[ \Conn_{\tL, flat}(R)=\Hom_{k-alg}(\Diff(A), R) \]
hence the pair $(\Diff(A),p)$ represents the functor $\Conn_{\tL, flat}$. Hence we may view the ring of differential operators as the solution to a universal problem. The functor $\Conn_{\tL, flat}$ is in some sense
a quotient functor of $\Conn_{\tL}$. When we view the associative rings $\Diff(A), \Uo(\tL)$ and $\Ur(\tL)$ as solutions to universal problems, we can prove that the constructions localize well using the uniqueness
of the representing object. If we have two pairs $(U_i, p_i), i=1,2$ representing the same functor $\Conn_{\tL, J}$, it follows there is a canonical isomorphism $(P_1,p_1)\cong (P_2,p_2)$. It may be easier
to prove that $(P_i,p_i)$ represent the same functor than to write down an explicit isomorphism $(P_1,p_1)\cong (P_2,p_2)$.

\begin{example} The ring $\Diff(E)$ is not almost commutative in general.\end{example}

Let $A$ be a commutative ring and let $E:=A^2$ be the free rank 2 $A$-module. It follows $\Diff(E)$ is the matrix ring of rank two matrices on $\Diff(A)$. It follows $\Diff^1(E)$ is the $A\otimes_k A$-module
of rank 2 matrices on $\Diff^1(A)=A\oplus \Der_k(A)$. One checks that it is not true in general that for $x,y\in \Diff^1(E)$ it follows $[x,y]\in \Diff^1(E)$, hence $\Diff(E)$ is not almost commutative in general.
The same holds for $\Diff(R)$ where $R$ is a non-commutative associative unital $k$-algebra. Since $\Diff^0(R):=R$ is follows $\Diff(R)$ is not almost commutative.

The canonical map $\rho:\tL \rightarrow \End_k(\Uo(\tL))$ is an $A\otimes_k A$-linear map with $\rho(D)=1$ hence $\rho$ is a connection. We get an induced map of $k$-algebras
\[ \Uo(\rho): \Uo(\tL) \rightarrow \End_k(\Uo(\tL)) .\]
There is the canonical filtration $\Uo(\tL)^i \subseteq \Uo(\tL)$ and we get for every $i \geq 0$ an induced map
\[ \Uo(\rho)^i: \Uo(\tL)^i \rightarrow \Diff^i(\Uo(\tL)) .\]
Hence the elements in $\Uo(\tL)$ act on $\Uo(\tL)$ as differential operators. We get an inclusion of rings 
\[ s:\Uo(\tL) \subseteq \Diff(\Uo(\tL)) .\]
Since $\Diff(\Uo(\tL))$ is not almost commutative in general, we cannot use $s$ to conclude that $\Uo(\tL)$ is an almost commutative ring.

\begin{lemma} The ring $\Uo(\tL)$ is almost commutative.
\end{lemma} 
\begin{proof} There is an embedding $s:\Uo(\tL) \subseteq \Diff(\Uo(\tL))$ and $s(u)(z):=uz$ is multiplication with the element $u\in \Uo(\tL)$. There is a filtration on $\Uo(\tL)$ defined as follows:
There is the canonical map $p:\T^*_k(\tL)\rightarrow \Uo(\tL)$. Define $\Uo(\tL)^0:=\{a\overline{D}:a\in A\}$. Define $\Uo(\tL)^i:=p(\T^*_k(\tL)_i)$ for $i \geq 1.$ It follows the filtration defined is compatible with 
the multiplication on $\Uo(\tL)$. If $aD\in \Uo(\tL)^0, u\in \Uo(\tL)^1$ we get an element $aDu=au\in \Uo(\tL)^1$. Let $b\in A$. We get
\[ [s(u), s(b)]= h(bD)-(bD)u=ub-bu=\tp(u)(b)D\in \Uo(\tL)^0.\]
Assume $u\in \Uo(\tL)^i, v\in \Uo(\tL)^j$ with $i+j=k$ and the hypothesis holds for $i+j=k-1$. Since $\Uo(\tL)\subseteq \Diff(\Uo(\tL))$ we may argue as follows: Let $a\in A$.
We may write
\[ [s(u), \phi_a]:=ua-au=\phi_1\]
and
\[ [s(v)m\phi_a]=va-av=\phi_2 \]
with $\phi_1\in \Diff^{i-1}(\Uo(\tL))$ and $\phi_2\in \Diff^{j-1}$
it follows $ua=au+\phi_1, va=av+\phi_2$. We get
\[ [s(u)s(v),\phi_a]=uva-auv=auv+\phi_1v+u\phi_2-auv=\phi_1v+u\phi_2.\]
We moreover get
\[ [s(v)s(u),\phi_a]=vua-avu=avu +\phi_2u+v\phi_1-avu=\phi_2u+v\phi_1.\]
We get
\[ [s(u)s(v)-s(v)s(u), \phi_a]= [\phi_1, v]+[u,\phi_2].\]
By induction it follows $[\phi_1,v]+[u,\phi_2]\in \Diff^{i+j-2}(\Uo(\tL))$ hence $[[s(u),s(v)],\phi_a]\in \Diff^{i+j-2}(\Uo(\tL))$ and hence $[s(u),s(v)]\in \Diff^{i+j-1}(\Uo(\tL))$.
It follows $s(\Uo(\tL))\subseteq \Diff(\Uo(\tL))$ is an almost commutative ring and since $s$ is an injective map it follows $\Uo(\tL)$ is almost commutative. The Lemma follows.
\end{proof}

\begin{example} Hochschild homology and cyclic homology of almost commutative PBW-algebras and finite dimensionality.\end{example}

In \cite{loday} the following is proved:
Let $U:=\{U_i\}$ be an almost commutative PBW-algebra containing $\mathbb{Q}$, and let $S:=\operatorname{Sym}_A^*(L)$ with $A:=U_0$ and $L:=U_1/U_0$.
It follows $S$ is a Poisson algebra with Poisson product $\{,\}$. We may define the complex $(\Omega^*_{S/k}, \delta)$ and mixed complex $(\Omega^*_{S/k}, \delta,d)$ in the sense of
\cite{loday}, Section 3.3.6. There is the following result:

\begin{theorem} \label{loday} There are isomorphisms
\begin{align}
&\label{h} \operatorname{HH}_*(U) \cong \H_*(\Omega^*_{S/k}, \delta) \\
&\label{ch} \operatorname{HC}_*(U) \cong \operatorname{HC}_*(\Omega^*_{S/k},\delta,d)
\end{align}
Here $\operatorname{HH}_*(U)$ is Hocschild homology of $U$ and $\operatorname{HC}_*(U)$ is cyclic homology of $U$.
\end{theorem}
\begin{proof} See \cite{loday} Theorem 3.3.9.
\end{proof}

Hence much is known on Hochschild and cyclic homology of $\Uo(\tL)$ in the case when $\Uo(\tL)$ is an almost commutative PBW-algebra containing the rational numbers. If the base field is of characteristic zero, 
we may construct  the groups $\Ext_{\Uo(\tL)}(V,W)$ and $\Tor_i^{\Uo(\tL)}(V,W)$ using Hochschild cohomology and homology:
\begin{align}
&\label{ext}\Ext^i_{\Uo(\tL)}(V,W)\cong \operatorname{HH}^*(\Uo(\tL), \Hom_k(W,V)) .\\
&\label{tor} \Tor_i^{\Uo(\tL)}(V,W) \cong \operatorname{HH}_*(\Uo(\tL), \Hom_k(V,W)) .
\end{align}
It may be we can prove results similar to Theorem \ref{loday} for the groups in \ref{ext} and \ref{tor} and to get results on finite dimensionality. If $A$ contains a field $k$ of characteristic zero and $T:=\Spec(A)$
is smooth and of finite type over $k$ it follows $S:=\operatorname{Sym}_A^*(L)$ is regular over $k$. Hence $U:=\Spec(S)$ is smooth over $k$. The complexes in Theorem \ref{loday}
are similar to the DeRham complex, and one should therefore expect the cohomology groups $\operatorname{HH}_*(U)$ and $\operatorname{HC}_*(U)$ to be finite dimensional in such cases.

A large class of cohomology and homology theories can be constructed using the Ext and Tor groups defined in this paper, and it may be we can use the methods sketched above to prove finite dimensionality of such cohomology and homology groups.

\begin{definition} Let $(\tL, \ta,\tp,[,],D)$ be a $\Dl$-Lie algebra and let $\rho: \tL \rightarrow \End_k(E)$ be a connection. Let $\Uo(\rho):\Uo(\tL)\rightarrow \End_k(E)$. Let $J\subseteq \Uo(\tL)$
be a 2-sided ideal. We say the \emph{$J$-curvature of $\rho$ is zero} if $\Uo(\rho)(J)=0$. We write $R^J_{\rho}=0$. An $\tL$-connection $(\rho,E)$ is \emph{$J$-flat} if $R^J_{\rho}=0$.
\end{definition}

\begin{example} Classical flat connections.\end{example}

Let $(\tL, \ta,\tp,[,],D)$ be a $\Dl$-Lie algebra and let $\rho: \tL \rightarrow \End_k(E)$ be a connection. Let $\Uo(\rho):\Uo(\tL)\rightarrow \End_k(E)$ be the induced map of $k$-algebras. 
Let $J\subseteq \Uo(\tL)$ be the 2-sided ideal generated by elements on the form $u\otimes v-v\otimes u -[u,v]$ for $u,v\in \tL$. It follows $R^J_{\rho}=0$ if and only if $\rho$ is a map
of $k$-Lie algebras. Hence $\rho$ is $J$-flat if and only if $\rho$ is flat in the classical sense.

\begin{theorem} \label{univfunctor} Let $\underline{Rings}$ denote the category of associative unital rings and morphisms. Definition \ref{univring} gives rise to two covariant functors
\begin{align}
&\label{ring1}\Uo: \ATiA \rightarrow \underline{Rings} \\
&\label{ring2}\Ur: \ATiA \rightarrow \underline{Rings}.
\end{align}
\end{theorem}
\begin{proof} Given a map of $\Dl$-Lie algebras $\phi: \tilde{L} \rightarrow \tilde{L}'$ define the following map
\[ \T(\phi): \T^*_k(\tilde{L})  \rightarrow \T^*_k(\tilde{L}') \]
by
\[ \T(\phi)(u_1\otimes \cdots \otimes u_i):=\phi(u_1)\otimes \cdots \otimes \phi(u_i)\in \tilde{L}'^{\otimes_k^i}.\]
We get in a canonical way a map of associative unital rings 
\[ \T(\phi): \T^*_k(\tL) \rightarrow \T^*_k(\tL').\]
It follows $\T(\phi)(D-1)=\T(\phi)(D)-\T(\phi)(1)=D'-1$,
\[  \T(\phi)(au-(aD)\otimes u)=a\phi(u)-a\phi(D)\otimes \phi(u)=a\phi(u)-(aD')\otimes \phi(u) \]
and
\[ \T(\phi)(ua-u\otimes (Da)=\phi(u)a -\phi(u)\otimes \phi(D)a=\phi(u)a-\phi(u)\otimes (D'a) \]
hence $\T(\phi)$ induce a canonical map of rings
\[\Uo(\phi):\Uo(\tilde{L})\rightarrow \Uo(\tilde{L}).\]
If $\phi, \psi$ are maps in $\ATiA$ it follows $\Uo(\phi \circ \psi)=\Uo(\phi) \circ \Uo(\psi)$ hence $\Uo(-)$ define by Lemma \ref{assunital} a covariant functor as claimed. A similar result holds
for $\Ur(-)$ and the Theorem follows.
\end{proof}

\begin{lemma} \label{lemmaideal} Let $(\tL, \ta, \tp,[,],D)$ be a $\Dl$-Lie algebra and let $E$ be a left $k$-module. Let $i: \tL \rightarrow \T^*_k(\tL)$ be the canonical injection map
and let $\tilde{\rho}:\T^*_k(\tL) \rightarrow \End_k(E)$ be a unital map of associative $k$-algebras. Let $\rho:=\tilde{\rho} \circ i$ and let $J_1$ be the 2-sided ideal from Definition \ref{ideals}.
It follows $\tilde{\rho}(J_1)=0$ if and only if $E$ is a left $A$-module, $\rho(D)=Id_E$ and $\rho$ is $A\otimes_k A$-linear map.
\end{lemma}
\begin{proof} Assume $\tilde{\rho}(J_1)=0$ and define $ae:=\rho(aD)(e)$ for $a\in A, e\in E$. Since $D-1\in J_1$ it follows $\rho(D)=Id_E$. Since $a(bD)-(aD)\otimes (bD)\in J_1$ it follows 
\[ (ab)e=\rho((ad)D)(e)=\rho(aD)\rho(bD)(e)=a(be).\]
We get
\[ (a+b)e=\rho((a+d)D)(e)=\rho(aD+bD)(e)=\rho(aD)(e)+\rho(bD)(e)=ae+be\]
and
\[ a(e+e')=\rho(aD)(e+e')=\rho(aD)(e)+\rho(aD)(e')=ae+ae'.\]
and
\[ 1e=\rho(1D)(e)=\rho(1)(e)=Id_E(e)=e\]
for all $a,b\in A, e,e'\in E$ and it follows $E$ is a left $A$-module. Since $au-(aD)\otimes u, ua- u\otimes (aD)\in J_1$ it follows 
\[ \rho(au)(e)=\rho(aD)\rho(u)(e)=a\rho(u)(e) \]
and
\[ \rho(ua)(e)=\rho(u)\rho(aD)(e)=\rho(u)(ae) \]
hence $\rho$ is $A\otimes_k A$-linear. Conversely one proves that if $E$ is a left $A$-module, $\rho(D)=Id_E$ and $\rho$ is an $A\otimes_k A$-linear map it follows $\tilde{\rho}(J_1)=0$ and the Lemma
follows.
\end{proof}

\begin{corollary}\label{modulestructure}Let $(\tL, \ta, \tp,[,],D)$ be a $\Dl$-Lie algebra and let $E$ be a left $k$-module. Let $J\subseteq \Uo(\tL)$ be a 2-sided ideal.
There is a one to one correspondence between left $\Uo_J(\tL)$-module structures on $E$ and
$\tL$-connections $\rho: \tL \rightarrow \End_k(E)$ with $\rho(D)=Id_E$ and $R^J_{\rho}=0$.
\end{corollary}
\begin{proof} By Lemma \ref{leftmodule} it follows a left $\Uo_J(\tL)$-module structure on $E$ corresponds to a map of associative $k$-algebras
\[ \tilde{\rho}: \Uo_J(\tL) \rightarrow \End_k(E).\]
Let $i: \tL \rightarrow \T^*_k(\tL)$ be the canonical injective map, let $p:\T^*_k(\tL) \rightarrow \Uo_J(\tL)$ be the canonical projection map and let $j:=p \circ i$.
It follows the induced map $\tilde{\rho}\circ p: \T^*_k(\tL) \rightarrow \End_k(E)$ is a map of $k$-algebras with $\tilde{\rho}\circ p(J_{\T})=0$. It follows the induced map $\rho:=\tilde{\rho}\circ j$
is an $A\otimes_k A$-linear map
\[ \rho: \tL \rightarrow \End_k(E) \]
with $\rho(D)=Id_E$ and $R^J_{\rho}=0$, hence $(E,\rho)$ is an element in $\Mod(\tL, Id)$ with $J$-curvature zero. 
The converse is similar hence we get a one to one correspondence and the Lemma follows.
\end{proof}

\begin{lemma}Let $(\tL, \ta, \tp,[,],D)$ be a $\Dl$-Lie algebra and let $\phi:(E,\rho_E)\rightarrow (F, \rho_F)$ ba a map of $\tL$-connections in $\Mod(\tL, Id)$.
It follows the map $\phi$ induce a map of $\Uo(\tL)$-modules $\phi: (E,\Uo(\rho_E))\rightarrow (F, \Uo(\rho_F))$. This means that for any element $x\in \Uo(\tL)$ and $e\in E$ it follows 
$ \phi(xe)=x\phi(e)$. Hence $\phi$ is an $\Uo(\tL)$-linear map.
\end{lemma}
\begin{proof} If $x \in \Uo(\tL)$ it follows $x=p(y)$ where $p: \T^*_k(\tL) \rightarrow \Uo(\tL)$ and $y:= \sum  x_{i_1}\otimes \cdots \otimes x_{i_l}$
with $x_{i_j}\in \tL$ for all $i_j$. By assumption it follows $\phi(\rho_E(x_u) e)=\rho_F(x_u)\phi(e)$ for all $x_u\in \tL$. We get
\[ \phi(xe)= \phi( \sum  x_{i_1}\otimes \cdots \otimes x_{i_l} e):=\]
\[ \phi( \sum  \rho_E(x_{i_1}) \circ \cdots \circ \rho_E(x_{i_l})(e)) = \sum  \rho_F(x_{i_1}) \circ \cdots \circ \rho_F(x_{i_l})(\phi(e)):= x\phi(e).\]
The Lemma follows.
\end{proof}

Given an object $(E,\rho)$ if $\Mod(\tL, Id)$ let 
\[  \Uo(\rho): \Uo(\tL)\rightarrow \End_k(E) \]
be the corresponding map of unital associative $k$-algebras constructed in Corollary \ref{modulestructure}. We get a map of objects
\[ F_{\tL}: \Mod(\tL, Id) \rightarrow \Mod(\Uo(\tL)) \]
defined by
\[F_{\tL}(E,\rho):=(E,\Uo(\rho)).\]
Any map $\phi:(E,\rho_E)\rightarrow (F,\rho_F)$ of $\tL$-connections induce a map of $\Uo(\tL)$-modules
\[ F_{\tL}(\phi): (E, \Uo(\rho_E)) \rightarrow (F,\Uo(\rho_F)) \]
defined by $F_{\tL}(\phi):=\phi$. For any pair of composable morphisms $\phi,\psi$ of connections it follows $F_{\tL}(\phi \circ \psi)=F_{\tL}(\phi) \circ F_{\tL}(\psi)$ hence $F_{\tL}$ is a functor.

Given a left $\Uo(\tL)$-module $(E,\tilde{\rho}_E)$ where $\tilde{\rho}_E:\Uo(\tL)\rightarrow \End_k(E)$ is a map of associative $k$-algebra we get an induced $\tL$-connection
\[ \rho_E: \tL \rightarrow \End_k(E) \]
with $\rho_E(D)=Id_E$. We may define the map of sets
\[ G_{\tL}:\Mod(\Uo(\tL) \rightarrow \Mod(\tL, Id) \]
by
\[ G_{\tL}(E,\tilde{\rho}_E):=(E,\rho_E).\]
Any map $\phi:(E,\tilde{\rho}_E)\rightarrow (F, \tilde{\rho}_F) $ of left $\Uo(\tL)$-modules induce a map of $\tL$-connections
\[ \phi:(E,\rho_E) \rightarrow (F,\rho_F)\]

\begin{definition} \label{modequivalence} Let $(\tL, \ta, \tp,[,],D)$ be a $\Dl$-Lie algebra.
We get by the discussion above a well defined covariant functor
\[ F_{\tL}: \Mod(\tL, Id) \rightarrow \Mod(\Uo(\tL)) \]
defined by
\[ F_{\tL}(\phi:(E,\rho_E)\rightarrow (F,\rho_F)):=\phi: (E,\Uo(\rho_E)) \rightarrow (F,\Uo(\rho_F)).\]
Similarly we get a covariant functor
\[ G_{\tL}:\Mod(\Uo(\tL)) \rightarrow \Mod(\tL, Id) \]
defined by
\[ G_{\tL}(\phi: (E,\tilde{\rho}_E) \rightarrow (F, \tilde{\rho}_F)):=\phi: (E,\rho_E)\rightarrow (F,\rho_F).\]
\end{definition}

Hence the functors $F_{\tL}$ and $G_{\tL}$ acts as the identity on morphisms  in $\Mod(\tL, Id)$ and $\Mod(\Uo(\tL))$. It follows $F_{\tL}$ and $G_{\tL}$ are exact functors.
Since $F_{\tL}\circ G_{\tL}(E, \tilde{\rho}_E)=(E,\tilde{\rho}_E)$ and $G_{\tL}\circ F_{\tL}(E,\rho)=(E,\rho)$,
it follows $F_{\tL}$ and $G_{\tL}$ are exact equivalences of categories between $\Mod(\tL, Id)$ and $\Mod(\Uo(\tL))$.
A similar result holds for the categories $\Conn(\tL, Id)$ and $\Mod(\Ur(\tL))$. There is an exact equivalence between them and the proof is similar to the above proof hence is left to the reader.

We sum up in a Theorem:

\begin{theorem} \label{mainDL} Let $(\tL,\ta, \tp,[,],D)$ be a $\Dl$-Lie algebra. The functor $F_{\tL}$ from Definition \ref{modequivalence} define an exact equivalence of categories
\[ F_{\tL}:\Mod(\tL, Id) \rightarrow \Mod(\Uo(\tL)) \]
with the following property: An object $(P, \rho_P)$ in $\Mod(\tL, Id)$ is a projective (resp. injective) object if and only if the corresponding object $F_{\tL}(P,\rho_P)$ is a projective (resp. injective) object 
in $\Mod(\Uo(\tL))$.
The pair $(\Uo(\tL), p^{\otimes})$ represents the functor 
\[ \Conn_{\tL}: \ALG \rightarrow \SETS \]
hence there is an isomorphism of sets
\[ \Hom_{k-alg}(\Uo(\tL), R) \cong \Conn_{\tL}(R) \]
for any object $R\in \ALG$. Hence the pair $(\Uo(\tL), p^{\otimes})$ is unique up to unique isomorphism. The ring $\Uo(\tL)$ has a canonical filtration $\Uo(\tL)^i$ for $i\geq 0$
and $\Uo(\tL)$ is an almost commutative ring.
\end{theorem}
\begin{proof} The first part of the Theorem follows from the above discussion. Assume $(P,\rho_P)$ is a projective object in $\Mod(\tL, Id)$. It follows that for any surjection of $\tL$-connections
\[ \phi: (E,\rho_E) \rightarrow (F, \rho_F) \]
and any map $f:(P, \rho_P) \rightarrow (F, \rho_F)$ there is a map of connections $\tilde{f}:(E,\rho_E) \rightarrow (F,\rho_F)$ with $\phi \circ \tilde{f}=f$.
Assume $F_{\tL}(P, \rho_P):=(P,\Uo(\rho_P))$ is projective object in $\Mod(\Uo(\tL))$. Let $\phi:E\rightarrow F$ be a surjection of $\Uo(\tL)$-modules and let $f:P\rightarrow F$ ba any map 
of $\Uo(\tL)$-modules. Applying the functor $G_{\tL}$ we get maps $\phi, f, \tilde{f}$ of connections in $\Mod(\tL,Id)$ with $\phi \circ \tilde{f}=f$, since the connection $(P,\rho_P)$ is a projective object in 
$\Mod(\tL, Id)$. It follows by definition of the functor $F_{\tL}$ there is an equality
$\phi \circ \tilde{f}=f$ as maps of $\Uo(\tL)$-modules. Hence $(P,\Uo(\rho_P))$ is a projective object in $\Mod(\Uo(\tL))$. 
The second statement follows from Proposition \ref{isofunctors} with $J=(0)$ and the Theorem is proved.
\end{proof}

%\begin{lemma} Let $R$ be an associative unital $k$ algebra with $k$ in the center $Z(R)$. It follows the ring of $k$-linear differential operators $\Diff(R)$ 
%is an almost commutative unital associative ring with $k \subseteq Z(\Diff(R))$.
%\end{lemma}
%\begin{proof} We prove $\Diff(R)$ is almost commutative. Let $\Diff^l(R) \subseteq \Diff(R)$ be the canonical filtration by degree of differential operators. Let $D\in \Diff^m(E)$ and
%$E\in \Diff^n(R)$ and assume $m=0, n=1$. We get $D(x):=ax$ for $x, a\in R$. By definition it follows $[E, \phi_b]=\phi_c$ for some $c\in R$. We
%get
%\[ [D\circ E, \phi_b]=D\circ E b-bE\circ D .\]
%By definition 
%\[ Eb =bE +\phi_c \]
%hence
%\[ D \circ Eb =D(bE+\phi_c)=(ab)E+a\phi_c.\]
%\end{proof}

%Note: It follows from Theorem \ref{mainDL} that projective objects in the cagetory $\Mod(\tL, Id)$ are connections $(P,\rho_P)$ that are direct summands of free $\Uo(\tL)$-modules.

\begin{definition} Let $(\tL, \ta, \tp,[,],D)$ be a $\Dl$-Lie algebra and let $J\subseteq \Uo(\tL)$ be a 2-sided ideal. Let $\Mod(\tL, Id, J)$ be the category
of $\tL$-connections $\rho: \tL \rightarrow \End_k(E)$ with $J$-curvature zero.
\end{definition}

\begin{corollary} \label{maincorrDL} Let $(\tL, \ta, \tp,[,],D)$ be a $\Dl$-Lie algebra and let $J\subseteq \Uo(\tL)$ be a 2-sided ideal. There is an exact equivalence of categories
\[ F: \Mod(\tL, Id, J) \cong \Mod(\Uo_J(\tL)) \]
with the property that an object $(P, \rho)\in \Mod(\tL, Id, J)$ is projective (resp. injective) if and only if $F(P,\rho)$ is projective (resp injective) in $\Mod(\Uo_J(\tL))$.
The pair $(\Uo_J(\tL), p^{\otimes}_J)$ represents the functor 
\[ \Conn_{\tL,J}: \ALG \rightarrow \SETS \]
hence the pair $(\Uo_J(\tL), p^{\otimes}_J)$ is unique up to unique isomorphism. The ring $\Uo_J(\tL)$ is an almost commutative ring.
\end{corollary}
\begin{proof} The proof is clear.
\end{proof}

Note: It follows from Corollary \ref{maincorrDL} that projective objects in the cagetory $\Mod(\tL, Id,J)$ are connections $(P,\rho_P)$ that are direct summands of free $\Uo_J(\tL)$-modules.
Because of Corollary \ref{maincorrDL} one wants in the case when $A$ is Noetherian and $\tL$ a finitely generated left $A$-module, to study the set of 2-sided ideals $J$ in the ring $\Uo(\tL)$ 
with the property that the quotient $\Uo_J(\tL)$ is Noetherian and almost commutative.

\begin{example} A surjection of $A\otimes_k A$-modules $\phi: \T^*_A(\tL) \rightarrow \Uo(\tL)$.\end{example}

Note: If $\{u_i\}_{i\in I}$ is a generating set for $\tL$ as left $A\otimes_k A$-module it follows the set $\{D, u_i\}_{i\in I}$ is a generating set for $\tL$ as left $A$-module:
Since any element $u \in \tL$ may be written as $u=\sum_i a_i(u_ib_i) $ we get
\[ u=\sum_i a_i(b_iu_i+\tp(u_i)(b_i)D= \sum_i (a_ib_i)u_i + bD\]
where $b:=\sum_i \tp(u_i)(b_i)$. 

Let $\tL$ be generated as left $A$-module by $u_1,\ldots, u_n$ where $u_1:=D$.
There is a canonical $A\otimes_k A$-module structure on $\T^*_A(\tL)$ where in the tensor product $\T^i_A(\tL)$ we have defined
\[  u_1\otimes \cdots \otimes u_ia\otimes u_{i+1}\otimes \cdots \otimes u_n= u_1\otimes \cdots \otimes u_i \otimes au_{i+1}\otimes \cdots \otimes u_n \]
and
\[  a(u_1\otimes \cdots \otimes u_i  \cdots \otimes u_n):= (au_1)\otimes \cdots \otimes u_i  \cdots \otimes u_n \]
and
\[  (u_1\otimes \cdots \otimes u_i  \cdots \otimes u_n)a:= u_1\otimes \cdots \otimes u_i  \cdots \otimes (u_na).\]
By induction the following holds:
There is a map of left $A\otimes_k A$-modules
\[ f: \T^*_A(\tL) \rightarrow \Uo(\tL) \]
defined by
\[ f(v_1\otimes \cdots \otimes v_n):=p^{\otimes}(v_1)\cdots p^{\otimes}(v_n) \in \Uo(\tL).\]
It is clear that $f$ is an $A\otimes_k A$-linear map. The abelian group $\Uo(\tL)$ is generated as left $A$-module by the set of products $u_{i_1}\otimes \cdots u_{i_l}$ for $i_a\in \{1,2,..,n\}$ and $l \geq 0$. This proves the surjectivity of the map $f$. The map $f$ is not a map of associative rings. Hence any 2-sided ideal $J\subseteq \Uo(\tL)$ lifts to a left $A\otimes_k A$-submodule $J_f\subseteq \T^*_A(\tL)$. In the case when $A$ is Noetherian,  we want to classify such lifts $J_f$ with the property that $\T^*_A(\tL)/J_f$ is a Noetherian $A\otimes_k A$-module. We also want to classify 2-sided ideals $I\subseteq \T^*_A(\tL)$ with Noetherian 
(or almost commutative) quotient $\T^*_A(\tL)/I$.

\begin{example} 2-sided ideals in $\T^*_A(\tL)$ and the $J$-curvature of a connection. \end{example}

When $J \subseteq \Uo(\tL)$ is a 2-sided ideal and $E$ is a left $\Uo_J(\tL)$ module, we get in a canonical way a left $A\otimes_k A$-module $J_f\subseteq \T^*_A(\tL)$, and the $A\otimes_k A$-module
$J_f$ is related to the $J$-curvature of the corresponding $\tL$-connection $\rho: \tL \rightarrow \End_k(E)$. Hence the study ot the set of 2-sided ideals in $\T^*_A(\tL)$ has applications to the study 
of the $J$-curvature of the connection $\rho$. The ring $\T^*_A(\tL)$ is a non-Noetherian ring in general.  If $\tL$ is a free module on the set $\{x_1,..,x_n\}$ and let $S(n)$ be the symmetric group on $n$ elements and let $G \subseteq S(n)$ be a subgroup.
Let  $J(G)$ be the 2-sided ideal generated by the set
\[ \{ x_1\otimes \cdots \otimes x_n-x_{\sigma(1)}\otimes \cdots \otimes x_{\sigma(n)}:\text{$i=1,..,n$ and $\sigma \in G$ } \} .\]
It follows $\T^*_A(\tL)/J(S(n))\cong \sym^*_A(\tL)\cong A[x_1,..,x_n]$ is the polynonial ring on $x_i$. Hence for this choice of $J$ it follows $\T^*_A(\tL)/J(S(n))$ is a Noetherian ring by Hilbert's basis theorem. 
One may ask if there is a more general form of the Hilbert basis theorem valid for quotients of the tensor algebra $\T^*_A(\tL)$, giving a classification of 2-sided ideals $J\subseteq \T^*_A(\tL)$ 
with Noetherian quotient $\T^*_A(\tL)/J$. 
%If we choose a strict subgroup $G\subseteq S(n)$, when will the quotient $\T^*_A(\tL)/J(G)$ be Noetherian? What about other 2-sided ideals?

In the case of the universal enveloping algebra $\U(L)$ of a semi simple Lie algebra $L$ over a field of characteristic zero, the \emph{space of primitive ideals} is a much studied object. See \cite{dixmier}
for some references. The study of the set of 2-sided ideals in enveloping algebras on the form $\U(A,L,f)$ is not well developed. Rings of differential operators $\Diff(A)$ have "few" 2-sided ideals: The Weyl algebra
may be realized as  the ring of differential operators on the polynomial ring, and the Weyl algebra is a simple ring - it has no nontrivial 2-sided ideals. It could be the set of 2-sided ideals in the universal ring $\Uo(\tL)$
in the case when $A$ is Noeteherian and $\tL$ is a finitely generated left $A$-module, is a "reasonably large space" and that it could be an interesting object to study. 
This study has by the above results applications to the study of the curvature of a connection.

\begin{example} \label{idealsinuniversal}  $\Dl$-Lie algebras and $A/k$-Lie-Rinehart algebras.\end{example}

Let $L:=\Der_k(A)$ and let $f\in \Z^2(\Der_k(A),A)$ be a 2-cocycle. Let $\tL:=L(f):=Az\oplus \Der_k(A)$ be the $\Dl$-Lie algebra induced by the pair $((\Der_k(A),Id) ,f)$.
 There is a canonical map
\[ \rho: Az \oplus \Der_k(A) \rightarrow \End_k(A) \]
defined by
\[ \rho(az+x):=\phi_a + x \]
where $\phi_a$ is multiplication with the element $a$. It follows $\rho$ is an $A\otimes_k A$-linear map. The map $\rho$ is in particular a $k$-linear map.
We get an induced map
\[ \T(\rho):\T^*_k(\tL ) \rightarrow \End_k(A) \]
defined by
\[\T(\rho)((a_1+x_1)\otimes \cdots \otimes (a_i+x_i))=(\phi_{a_1}+x_1)\circ \cdots \circ (\phi_{a_i}+x_i)\in \Diff^i(E).\]
of associative unital $k$-algebras. If $k$ is a field of characteristic zero and $A$ a regular $k$-algebra of finite type it follows $\T(\rho)$ is a surjective map of $k$-algebras with the property that
$\T(\rho)(J_1)=0$ where $J_1$ is the ideal from Definition \ref{ideals}. We get an induced exact sequence
\[ 0 \rightarrow J \rightarrow \Uo(\tL ) \rightarrow \Diff(A) \rightarrow 0 \]
where $J\subseteq \Uo(\tL)$ is a 2-sided ideal.
Let $J_3$ be the 2-sided ideal in $\T^*_k(\tL)$ generated by the elements $z-1, au-(az)\otimes u, ua-u\otimes az$ and $u\otimes v-v\otimes u-[u,v]$ for $a\in A, u,v\in \tL$.
If we define $\U(A,L,f):=\Uo(\tL)/J_3$ it follows left $\U(A,L,f)$-modules correspond to $L$-connections of curvature type $f$. We get an exact sequence
\[ 0 \rightarrow I \rightarrow \Uo(\tL) \rightarrow \U(A,L,f) \rightarrow 0\]
of associative unital rings.

\begin{theorem} \label{annihilatorideal} Let $k$ be a field of characteristic zero and let $A$ be a regular $k$-algebra of finite type. Let $L:=\Der_k(A)$ and $f\in \Z^2(L,A)$. Let $L(f)$ be the $\Dl$-Lie algebra
associated to $(L,f)$. We get exact sequences of rings 
\[ 0 \rightarrow J \rightarrow \Uo(\tL) \rightarrow \Diff(A) \rightarrow 0 \]
and
\[ 0 \rightarrow I \rightarrow \Uo(\tL) \rightarrow \U(A,L,f) \rightarrow 0.\]
Hence if we view an $\tL$-connection $(E,\rho)$ as a left $\Uo(\tL)$-module it follows $\rho$ has curvature type $f$ if and only if $I\subseteq ann(E,\rho)$, where $ann(E,\rho)$ is the annihilator ideal in 
$\Uo(\tL)$ of the $\Uo(\tL)$-module $E$. Similarly, $\rho$ is a flat connection if and only if $J\subseteq ann(E,\rho)$.
\end{theorem}
\begin{proof} The proof follows from the discussion above.
\end{proof}

By Theorem \ref{annihilatorideal}, it follows the annihilator ideal $ann(E,\rho)$ determines the curvature $R_{\rho}$ of the connection $\rho$. Hence if we want to study the curvature
$R_{\rho}$ we need to know the structure of the set of 2-sided ideals in the associative ring $\Uo(\tL)$.

\begin{example} \label{almostcommring} Almost commutative unital associative rings and the universal ring.\end{example}

Let 
\[ k \subseteq U_0 \subseteq U_1 \subseteq \cdots \subseteq U_i \subseteq \cdots \subseteq U \]
be a filtered associative unital ring $U$, where the multiplication is almost commutative. This means for any element $x\in U_i, y\in U_j$ it follows $[x,y]\in U_{i+j-1}$. It follows
the associated graded ring $Gr(U)$ is commutative. Assume $k\subseteq Z(U)$ is in the center of $U$. Let $\tL:=U_1$, $L:=U_1/U_0$ and $D:=1\in A$. We get an exact sequence
\begin{align}
&\label{AKLIE} 0 \rightarrow A \rightarrow \tL \rightarrow L \rightarrow 0 
\end{align}
and a canonical structure $\alpha:L\rightarrow \Der_k(A)$ and $\ta: \tL\rightarrow \Der_k(A)$ of $A/k$-Lie-Rinehart algebras on $\tL,L$. The sequence \ref{AKLIE} is an exact sequence of Lie-Rinehart algebras.
If we let $D:=1\in A$ it follows the $A\otimes_k A$-module $\tL$ is in a canonical way a pre-$\Dl$-Lie algebra $(\tL, \ta,[,],D)$ with the inclusion map $i:\tL \rightarrow U$ an $A\otimes_k A$-linear map. 
Hence there is a isomorphism of sets
\[ \Conn_{\tL}(U)\cong\Hom_{k-alg}(\Uo(\tL), U) \]
and we get a canonical map
\[ \Uo(i): \Uo(\tL) \rightarrow U \]
defined by
\[ \Uo(i)(u_1\otimes \cdots \otimes u_j):=i(u_1)\cdots i(u_j)\in U\]
commuting with the inclusion map $i: \tL \rightarrow U$. Hence if $U$ is generated as left $A$-module by $\tL$ and $A$-linear combinations of 
powers of elements of the form $u_1^{p_1}\cdots u_l^{p_l}$ with  $u_j\in \tL$,  it follows there is a 2-sided ideal $J\subseteq \Uo(\tL)$ and an isomorphism
\[ \Uo(\tL)/J\cong U .\]
 Hence $\Dl$-Lie algebras, pre-$\Dl$-Lie algebras and  the universal ring appears naturally if we want to study almost commutative rings and rings of generalized differential operators.

Note: If $A$ a regular commutative ring over a field of characteristic zero and $L$ an invertible $A$-module, it follows $\Diff(L)$ is an almost commutative ring generated by $\Diff^1(L)$. It follows there is a
pre-$\Dl$-Lie algebra $(\tL, \ta,[,],D)$ and an isomorphism $\Uo_J(\tL)\cong \Diff(L)$ of filered associative rings. Hence a left $\Diff(L)$-module corresponds to a $\tL$-connection $\rho$ with $J$-curvature zero.

\begin{lemma} \label{quotient} Let $A$ be a not neccessarily unital  associative ring and let $I,J$ be 2-sided ideals in $A$. Let $p_I:A\rightarrow A/I$ and $p_J:A\rightarrow A/J$ be the canonical projection maps.
There is a canonical isomorphism 
\[ A/I/p_I(J)\cong A/J/p_J(I) \]
of associative rings.
\end{lemma}
\begin{proof} Let $p_I:A\rightarrow A/I$ be the canonical projection map. It follows $p_I^{-1}(p_I(J))=I+J$. Clearly $I+J\subseteq p_I^{-1}(p(J))$. Assume $x\in p_I^{-1}(p(J))$. It follows $p_I(x)\in p_I(J)$.
Hence there is an element $y\in J$ with $p_I(x)=p_I(y)$ and it follows $p_I(x-y)=0$ hence $x-y\in ker(p_I)=I$. It follows $x=y+u$ with $u\in I$ and hence $x\in I+J$. It follows we get a canonial isomorphism
\[ A/I+J \cong A/I/p_I(J) \]
and this gives rise to a canonical isomorphism 
\[ A/I/p_I(J) \cong A/I+J \cong A/J/p_J(I) \]
of associative rings. The Lemma follows.
\end{proof}

\begin{lemma} \label{almostcommquot}Let $(U, U_i)$ be an almost commutative ring and let $I\subseteq U$ be a 2-sided ideal. Let $V:=U/I$ be the quotient. It follows $V$ is an almost commutative ring.
Let $W\subseteq U$ be an associative subring and let $W_i:=W\cap U_i$. It follows $(W,W_i)$ is an almost commutative ring.
\end{lemma}
\begin{proof} There is a canonical projection map 
\[ p: U \rightarrow V \]
and we define a filtration on $V$ with $V_i:=p(U_i)$. We get an increasing filtration $V_i$ on $V$ and it follows  $Gr(V)$ is commutative: If $x \in V_i:=p(U_i)$ and $y\in V_j:=p(U_j)$ 
with $x=p(u), y=p(v)$ it follows $uv-vu\in U_{i+j-1}$ hence $xy-yx\in V_{i+j-1}$. The second statement follows similarly and the Lemma is proved.
\end{proof}

We state a general result on properties of modules on and associative ring $A$.  Note: The following holds for non-unital rings as well.

\begin{lemma} \label{assnoetherian} Let $A$ be an associative unital ring and let $M' \subseteq M$ be a submodule with $p:M\rightarrow M/M'$ the canonical projection map.
Assume $N_1\subseteq N_2\subseteq M$ are sub modules. The module $M$ is Noetherian if and only if all sub modules $M' \subseteq M$ are finitely generated.
\begin{align}
&\label{No1}\text{If $p(N_1)=p(N_2)$ and $N_1\cap M'=N_2\cap M'$ it follows $N_1=N_2$.}\\
&\label{No2}\text{If $M'$ and $M/M'$ are Noetherian it follows $M$ is Noetherian.}\\
&\label{No3}\text{If $M\neq 0$ is Noetherian it follows $A/ann_A(M)$ is a Noetherian ring.}
\end{align}
\end{lemma}
\begin{proof} In \cite{mat} Theorem 3.1 and 3.5 the Lemma is stated and proved for commutative unital rings. Note that the  proof in \cite{mat} is valid for an arbitrary associative non-unital ring.
\end{proof}

Consider the following 2-sided ideal in $\T^*_k(\tL)$:
\[ I:=\{ u\otimes v -v\otimes u-[u,v]:\text{ such that $u,v\in \tL$} \} \subseteq \T^*_k(\tL) .\]
Let $I_i:=\T^*_k(\tL)_i \cap I$. It follows $I_i\subseteq I$ is a filtration of the 2-sided ideal $I$ that is compatible with the multiplication on $\T^*_k(\tL)$.

\begin{lemma} \label{lemmaalmost} Let $\T :=\T^*_k(\tL)$ and let $\T_j:=\T^*_k(\tL)_j$. For all $i,j\geq 0$ let $x:=x_1\otimes \cdots \otimes x_j, y:=y_1\otimes \cdots \otimes y_j$ with $x_a,y_b\in \tL$ for all
$1\leq a \leq i$ and $1\leq b \leq j$. The following holds: We may write
\[ x\otimes y - y \otimes x = \eta + \omega \]
with $\eta\in \T_{i+j-1}$ and $\omega\in I_{i+j}$.
\end{lemma}
\begin{proof} The claim is clearly true for $i=j=0$ or $i=0, j=1$. Assume $i=j=1$. We get
\[ x\otimes y - y\otimes x = [x,y]+(x\otimes y -y\otimes x -[x,y])=\eta +\omega \]
with $\eta=[x,y]$ and $\omega = x\otimes y -y\otimes x -[x,y]$. Here $\eta \in \T_1$ and $\omega \in I_{2}$. Hence the claim is true for $i=j=1$. Assume the claim is true for $i+j\leq k$ with $k\geq 1$ an integer.
Assume $i+j=k+1$. Let
\[ x:=x_1\otimes \cdots \otimes x_i \in \T_i \]
and
\[ y:=y_1\otimes \cdots \otimes y_j \in \T_j .\]
We get by the induction hypothesis
\[ x\otimes y = x \otimes y_1\otimes \cdots \otimes y_{j-1}\otimes y_j = (y_1\otimes \cdots \otimes y_{j-1}\otimes x +\eta_1 +\omega_1)\otimes y_j =\]
\[ y_1\otimes \cdots \otimes y_{j-1}\otimes x \otimes y_j +\eta_1 \otimes y_j +\omega_1\otimes y_j \]
with $\eta_1\otimes y_j \in \T_{i+j-1}$ and $\omega_1\otimes y_j \in I_{i+j}$.
Again by induction we get
\[ y_1\otimes \cdots \otimes y_{j-1}\otimes x \otimes y_j =\]
\[ y_1\otimes \cdots \otimes y_{j-1}\otimes (y_j \otimes x + \eta_2 +\omega_2) =\]
\[ y\otimes x + y_1\otimes \cdots \otimes y_{j-1}\otimes \eta_2 +y_1\otimes \cdots \otimes y_{j-1}\otimes \omega_2 \]
with $y_1\otimes \cdots \otimes y_{j-1}\otimes \eta_2\in \T_{i+j-1}$ and $y_1\otimes \cdots y_{j-1}\otimes \omega_2\in I_{i+j}$.
We get
\[ x\otimes y = y\otimes x + \eta + \omega \]
with
\[\eta=y_1\otimes \cdots \otimes y_{j-1}\otimes \eta_2+\eta_1\otimes y_j \in \T_{i+j-1} \]
and
\[y_1\otimes \cdots \otimes y_{j-1}\otimes \omega_2 +\omega_1 \otimes y_j \in I_{i+j}.\]
The Lemma follows.
\end{proof}

\begin{lemma} \label{almostlemma} Use the notation in Lemma \ref{lemmaalmost}. For any elements $u\in \T_i, v\in \T_j$ we may write
\[ u\otimes v-v\otimes u = \eta +\omega \]
with  $\eta\in \T_{i+j-1}$ and $\omega \in I_{i+j}$.
\end{lemma}
\begin{proof} We may write $u=x_1+x$ with $x_1\in \T_{i-1}$ and $x\in \tL^{\otimes_k i}$ and similar $v=y_1+y$ with $y_1\in \T_{j-1}, y\in \tL^{\otimes_k j}$.
We get
\[ u\otimes v-v\otimes u =(x_1+x)\otimes (y_1+y)-(y_1+y)\otimes (x_1+x)=\]
\[x_1\otimes y_1+x_1\otimes y + x\otimes y_1 +x\otimes y -( y_1\otimes x_1+y_1\otimes x + y\otimes x_1+y\otimes x) =\]
\[ x_1\otimes y_1+x_1\otimes y + x\otimes y_1-(y_1\otimes x_1+y_1\otimes x + y\otimes x_1) +x\otimes y -y\otimes x =\]
\[ z+x\otimes y -y\otimes x \]
with
\[ z:= x_1\otimes y_1+x_1\otimes y + x\otimes y_1-(y_1\otimes x_1+y_1\otimes x + y\otimes x_1)  \in \T_{i+j-1}.\]
By Lemma \ref{lemmaalmost} we get
\[ x\otimes y -y\otimes x = \eta+\omega \]
with $\eta\in\T_{i+j-1}$ and $\omega \in I_{i+j}$. It follows 
\[ u\otimes v-v\otimes u= \eta + z+ \omega \]
where $\eta+z\in \T_{i+j-1}$ and $\omega\in I_{i+j}$ and the Lemma follows.
\end{proof}

\begin{definition} Use the notation in Lemma \ref{lemmaalmost}. Let $I\subseteq \T$ be the 2-sided ideal $I:=\{ u\otimes v -v\otimes u -[u,v]: u,v\in \tL\}$. Define
$\U(\tL):= \T/I =\T^*_k(\tL)/I$. Let $p:\T \rightarrow \U(\tL)$ be the canonical projection map and let $\U_i(\tL):=p(\T^*_k(\tL)_i)$ for $\i \geq 0$.
There are canonical projection maps $p^{\otimes}: \T\rightarrow \Uo(\tL)$ and $p^{\rho}:\T \rightarrow \Ur(\tL)$. Define
\[ \tUo(\tL):= \Uo(\tL)/p^{\otimes}(I) \]
and
\[ \tUr(\tL):=\Ur(\tL)/p^{\rho}(I).\]
Since $\tUo(\tL)$ and $\tUr(\tL)$ are quotients of $\T$ we get filtrations $\tUo_i(\tL) \subseteq \tUo(\tL)$ and $\tUr_i(\tL) \subseteq \tUr(\tL)$
for all integers $i\geq 1$.
\end{definition}

We get a filtration on the associative ring $\U(\tL)$:
\[ \U_1(\tL) \subseteq \U_2(\tL) \subseteq \cdots \subseteq \U_i(\tL) \subseteq \U(\tL) \]
compatible with the multiplication: For any elements $x\in \U_i(\tL), y\in \U_j(\tL)$ it follows $xy\in \U_{i+j}(\tL)$.
A similar property holds for the filtrations $\tUo_i(\tL)$ and $\tUr_i(\tL)$.

\begin{proposition} \label{almostcomm} The associative ring $\U(\tL)$ is an almost commutative ring.
\end{proposition}
\begin{proof} Use the notation from Lemma \ref{lemmaalmost}. There is a canonical projection map
\[ p:\T\rightarrow \U(\tL) \]
with $\U_i(\tL):=p(\T_i)$. Let $ x\in \U_i(\tL), y\in \U_j(\tL)$ with $x=p(u), y=p(v)$. It follows from Lemma \ref{almostlemma} we may write
\[ u\otimes v-v\otimes u= \eta + \omega \]
with $\eta\in \T_{i+j-1}, \omega \in I_{i+j}$. It follows $xy-yx =p(u\otimes v-v\otimes u)=p(\eta+\omega)=p(\eta)\in \U_{i+j-1}(\tL)$ hence $\U(\tL)$ is almost commutative.
The Proposition follows.
\end{proof}

\begin{corollary} The rings $\tUo(\tL)$ and $\tUr(\tL)$ are almost commutative associative unital rings.
\end{corollary}
\begin{proof} By Lemma \ref{quotient} we may do the following: There is by definition  inclusions $I, J_1 \subseteq \T^*_k(\tL)$ and there is a canonical quotient map
\[ p^{\otimes }: \T^*_k(\tL)  \rightarrow \Uo(\tL) .\]
By definition we have
\[ \tUo(\tL):= \Uo(\tL)/p^{\otimes }(I).\]
By Lemma \ref{quotient} we may write
\[ \tUo(\tL)\cong \U(\tL)/p(J_1) \]
where
\[ p: \T^*_k(\tL) \rightarrow \U(\tL) \]
is the canonical projection map. It follows $\tUo(\tL)$ is a quotient of $\U(\tL)$ which by Proposition \ref{almostcomm} is almost commutative. It follows
$\tUo(\tL)$ is almost commutative. A similar argument shows $\tUr(\tL)$ is almost commutative. The rings $\tUo(\tL)$ and $\tUr(\tL)$ are quotients of associative unital rings by 2-sided ideals. Hence
it follows $\tUo(\tL)$ and $\tUr(\tL)$ are almost commutative associative unital rings. The Corollary follows.
\end{proof}

Let in the following $\tL:=A\{u_1,..,u_n\}$ be a finite generating set of the $\Dl$-Lie algebra $\tL$ as left $A$-module. Make the following definitions:

\begin{definition} Let $B_i:=A\{ u_{j(1)}\otimes \cdots \otimes u_{j(i)}\text{such that $u_{j(a)}\in \{u_1,..,u_n\}$} \}$. Let $B^i:=B_1\oplus \cdots \oplus B_i \subseteq \T^*_k(\tL)$.
\end{definition}

\begin{lemma} Let $(\tL,\ta, \tp,[,],D)$ be a $\Dl$-Lie algebra with generating set $\{u_1,..,u_n\}$ as left $A$-module. The following holds: For any element $x\in \T^*_k(\tL)_i$ we may write
$x=x_1+x_2$ where $x_1\in B^i \subseteq \T^*_k(\tL)_i$ and $x_2\in I_i \subseteq \T^*_k(\tL)_i$.
\end{lemma}
\begin{proof} The claim is obvious for $i=0,1$. Let $i=2$ and let $au_i\otimes bu_j\in \tL^{\otimes ^2_k}$. We get
\[ au_i\otimes bu_j=a( u_i\otimes bu_j-bu_j\otimes u_i-[u_i,bu_j] +bu_j\otimes u_i+[u_i,bu_j])=\]
\[ abu_j\otimes u_i +a[u_i,bu_j]+ a\omega \]
where
\[ \omega :=u_i\otimes bu_j-bu_j\otimes u_i -[u_i,bu_j])\in I_2.\]
Define $v_1:=a[u_i,bu_j]\in \tL$ and $x_2:=a\omega \in I_2$. We may write $v_1=\sum c_iu_i$ and define $x_1:=abu_j\otimes u_i +v_1\in B^2$.
We get
\[ au_i\otimes bu_j =x_1+x_2\]
with $x_1\in B^2$ and $x_2\in I_2$ and the claim is true for $i=2$. Assume the claim is true for $x\in \T^*_k(\tL)_{i-1}$. Let $x\in \T^*_k(\tL)_i$.
we may write $x=u+x_1$ witth $u\in \T^*_k(\tL)_{i-1}$ and $x_1\in \tL^{\otimes^i_k}$. By induction we get $u=u_1+u_2$ where $u_1\in B^{i-1}$ and $u_2\in I_{i-1}\subseteq I_i$.
we may write $x_1$ as s sum of elements on the form $a_1v_1\otimes \cdots a_iv_i$ with $a_j\in A$ and $v_j\in\{u_1,..,u_n\}$. 
We get
\[ a_{i-1}v_{i-1}\otimes a_iv_i=a_{i-1}(v_{i-1}\otimes a_iv_i+a_iv_i\otimes v_{i-1}-[v_{i-1},a_iv_i])+ \]
\[ a_{i-1}(a_iv_i\otimes v_{i-1}+[v_{i-1},a_iv_i])=\]
\[ a_{i-1}a_iv_1\otimes v_{i-1}+a_{i-1}w_i+a_{i-1}[v_{i-1},a_iv_i]\]
with 
\[ w_i:= v_{i-1}\otimes a_iv_i+a_iv_i\otimes v_{i-1}-[v_{i-1},a_iv_i] \in I_2.\]
It follows 
\[x_1=a_1v_1\otimes \cdots \otimes a_{i-2}v_{i-2}\otimes a_{i-1}(a_iv_i\otimes v_{i-1}+w_i+[v_{i-1},a_iv_i]).\]
Hence
\[x_1=a_1v_1\otimes \cdots a_{i-2}v_{i-2}\otimes a_{i-1}a_iv_i\otimes v_{i-1}+ \]
\[ a_1v_1\otimes \cdots a_{i-2}v_{i-2}\otimes a_{i-1}w_i + a_1v_1\otimes \cdots a_{i-2}v_{i-2}\otimes a_{i-1}[v_{i-1},a_iv-i]=\]
\[ a_1v_1\otimes \cdots a_{i-2}v_{i-2}\otimes a_{i-1}a_iv_i\otimes v_{i-1}+\omega_1+\eta_1 \]
where
\[ \omega_1:= a_1v_1\otimes \cdots a_{i-2}v_{i-2}\otimes a_{i-1}w_i \]
and
\[ \eta_1:= a_1v_1\otimes \cdots a_{i-2}v_{i-2}\otimes a_{i-1}[v_{i-1},a_iv-i].\]
It follows $\omega_1\in I_i$ and $\eta_1\in \T^*_k(\tL)_{i-1}$.
We may write as follows:
\[ a_1v_1\otimes \cdots a_{i-2}v_{i-2}\otimes a_{i-1}a_iv_i\otimes v_{i-1}= \eta \otimes v_{i-1} \]
with $\eta'\in \T^*_k(\tL)_{i-1}$ hence it follows $\eta'=\eta'_1+\eta'_1$ where $\eta'_1\in B^{i-1}$ and $\eta'_2\in I_{i-1}$. It follows
\[ a_1v_1\otimes \cdots a_{i-2}v_{i-2}\otimes a_{i-1}a_iv_i\otimes v_{i-1}= \eta'_1\otimes v_{i-1}+\eta'_2\otimes v_{i-1}\]
with $\eta'_1\otimes v_{i-1}\in B^i$ and $\eta'_2\otimes v_{i-1}\in I_i$. It follows
\[ x_1=\eta'_1\otimes v_{i-1}+\eta'_2\otimes v_{i-1} +\omega_1 +\eta_1 \]
with $\eta_1\in \T^*_k(\tL)_{i-1}$. we may write $\eta_1=a_1+a_2$ with $a_1\in B^{i-1}$ and $a_2\in I_{i-2}$. It follows
\[x_1= \eta + \omega \]
where
\[\eta:=\eta'_1\otimes v_{i-1}+a_1\in B^i\]
and
\[\omega:= \eta'_2\otimes v_{i-1}+a_2 \in I_i\]
For a general $x_1\in \tL^{\otimes i}$ we may do sometthing similar: $x_1=x'_1+x'_2$ with $x'_1\in B^i$ and $x'_2\in I_i$. It follows
\[ x=u_1+x'_1+ u_2+x'_2 \]
with $u_1+x'_1\in B^i$ and $u_2+x'_2\in I_i$ and the claim of the Lemma follows.
\end{proof}

\begin{corollary} \label{fingen} Let $(\tL, \ta,tp,[,],D)$ be a $\Dl$-Lie algebra where  $A$ is a Noetherian ring and where $\tL$ is finitely generated as left $A$-module. 
Let $\T:=\T^*_k(\tL)^1$ and let $I:=\{u\otimes v-v\otimes u-[u,v]:u,v\in\tL\} \subseteq \T^*_k(\tL)$. Let $p:\T\rightarrow \U(\tL):=\T/I$.
Let $\U_i(\tL)\subseteq \U(\tL)$ be the canonical filtration. It follows the quotient $\U_i(\tL)/\U_{i-1}(\tL)$ is a finitely generated left $A$-module for all $i\geq 1$
\end{corollary}
\begin{proof} Let $x  \in \U_i(\tL)/\U_{i-1}(\tL)$ be an element. We may write $x=x_1+x_2$ where $x_1=p(u)$ with $u\in B^i$ and $x_2=p(v)$with $v\in I_i$.
We get 
\[ x=p(u+v)=p(u)=\sum_{(j(1),..,j(i))}a_{( j(1),..,j(i) )}u_{j(1)}\bullet \cdots \bullet u_{j(i)} \]
hence $\U_i(\tL)/\U_{i-1}(\tL)$ is generated by the finite set 
\[ \{u_{j(1)}\bullet \cdots \bullet u_{j(i)}\text{ for $j(a)\in \{1,..,n\}$  } \}.\]
 The claim follows.
\end{proof}

\begin{theorem} \label{mainnoetherian} Let $(\tL, \ta, \tp,[,],D)$ be a $\Dl$-Lie algebra where $A$ is Noetherian and $\tL$ is finitely generated as left $A$-module. It follows the rings $\tUo(\tL)$ and $\tUr(\tL)$
are Noetherian rings.
\end{theorem}
\begin{proof} Let $J_1,J_2\subseteq \T^*_k(\tL)$ be the ideals defining the rings  $\Uo(\tL)$ and $\Ur(\tL)$. Let $p:\T:=\T^*_k(\tL)\rightarrow \T/I=\U(\tL)$ where
$I$ is the 2-sided ideal generated by the elements $u\otimes v-v\otimes u-[u,v]$ for $u,v\in \tL$. Let $\tUo_i(\tL) \subseteq \Uo(\tL)$ be the canonical filtration for $i\geq 1$.
We get since the map $p$ is surjective a canonical surjective map of left $A$-modules
\[ \Uo_i(\tL)/\Uo_{i-1}(\tL)\rightarrow \tUo_i(\tL)/\tUo_{i-1}(\tL) \]
and since $\Uo_i(\tL)/\Uo_{i-1}(\tL)$ is a finitely generated left $A$-module it follows the module $\tUo_i(\tL)/\tUo_{i-1}(\tL)$ is finitely generated as left $A$-module. The associated graded ring
$Gr(\tUo(\tL))$ is generated by $\tUo_1(\tL)$ as $A$-algebra. Hence $Gr(\tUo(\tL))$ is a Noetherian ring. By \cite{ginzburg}, Proposition 1.1.6 it follows $\Uo(\tL)$ is a Noetherian ring.
A similar argument proves $\Ur(\tL)$ is a Noetherian ring and the Theorem follows.
\end{proof}

Recall the following: Let $\Mod(\tL, Id)$ denote the category of pairs $(\rho,E)$ where $\rho$ is a $B:=A\otimes_k A$-linear map $\rho:\tL \rightarrow \End_k(E)$ such that $\rho(D)=Id_E$ and where morphisms are defined as follows: Given two element $(\rho,E), (\rho',E')\in \Mod(\tL,Id)$, a morphism $\phi:(\rho,E)\rightarrow (\rho',E')$ is an $A$-linear map $\phi: E\rightarrow E'$ such that for all elements $u\in \tL$ 
it follows $\rho'(u)\circ \phi=\phi \circ \rho(u)$. 

Let $\Conn(\tL, Id)$ denote the category of pairs $(\rho,E)$ where $\rho:\tL \rightarrow \End_k(E)$ is a left $A$-linear map such that the following holds for all $u\in \tL, a\in A$ and $e\in E$:
\[ \rho(u)(ae)=a\rho(u)(e)+\tp(u)(a)e.\]
A morphism $\phi:(\rho,E)\rightarrow (\rho',E')$ in $\Conn(\tL, Id)$ is an $A$-linear map $\phi:E\rightarrow E'$ such that for any element $u\in \tL$ it follows $\rho'(u)\circ \phi=\phi \circ \rho(u)$.

\begin{theorem} \label{mainequiv} There is an exact equivalence of categories
\begin{align}
&\label{F2} F_1: \Conn(\tL, Id) \cong \Mod(\Ur(\tL))
\end{align}
with the property that $F_2$ preserves injective and projective objects.
\end{theorem}
\begin{proof} Ths proof is similar to the proof for $\Uo(\tL)$ and is left to the reader as an exercise.
\end{proof}

\begin{definition} \label{cohomology}  Let $(\tL,\ta,\tp,[,],D)$ be a $\Dl$-Lie algebra and let $(L,\alpha)$ be a classical Lie-Rinehart algebra.
Let $f\in \Z^2(\Der_k(A),A)$ is a 2-cocycle.  Let furthermore $(\rho_U, U),(\rho_V,V)$ be objects in $\Mod(\tL, Id)$ and let $(\rho_W,W),(\rho_Z,Z)$ be objects in $\Conn(L(\fal), Id)$. 
By Theorem \ref{mainequiv} there are exact equivalences of categories
\[ \Mod(\tL, Id)\cong \Mod(\Uo(\tL)) \]
and
\[ \Conn(L(\fal), Id)\cong \Mod(\Ur(L(\fal))) \]
preserving injective and projective objects.
Since the categories $\Mod(\Uo(\tL))$ and $\Mod(\Ur(L(\fal)))$ have enough injectives we may define the $\Ext$ and $\Tor$-groups of $U,V,W,Z$. Let in the following
$\Uo:= \Uo(\tL)$ and $\Ur:=\Ur(L(\fal))$. We may define the groups
\[ \Ext^i_{\Uo}(U,V), \Tor_i^{\Uo}(U,V) \]
and
\[  \Ext^i_{\Ur}(Z,W), \Tor_i^{\Ur}(Z,W) .\]
Let 
\[ \H_i(\tL, V):= \Tor_i^{\Uo}(A,V) \]
\[ \H^i( \tL, V):= \Ext^i_{\Uo}(A,V).\]
and
\[ \H^i(L(\fal), W):= \Ext^i_{\Ur}(A,W).\]
\[ \H_i(L(\fal), W):= \Tor_i^{\Ur}(A,W) \]
\end{definition}

Note: $\Ext$ and $\Tor$-groups over an associative ring $R$ are modules over the center $Z(R)$. Hence the cohomology and homology groups defined in \ref{cohomology} are $k$-modules and not $\Uo(\tL)$ or
$\Ur(\tL)$-modules. Hence with Definition \ref{cohomology} the cohomology and homology groups do not have naturally defined connections. Classically when considering families of varieties we get induced 
connections on higher direct image sheaves - Gauss-Manin connections - and the problem of defining Gauss-Manin connections for the $\Ext$ and $\Tor$-groups in Definition \ref{cohomology} 
will be investigated in a coming paper.

\begin{example} Cohomology of modules over almost commutative rings and connections.\end{example}

If $U$ is an almost commutative associative unital ring as in example \ref{almostcommring} and $E$ is a left $U$-module, there is an isomorphism
$U\cong \Uo_J(\tL)$ where $J\subseteq \Uo(\tL)$ is a 2-sided ideal. Hence the module $E$ may be viewed as an $\tL$-connection
\[ \rho: \tL \rightarrow \End_k(E) \]
with $J$-curvature equal to zero. The $\Ext$ and $\Tor$-groups of $E$ as left $U$-module may be interpreted as $\Ext$ and $\Tor$-groups of the corresponding $\tL$-connection 
$(E,\rho)$. Almost commutative rings and their modules is a much studied topic in non-commutative ring theory. A module $E$ on an almost commutative ring $U$ becomes more "geometric" if we view
$E$ as a connection $\rho: \tL\rightarrow \End_k(E)$ where $\tL:=U_1$. In the case when the left $U$-module $E$ is a finite rank projective $A$-module where $A:=U_0$, 
it follows the pair $(E,\rho)$ may be viewed as an algebraic connection on an algebraic vector bundle on the affine scheme $\Spec(A)$.

\begin{example} Hochschild cohomology of left and right $\Uo(\tL)$ and $\Ur(\tL)$-modules. \end{example}

If $B$ is an algebra over a field or a commutative ring $R$ such that $B$ is projective as left $R$-module, it follows the $\Ext$ and $\Tor$-groups may be calculated by the Hochschild cohomology
and homology groups of certain left and right modules. If $(E,\rho)$ and $(E',\rho')$ are left $\Uo(\tL)$-modules it follows from \cite{weibel}, Lemma 9.1.9 that
\[ \Ext^i_{\Uo(\tL)/k}((E,\rho),(E',\rho'))\cong \H^i(\Uo(\tL), \Hom_k((E,\rho),(E',\rho')) \]
where $\Hom_k((E,\rho),(E',\rho))$ is the left and right $\Uo(\tL)$-module of $k$-linear maps between $(E,\rho)$ and $(E',\rho')$. Here $\H^i$ is the Hocschild cohomology of the $\Uo(\tL)$-bimodule
$\Hom_k((E,\rho),(E',\rho'))$. Hence there is an explicit complex calculating the $\Ext$-groups from Definition \ref{cohomology}.

\begin{example}\label{noetherianconn}  Cohomology and homology of a Noetherian module \end{example}

\begin{definition} \label{annihilator} Let $(\tL, \ta, \tp, [,],D)$ be a $\Dl$-Lie algebra and let $(\rho,E)\in \Mod(\tL, Id)$.
Let $a^{\otimes}(\rho,E)\subseteq \Uo(\tL)$ be the annihilator ideal of $(\rho,E)$ and define $\Uo_E(\tL):=\Uo(\tL)/a^{\otimes}(\rho,E)$. Define similarly 
$\Ur_E(\tL):=\Ur(\tL)/a^{\rho}(\rho,E)$ where $a^{\rho}(\rho,E)$ is the annihilator ideal of $(\rho,E)$ as $\Ur(\tL)$-module.
\end{definition}

\begin{proposition} \label{univnoetherian} Assume $A$ is a Noetherian ring and $(\rho,E)\in \Mod(\tL,Id)$ a connection with the property that $E$ is a finitely generated $A$-module. It follows $(\rho,E)$ is a Noetherian
$\Uo(\tL)$-module. The ring $\Uo_E(\tL)$ is Noetherian.
\end{proposition}
\begin{proof} Let $(\rho',E')\subseteq (\rho,E)$ be an $\Uo(\tL)$-sub module. It follows $E'\subseteq E$ is a sub-$A$-module and since $A$ is Noetherian it follows $E'$ is a finitely generated $A$-module.
It follows $(\rho',E')$ is a finitely generated $\Uo(\tL)$-module, hence by Lemma \ref{assnoetherian} it follows $(\rho,E)$ is a Noetherian $\Uo(\tL)$-module. 
Again from Lemma \ref{assnoetherian}, \ref{No3} it follows the ring $\Uo_E(\tL):=\Uo(\tL)/a^{\otimes}(\rho,E)$ is Noetherian. The Proposition follows.
\end{proof}

The rings $\Uo(\tL)$ and $\Ur(\tL)$ from Theorem \ref{mainequiv} are non-Noetherian in general but from Proposition \ref{univnoetherian} we can in many cases reduce to the Noetherian 
case when studying cohomology and homology of a connection $(\rho,E)$ over a Noetherian ring $A$ when $E$ is a finitely generated $A$-module,
We may study $(\rho,E)$ as $\Uo_E(\tL)$ or $\Ur_E(\tL)$-module and the  rings $\Uo_E(\tL)$ and $\Ur_E(\tL)$ are by Proposition \ref{univnoetherian}  Noetherian rings.
There is moreover an explicit complex calculating the $\Ext$ and $\Tor$-groups in Definition \ref{cohomology}.

Note: If $A$ is Noetherian and we are given a finite family $(\rho_i, E_i)_{i\in I}$ of $\tL$-connections with $E_i$ a finitely generated $A$-module for all $i$ and we let $E:=\oplus E_i$
it follows $ann(\rho, E)\subseteq ann(\rho_i ,E_i)$ for all $i$. The $A$-module $E$ is a finitely generated $A$-module hence $\Uo_E(\tL)$ and $\Ur_E(\tL)$ are by Proposition \ref{univnoetherian}
Noetherian and $E_i$ are left $\Uo_E(\tL)$ and $\Ur_E(\tL)$-modules for all $i$. Hence when studying cohomology or homology a finite set of connections $E_i$ we may always work over a fixed Noetherian ring $\U$.

\begin{example} Differential operators, connections and projective modules.\end{example}

Recall the notion of a \emph{projective basis} for a finitely generated projective $A$-module $E$, where $A$ is a commutative unital ring over a base ring $k$ (see \cite{maa13}).
A set of $r$ elements $x_1,..,x_r \in E^*$ and $e_1,..,e_r\in E$ satisfying the formula
\[ \sum_i x_i(e)e_i=e \]
for all $e\in E$ is a projective basis. The $A$-module $E$ has a projective basis if and only if it is finitely generated and projective as $A$-module. If $A$ is a finitely generated algebra over a field
and $E$ is a finitely generated and projective $A$-module,  a projective basis for $E$ may be calculated using a \emph{Gr\"{o}bner basis}. 

\begin{definition} Let $k\rightarrow A$ be an arbitrary map of unital commutative rings and let $E$ be a left $A$-module. Define $\Diff^0(A):=A$ and $\Diff^0(E):=\End_A(E)$. An operator
$D \in \End_k(E)$ is a differential operator of order $\leq l$ if for all sequences of $l+1$ elements $a_1,..,a_{l+1} \in A$ the following holds:
\[ [\cdots [D,a_1 Id_E]\cdots ], a_{l+1}Id_E]=0 .\]
Here $[,]$ is the Lie product on $\End_k(E)$. We let $\Diff^l_k(E)$ (or for short $\Diff^l(E)$ be the set of all differential operators of order $\leq l$. Let $\Diff(E):=\cup_{l\geq 0}\Diff^l(E)$.
\end{definition}

It follows the $k$-vector space  $\Diff(E)$ has a filtration of $k$-vector spaces
\[ \Diff^0(E)=\End_A(E) \subseteq \Diff^1(E)  \subseteq \cdots \subseteq \Diff^l(E)  \subseteq \cdot \Diff(E).\]
The composition of operators gives for all integers $l,l'\geq 0$ an associative product
\[ \Diff^l(E) \times \Diff^{l'}(E) \rightarrow \Diff^{l+l'}(E) \]
making $\Diff(E)$ into an associative ring. Similar properties hold for the ring $A$: $\Diff(A)$ is an associative ring containing $A$. The ring $A$ does not lie in the centre of $\Diff(A)$.
The product on $\Diff(A)$ respects the filtration, hence we get a well defined product
\[  \Diff^l(A) \times \Diff^{l'}(A) \rightarrow \Diff^{l+l'}(A) .\]
It follows $\Diff^l(A)$ and $\Diff^l(E)$ are left and right $A$-modules. It follows $\Diff^l(A)$ and $\Diff^l(E)$ are left $A\otimes_k A$-modules for all integers $l\geq 0$.
Given an operator $D\in \Diff^l(A)$ or $\Diff^l(E)$ the action of $A\otimes_k A$ is defined as follows:
\[ a\otimes b.D:= \phi_a \circ D \circ \phi_b \]
where $\phi_a$ is multiplication with the element $a$.

Define the following map:

\[ \rho: \Diff_k(A) \rightarrow \End_k(E) \]
by
\[ \rho(D)(e) := \sum_i D(x_i(e))e_i .\]

\begin{lemma} \label{hconn} The map $\rho$ is a map of $A\otimes_k A$-modules. It induces a map
\[ \rho^l:\Diff^l(A) \rightarrow \Diff^l(E).\]
\end{lemma}
\begin{proof} The action of $A\otimes_k A$ on $\Diff(E)$ is as follows $a\otimes b \psi:= \phi_a \circ \psi \circ \phi_b$ where $\phi_a$ is multiplication with $a$. 
One checks $\rho(D)\in \End_k(E)$. Moreover $\rho(D+D')=\rho(D)+\rho(D')$.
Let $a\otimes b \in A\otimes_k A$. We get
\[ \rho(a\otimes b.D)(e):= \sum_i aD(bx_i(e)) e_i=\]
\[ \sum_i aD(x_i(be))e_i =a\rho(D)(be) := (a\otimes b.\rho(D))(e).\]
Hence $\rho(a\otimes b.D)=a\otimes b.\rho(D)$. It follows $\rho$ gives an $A\otimes_k A$-linear map as claimed.
Assume $D\in \Diff^l(A)$. We need to prove that $\rho(D)\in \Diff^l(E)$. Let $a_1,..,a_{l+1}\in A$ be $l+1$ arbitrary elements. We get
\[  [\cdots [\rho(D), a_1Id_E] \cdots ] , a_{l+1}Id_E](e)=\]
\[ \sum_i [\cdots [D, a_1Id_A]\cdots ] , a_{l+1}Id_A](x_i(e))e_i \]
and since $D\in \Diff^l(A)$ it follows 
\[ [\cdots [D, a_1Id_A]\cdots ], a_{l+1}Id_A]=0 .\]
It follows
\[  [\cdots [\rho(D) , a_1Id_E] \cdots ] , a_{l+1}Id_E](e)=0 \]
for all $e\in E$. It follows $\rho(D)\in \Diff^l(E)$. The Lemma follows.
\end{proof}

We get for any pair of integers $k,l\geq 0$ a map of $k$-vector spaces
\[ R^{k,l}_{\rho}: \Diff^k(A)\otimes_k  \Diff^l(A) \rightarrow \Diff^{k+l}(E) \]
 defined by
\[ R^{k,l}_{\rho}(D,D'):=  \rho^{k+l}(D\circ D')-   \rho^k(D)\rho^l(D') .\]

\begin{definition} Let the map $\rho^l$ from Lemma \ref{hconn} be an $l$-connection on $E$. Let $\rho$ be an $\infty$-connection on $E$. Let $R^{k,l}_{\rho}$ be the $(k,l)$-curvature of $\rho$.
\end{definition}

Note: An ordinary connection $\nabla: \Der_k(A)\rightarrow \Diff^1(E)$ has curvature
\[ R_{\nabla} \in \Hom_A( \wedge^2 \Der_k(A),\End_A(E)).\]
The $(k,l)$-connection $R^{k,l}_{\rho}$ does not satisfy a similar property. 

In the case when $k=l=1$ we get a map of $A\otimes_k A$-modules
\[ \rho^1: \operatorname{D}^1(A,0) \rightarrow \Diff^1(E) .\]
We also get a map of left $A$-modules
\[ \nabla: \Der_k(A) \rightarrow \End_k(E).\]
The map $\nabla$ is a connection on $E$.

Note: By the paper \cite{maa0} it follows the connections $\rho^1$ and $\nabla$  are non-flat in general. There is an explicit formula for the curvature $R_{\rho^1}$ of $\rho^1$  and $\nabla$ in terms of an idempotent $\phi$ for the module $E$. One uses the projective basis $x_i,e_j$ for $E$  to define an idempotent and the formula for the curvature involves the idempotent $\phi$. See Theorem 2.14 in \cite{maa0} for a proof of the formula and some explicit examples. Given a projective basis $x_1,..,x_r,e_1,..,e_r$ for $E$ we get a surjection $p: A\{u_1,..,u_r\} \rightarrow E$ defined by $p(u_i)=e_i$. It follows $A^r/ker(p)\cong E$. Define the following matrix $\phi \in \End_A(A^r)$:
\[
\phi:=\begin{pmatrix}  x_1(e_1) & x_1(e_2) & \cdots & x_1(e_r) \\
              x_2(e_1)  & x_2(e_2) & \cdots  & x_2(e_r) \\
             \vdots & \vdots  & \cdots           & \vdots \\
             x_r(e_1) & x_r(e_2) & \cdots & x_r(e_r) 
\end{pmatrix} 
\]
Given two derivations $\delta,\eta \in \Der_k(A)$ we may consider the matrix $\delta(\phi):= (\delta(x_i(e_j)) )\in \End_A(A^r)$. The Lie product $[\delta(\phi),\eta(\phi) ]\in \End_A(A^r)$. By \cite{maa0} it follows
\[   [\delta(\phi),\eta(\phi)](ker(p)) \subseteq ker(p) \]
hence the matrix $[\delta(\phi),\eta(\phi)]$ induces an endomorphism of $E$. 
\begin{theorem}\label{curvform} The following holds:
\begin{align}
&\label{curvature} R_{\nabla}(\delta,\eta) =[\delta(\phi),\eta(\phi)] \in \End_A(E).
\end{align}
\end{theorem}
\begin{proof}
For a proof see \cite{maa0}, Theorem 2.14. 
\end{proof}

From formula \ref{curvature}  we see that given a projective basis $x_i,e_j$ it follows the corresponding connection $\nabla$ is seldom flat since the Lie product $[\delta(\phi),\eta(\phi)]$ is seldom zero as an element 
of $\End_A(E)$.

\begin{lemma}  By Lemma \ref{curvrho} it follows the map $\rho^1: \operatorname{D}^1(A,0)\rightarrow \Diff^1(E)$ 
has curvature $R_{\rho^1}$ defined by $R_{\rho^1}(u,v)=\rho^1([u,v])-[\rho^1(u),\rho^1(v)]$ for $u,v \in \operatorname{D}^1(A,0)$.
It follows 
\[ R_{\rho^1}(u,v)=R^{1,1}_{\rho}(u,v)-R^{1,1}_{\rho}(v,u).\]
\end{lemma}
\begin{proof} We get since $\rho$ is defined for operators in $\Diff^2(A)$ the following calculation:
\[ R_{\rho^1}(u,v)=\rho^1([u,v])-[\rho^1(u),\rho^1(v)]=\rho^2(uv-vu)-(\rho^1(u)\rho^1(v)-\rho^1(v)\rho^1(u))=\]
\[ \rho^2(uv)-\rho^1(u)\rho^1(v)-(\rho^2(vu)-\rho^1(v)\rho^1(u))=R^{1,1}_{\rho}(u,v)-R^{1,1}_{\rho}(v,u).\]
The Lemma follows.
\end{proof}

Hence $R_{\rho}$ is determined by $R^{1,1}_{\rho}$.

\begin{lemma} \label{rings} The map $\rho$ is a morphism of rings if and only if $R^{k,l}_{\rho}=0$ for all pairs of integers $k,l\geq 0$.
If $\rho$ is a map of rings it follows the connection $R_{\rho}$ is a flat connection.
\end{lemma}
\begin{proof} Assume $\rho$ is a map of rings. It follows for $D\in \Diff^k(A), D'\in \Diff^l(A)$ we get
\[ R^{k,l}_{\rho}(D,D')=\rho^k(D)\rho^l(D')-\rho^{k+l}(D\circ D')=\rho(D)\rho(D')-\rho(D\circ D')=0.\]
It follows $R^{k,l}_{\rho}=0$ for all $k,l$. The converse is proved in a similar fashion.
If $\rho$ is a map of rings it follows $R^{1,1}_{\rho}=0$ hence by Lemma \ref{curvrho} it follows $R_{\rho}$ is a flat connection.
The Lemma follows.
\end{proof}

Note: A map of rings $\rho: \Diff(A)\rightarrow \Diff(E)$ is sometimes referred to as a \emph{stratification} in the litterature (see \cite{ogus}). Hence the maps $R^{k,l}_{\rho}$ are \emph{obstructions} for  $(E,\rho)$ to
be a stratification. The following may happen: For a given choice of projective basis $x_i,e_j$ for $E$, it might be the corresponding $\infty$-connection $(E,\rho)$ is not a stratification. It might still
be there exists another projective basis $x_i',e_j'$ for $E$ such that the corresponding $\infty$-connection $(E,\rho')$ is a stratification. To determine if a module $E$ has a stratification is a difficult problem
in general. It is well known that if $A$ is a finitely generated regular algebra over a field $k$ of characteristic zero and $E$ is a coherent $A$-module it is neccessary that $E$ is projective for $E$ 
to have a stratification $\rho$. Given a stratification $\rho$ on a coherent module $E$ we get a connection $\nabla:\Der_k(A)\rightarrow \End_k(E)$ and one uses the connection $\nabla$ to prove that
$E$ is locally free, hence projective. The proof does not use the flatness of the connection $\nabla$.

Note: One would like to realize the category of  $l$-connections $(E,\rho^l)$ and morphisms of $l$-connections as a module category over "some universal enveloping algebra" $\U^{ua}(\Diff^l(A))$ 
of the module of $l$'th order differential operators $\Diff^l(A)$ as done for $1$-connections, and to define the notion of cohomology and homology of an $l$-connection $(E,\rho^l)$ as done in Definition \ref{cohomology}
for $1$-connections.  

When $A$ is a finitely generated and regular commutative ring over a field $k$ of characteristic zero, it follows $\Diff(A)$ is generated by $\operatorname{D}^1(A,0)$: Every higher order differential operator $D$ 
is a sum of products of first order differential operators. Hence to give a ring homomorphism $\rho: \Diff(A)\rightarrow \Diff(E)$ is equivalent to give a flat connection 
$\nabla: \operatorname{D}^1(A,0) \rightarrow \Diff^1(E)$. This property does not hold in characteristic $p>0$.

\begin{example} \label{cherntype} Rings of differential operators, annihilator ideals and polynomial relations between Chern classes of a connection. \end{example}

The universal algebras $\Uo(L(\fal))$ and $\Ur(L(\fal))$ may have applications in the theory of characteristic classes. Recall the following results from from \cite{maa1}:

\begin{lemma} \label{lemmaprod} Let $A$ be a commutative ring containing the field of rational numbers $k$ and let  $(\rho, E)$ be an $L$-connection of curvature type $f$ where  $f\in Z^2(L,A)$ is a 2-cocycle
and where $E$ is a projective $A$-module of rank $rk(E)$. 
The following formula holds in $C^{2k}(L,\End_A(E))$:
\[ R_{\rho}^k(x_1,..,x_{2k})=f^k(x_1,..,x_{2k})I_E \]
where $I_E \in \End_A(E)$ is the identity endomorphism of $E$.
\end{lemma}
\begin{proof} See \cite{maa1}, Lemma 5.14. 
\end{proof}

The graded ring $\H^{2*}(L,A):=\oplus_{k=0,\ldots, l} \H^{2k}(L,A)$ is a commutative ring. Let $k[x_1,..,x_l]$ be the polynomial ring in the independent variables $x_1,..,x_l$. Given a polynomial $P(x_1,..,x_l)$ and a 
connection $(\rho, E)$ where $E$ is a projective $A$-module of finite rank we may evaluate the polynomial $P$ in the Chern classes $c_k(E)$ to get a cohomology class
\[ P(c_1(E),c_2(E),..,c_l(E)) \in \H^{2*}(L,A).\]

Given a set of connections $(\rho_i, E_i)_{i=1,\ldots, k}$ where $E_i$ is a finitely generated and projective $A$-module we may for any polynomial $P(x_1,\ldots, x_k)\in \mathbb{Q}[x_1,\ldots, x_k]$
consider the class

\[ P(c_1(E_1),\ldots ,c_1(E_k)) \in \H^{2*}(L,A).\]

When we vary the polynomial $P$ and the number $k$ of independent variables, 
we get a subring $R(c_1)$ of the cohomology ring $\H^{2*}(L,A)$ - the subring (over the field of rational numbers) generated by the first Chern classes $c_1(E) \in H^2(L,A)$ of
all finitely generated and projective $A$-modules $E$.

Define for $k\geq 2$ the polynomial 

\[ P_k(x_1,..,x_l):=rk(E)^{k-1}x_k-x_1^k \in k[x_1,..,x_l]. \]

Lemma \ref{lemmaprod} has consequences for the Chern class $c_k(E)$.

\begin{corollary} \label{det} Let $A$ be a commutative ring containing the field of rational numbers $k$ and let  $(\rho, E)$ be an $L$-connection of curvature type $f$  where $f\in Z^2(L,A)$ is a 2-cocycle
and where $E$ is a projective $A$-module of rank $rk(E)$. 
The following holds: 
\[ P_k(c_1(E),..,c_l(E))=0\]
for all $k\geq 2$.
\end{corollary}
\begin{proof} By Lemma \ref{lemmaprod} we get
\[ c_k(E) = tr(R_{\rho}^k)=tr(f^k I_E)=rk(E)f^k \in \H^{2k}(L,A).\]
By definition
\[ c_1(E)=tr(R_{\rho})=rk(E) f \in \H^2(L,A).\]
It follows
\[ c_1(E)^k=rk(E)^k f^k \]
hence
\[   P_k(c_1(E),..,c_l(E)):=rk(E)^{k-1}c_k(E) - c_1(E)^k =0 \]
and the Corollary follows.
\end{proof}

Hence from Corollary \ref{det} the following holds: For a connection $(\rho,E)$ where $\rho$ has curvature type $f$ for a 2-cocycle $f\in\Z^2(L,A)$ we get an equality
\[ c_k(E)=\frac{1}{rk(E)^{k-1}}c_1(E)^k .\]
Hence the $k$'th Chern class $c_k(E)$ is determined by the first Chern class and hence $c_k(E)\in R(c_1)$ for all $k \geq 1$.

Hence the annihilator ideal $ann(\rho, E)$ can be used to detect if there are polynomial relations between the Chern classes of $E$. If the connection $(\rho, E)$ has curvature type $f\in \Z^2(L,A)$, this 
puts strong conditions on the Chern classes of $E$: Hence the annihilator ideal $ann(\rho, E)$ detects if the Chern class $c_k(E)$ is interesting from the point of view of
Hodge theory. The Hodge conjecture is known for $\H^2$ hence the ring $R(c_1)$ is well known when $\H^2(L,A)$ calculates the 2'nd singular cohomology of $X_{\C}$ with complex coefficients.

\begin{example} 2-sided ideals in the universal ring $\Uo(L(\fal))$ and Morita equivalence. \end{example}

Let $(\tL, \ta,\tp,[,],D)$ be a $\Dl$-Lie algebra and let $\rho:\tL\rightarrow \End_k(E)$ be an object in $\Mod(\tL,Id)$. 
%Let $B:=A\otimes_k A$. 
%\[ \Hom_{B}(\tL,\End_k(E)) \cong \Hom_{B-alg}(\T^*_B(\tL), \End_k(E)) \]
%where $\T^*_B(\tL)$ is the tensor algebra of $\tL$ as $B$-module. Here $\End_k(E)$ has the canonical $B$-module structure. Hence 
There is an equivalence of categories

\[ \Mod(\tL, Id) \cong \Mod(\Uo(\tL))).\]
Let in the following $f\in \Z^2(\Der_k(A),A)$ and let $(L,\alpha)$ be an $A/k$-Lie-Rinehart algebra. Let $(L(\fal), \alpha_f, \pi_f,[,],z)$ be the corresponding $\Dl$-Lie algebra.
Assume $(\rho,E)$ is an object in $\Mod(L(\fal), Id)$ with the property that $R_{\rho}(u,v)=0$ for all $u,v\in L(\fal)$. It follows $(\rho,E)$ is an object in $\Mod(\tUo(L(\fal)))$ since
$\tUo(L(\fal))$ is the quotient of $\Uo(L(\fal))$ by the 2-sided ideal generated by elements on the form $u\otimes v-v\otimes u-[u,v]$ for $u,v\in L(\fal)$. The category
$\Mod(L(\fal), Id)$ is equivalent to the category of $L$-connection of curvature type $\fal$ hence there is an equivalence of categories
\[ \Mod(\tUo(L(\fal)) \cong \Mod(\U(A,L,\fal)) .\]
Hence the two associative rings $\tUo(L(\fal))$ and $\U(A,L,\fal)$ are Morita equivalent. They are not isomorphic in general but they both have isomorphic centres equal to the base ring $k$.

\begin{example} \label{familyconn} Families of 2-sided ideals in $\Ur(L(0))$ and families of connections\end{example}

Let Consider the abelian extension $L(0)$ of $L$ by $Az$ and the canonical map $\ta: L(0)\rightarrow \Der_k(A)$ defined by $\ta(az+x):=\alpha(x)$. It follows $(L(0),\ta)$ is an $A/k$-Lie-Rinehart
algebra. There is an equivalence of categories
\[  \Conn(L(0),Id) \cong \Mod(\Ur(L(0)))\]
where $\Conn(L(0),Id)$ is the category of connections $\tilde{\rho}:L(0)\rightarrow \End_k(E)$ such that $\tilde{\rho}(z)=Id_E$. The category $\Conn(L(0),Id)$ is equivalent to the category
$\Conn(L)$ of ordinary connections  $\rho:L\rightarrow \End_k(E)$ and morphisms. We get an equivalence of categories
\[ \Conn(L) \cong \Mod(\Ur(L(0))) .\]
For any 2-cocycle $f\in \Z^2(L,A)$ there is an associated 2-cocycle $\tilde{f} \in \Z^2(L(0),A)$ and  we may consider the 2-sided ideal
\[  I(f):=\{u\otimes v -v\otimes u-[u,v]-\tilde{f}(u,v)z \text{ such that $u,v\in L(0)$.} \} \subseteq \Ur(L(0)) .\]
Here $\tilde{f}(u,v):=f(x,y)$.

Let $\Ur_f(L(0)):=\Ur(L(0))/I(f)$ be the quotient. There is an equivalence of categories between the category $\Mod(\Ur_f(L(0)))$ and the category of connections $(\tilde{\rho},E)$ in $\Mod(L(0),Id)$ such that
\[  \tilde{\rho}(u)\tilde{\rho}(v)-\tilde{\rho}(v)\tilde{\rho}(u)-\tilde{\rho}([u,v])-\tilde{f}(u,v)Id_E=0 \text{ in $\End_k(E)$}.\]
It follows 
\[ R_{\tilde{\rho}}(u,v)=\tilde{f}(u,v)Id_E \]
hence the induced connection $\rho:= \tilde{\rho} \circ i$ on $L$ has curvature type $f$. It follows there is an equivalence of categories 
\[  \Mod(\Ur_f(L(0))) \cong \Mod(\U(A,L,f)) \]
hence $\Ur_f(L(0))$ and $\U(A,L,f)$ are Morita equivalent rings. They are not isomorphic in general but have $k$ as centre.
Hence we may construct the categories $\Mod(\U(A,L,f))$ using quotients of one fixed ring $\Ur(L(0))$ by the family of 2-sided ideals
$I(f)$ for $f\in \Z^2(L,A)$. In the case when $A$ is noetherian and $L$ a finitely generated an projective $A$-module it follows $\U(A,L,f)$ is Noetherian. Since being Noetherian
is Morita invariant it follows the ring $\Ur_f(L(0))$ is Noetherian for any 2-cocycle $f\in \Z^2(L,A)$.

We observe that an ordinary $L$-connection $(\rho,E)$ gives rise to a connection $\tilde{\rho}:L(\fal)\rightarrow \End_k(E)$. The map $\tilde{\rho}$ is an $A\otimes_k A$-linear map.
It follows $(\rho,E)$ is a left $\Uo(L(\fal))$-module.
If the annihilator ideal $ann(\rho, E)$ contains the ideal $J(\fal)$ generated by elements on the form $u\otimes v-v\otimes u-[u,v]$ for $u,v\in L(\fal)$ 
it follows the connection $\rho$ has curvature type $f^{\alpha}$. It follows from Corollary \ref{det} 
the $k$-th Chern class $c_k(E)$ is determined by $c_1(E)$. Hence the structure of the set of 2-sided ideals in $\Uo(L(\fal))$ can be used to study properties of the Chern classes $c_k(E)$. 

\begin{example} Moduli spaces of $\Gamma$-modules and connections\end{example}

When studying moduli spaces of connections many authors use the Hilbert scheme
and Quot scheme to construct parameter spaces of connections and these spaces are large and complicated. The rings $\Uo(\tL)$ and $\Ur(\tL)$ are non-Noetherian in general but they have as
shown in this paper many Noetherian quotients. It could be one gets an alternative to the study of the Hilbert and Quot schemes by studying parameter spaces of 2-sided ideals in Noetherian
quotients of the rings $\Uo(\tL)$ and $\Ur(\tL)$. Instead of studying large parameter spaces of pairs $(\rho,E)$ where $\rho$ is an $\tL$-connection, 
we study the parameter space of annihilator ideals $ann(\rho,E)$ in a Noetherian quotient one of the rings $\Uo(\tL)$ and $\Ur(\tL)$. 

In \cite{simpson} the author constructs for any smooth projective
complex variety $X$, any sheaf of filtered algebras $\Gamma$ on $X$ and any nummerical polynomial $P$, a quasi projective scheme $\M(X, \Gamma, P)$ parametrizing  
semi stable $\Gamma$-modules with Hilbert polynomial $P$. The construction of $\M(X, \Gamma,P)$ uses the Hilbert scheme and GIT quotients. 
Assume we are given a parameter space $\M(d,\Gamma,P)$ parametrizing locally trivial $\Gamma$-modules $(\rho, \mathcal{E})$ with Hilbert polynomial $P$, such that $\mathcal{E}$ is a locally trivial $ \O_X$-modules
of rank $d$. Given two isomorphic $\Gamma$-modules $(\rho,\mathcal{E})$ and  $(\rho',\mathcal{E}')$ where $\mathcal{E}$ and $\mathcal{E}'$ are isomorphic locally free $\O_X$-modules corresponding to different points in the parameter space $\M(d,\Gamma,P)$. The sheaves of annihilator ideals $ann(\rho, \mathcal{E})$ and $ann(\rho', \mathcal{E}')$ in $\Gamma$ will be equal. Hence in the parameter space
$\M(d,\Gamma,P)$ we get two different points corresponding to $(\rho,\mathcal{E})$ and $\rho',\mathcal{E}')$. In the parameter space of sheaves of 
2-sided ideals we get one point corresponding to $ann(\rho,\mathcal{E})=ann(\rho', \mathcal{E}')$. Hence we should expect the parameter space of sheaves of 2-sided ideals in $\Gamma$ to have fewer points
than the parameter space $\M(d,\Gamma,P)$.  In the case of a holomorphic Lie algebroid $\mathcal{L}$ on a complex projective manifold $X$, 
it follows from \cite{tortella} that the moduli spaces $\M_{\mathcal{L},Q}(P)$ are in many cases empty. Hence one should take care when studying such moduli spaces in general: If one is unable to write down explicit
non-trivial examples, this may indicate they are empty. In the affine situation as shown in this paper, it is relatively easy to write down explicit non-trivial examples of the theory as shown in Theorem \ref{curvform}.
One may moreover implement computer algorithms calculating such examples.

%If $ann(E,\nabla) \subseteq \T^*_B(L(\fal))$ is a maximal ideal it follows $(E,\nabla)$ is a simple module. One wants to classify maximal 2-sided ideals in $\T^*_B(L(\fal))$ arising as annihilator ideals of connections
 %$(E,\nabla)$ where $E$ is a finitely generated and projective as left $A$-module.

%Hence the annihilator ideal $ann(E)=ann(\rho, E)$ in $\Uo(L(\fal))$ of the connection $(\rho, E)$ is related to 
%properties of the Chern classes $c_k(E)$ of $(\rho, E)$ and the Chern character $Ch: \K_0(L) \rightarrow \H^{2*}(L,A)$.

If $A$ is finitely generated and regular algebra over the complex numbers $\C$ and $L:=\Der_{\C}(A)$ it follows $\H^i(L,A)\cong \H^i_{sing}(X_{\C}, \C)$ is singular cohomology of $X_{\C}$ with complex coefficients. 
Here $X_{\C}$ is the underlying complex manifold of $X:=\Spec(A)$.  The Hodge conjecture is known to hold for $\H^2_{sing}(X_{\C}, \C)$ hence if one wants to study cohomology classes in $\H^{2*}_{sing}(X_{\C} ,\C)$  that are not coming from $\H^2_{sing}(X_{\C}, \C)$ under the cup product, one needs to study other types of connections.  One has to study left  $\Uo(L(\fal))$  modules $(\rho,E)$ that are finitely generated and projective over $A$, such that $ann(\rho, E)$ does not contain the ideal $J(\fal)$ for a 2-cocycle $\fal$ - we get a computable criteria on $(\rho, E)$ which can be used to determine if the Chern class $c_k(E)$ is "interesting". Hence the structure of the set of 2-sided ideals in  $\Uo(L(\fal))$ is related to the study of algebraicity of cohomology classes in singular cohomology.  

For the universal enveloping algebra $\U(\lg)$ of a finite dimensional semi simple Lie algebra $\lg$ it follows the set of 2-sided ideals in $\U(\lg)$ corresponding to finite dimensional irreducible $\lg$-modules is a discrete set. 
The enveloping algebras $\U(\lg)$ and $\U(A,L,f)$  are Morita equivalent to quotients of the universal algebra $\Uo(L(\fal))$. One wants to study the correspondence between 2-sided ideals in 
$\Uo(L(\fal))$ and cohomology classes in $\H^{2*}(L,A)$ and use this correspondence to give interesting examples of algebraic and non-algebraic classes in the singular cohomology of a complex algebraic manifold.

\section{The universal ring is an almost commutative Noetherian ring}

In this section we construct in Theorem \ref{mainuniversal} for any $\Dl$-Lie algebra  $(\tL, \ta, \tp,[,],D)$ and any $\tL$-connection $(\rho,E)$ the universal ring $\tUo(\tL,\rho)$ of $(\rho,E)$. In the case 
when $A$ is Noetherian and $\tL,E$ finitely generated as left $A$-modules it follows the associative unital ring $\tUo(\tL, \rho)$ is an almost commutative Noetherian sub ring of $\Diff(E)$ - the ring of differential
operators  on $E$. We prove a similar result for the ring $\Uo_E(\tL)$ - it is an almost commutative sub-ring of $\tUo(\tL, \rho)$.
The non-flat connection $(\rho,E)$ is a finitely generated left $\tUo(\tL,\rho)$-module and we may use $\tUo(\tL, \rho)$ to construct the characteristic variety $\SS(\rho,E)$ of $(\rho,E)$.
We may use the variety $\SS(\rho,E)$ to define holonomicity for non-flat connections. Previously this was defined for flat connections (see Example \ref{characteristic}).

%In this section we construct for a connection $\rho:\tL\rightarrow \End_k(E)$ where $(\tL,\ta;\tp;[,],D)$ is a $\Dl$-Lie algebra and $E$ is a finitely generated  $A$-module 
%the universal ring $\U(\tL,\rho)$ of the connection $(E,\rho)$.
%The universal ring $\U(\tL,\rho)\subseteq \Diff(E)$ is in the case when $L$ is a finitely generated $A$-module an almost commutative Noetherian sub ring of $\Diff(E)$ - 
%the ring of differential operators on $E$. The ring $\U(\tL,\rho)$ is constructed as a quotient
%of the tensor algebra $\T^*_B(\End(\tL,E))$ of the $\Dl$-Lie algebra $\End(\tL,E)$ constructed in \cite{maa141}.  
%The subring $\U(\tL,\rho) \subseteq \Diff(E)$ gives many  examples of Noetherian almost commutative quotients
%of the tensor algebra  $\T^*_B(\End(\tL,E))$. The ring $\T^*_B(\End(\tL,E))$ is not Noetherian in general. The connection $\rho$ is non-flat in general.
%In the classical situation for a flat connection $(E,\nabla)$ on a Lie-Rinehart algebra $(L,\alpha)$ we may construct the universal enveloping algebra 
%$\U(A,L)$ of Rinehart (see \cite{rinehart}) and one uses the PBW-theorem to prove $\U(A,L)$ is Noetherian when $A$ is Noetherian and $L$ is a finitely generated and projective $A$-module.

\begin{proposition} Let $f \in \Z^2(\Der_k(A),A)$ be a 2-cocycle and let $(L,\alpha)$ be an $A/k$-Lie-Rinehart algebra. Let $\alpha_f:L(\fal)\rightarrow \D$ be the corresponding
 $\Dl$-Lie algebra. Let $\rho: L(\fal)\rightarrow \End_k(E)$ be an object in $\Mod(L(\fal))$. It follows $\rho$ is $A\otimes_k A$-linear map. 
We get an induced map
\[ \T(\rho): \Uo(L(\fal))\rightarrow \Diff(E) \]
and
\[ \T(\rho)_i: \Uo_i(L(\fal)) \rightarrow \Diff^i(E) \]
for all $i \geq 1$. 
\end{proposition}
\begin{proof} Let $u:=az+x,v:=bz+y \in L(\fal)$ and let $c\in A$. Since $\rho$ is $A\otimes_k A$-linear we get the following: 
\[ \rho(cu)=c\rho(u)\text{ and }\rho(u+v)=\rho(u)+\rho(v).\]
Moreover
\[ \rho(uc)=\rho(u)c.\]
We get for $e\in E$ the following:
\[ \rho(u)(ce)=\rho(uc)(e)=\rho(cu+\alpha_f(u)(c)z)(e)=c\rho(u)(e)+\alpha_f(u)(c)\rho(z)(e).\]
We get
\[ \rho(z)(ce)=\rho(zc)(e)=\rho(cz)(e)=c\rho(z)(e) \]
 hence  $\rho(z)\in \End_A(E).$
It follows
\[ [\rho(u),cId_E]=\rho(u)c-c\rho(u)=\alpha_f(u)(c)Id_E\in \End_A(E) \]
hence
\[ \rho(u)\in \Diff^1(E). \]
It follows we get an induced map
\[ \T(\rho): \Uo(L(\fal)) \rightarrow \Diff(E) \]
and
\[ \T(\rho)_i :\Uo_i(L(\fal)) \rightarrow \Diff^i(E) \]
and the Lemma follows.
\end{proof}

%\begin{lemma} Let $k\rightarrow A$ be a map of commutative unital rings. Let $\Diff_k(E):=\Diff(E)$ be the ring of k-linear differential operators on $E$. The ring $\Diff(E)$ is almost commutative.
%\end{lemma}
%\begin{proof} We mus prove that for any differential operators $x\in \Diff^i(E)$ and $y\in \Diff^j(E)$ it follows $xy-yx \in \Diff^{i+j-1}(E)$. By definition it follows for any element $a \in A$ that
%\[ [x,aId_E]\in \Diff^{i-1}(E) \]
%and
%\[ [y,aId_E]\in \Diff^{j-1}(E).\]
%It follows we may write $xa=ax+ \eta$ with $\eta \in \Diff^{i-1}(E)$ and $ya=ay+\eta'$ with $\eta'\in \Diff^{j-1}(E)$. We get
%\[  [xy-yx, aId_E]=xya-yxa-a(xy-yx)=\]
%\[ x(ay+\eta')-y(ax+\eta)-a(xy-yx)=\]
%\[ a(xy-yx)+\eta y +x\eta'-\eta'x -y\eta-a(xy-yx)=\eta y +x\eta'-\eta'x-y\eta \in \Diff^{i+j-1}(E).\]
%The Lemma follows.
%\end{proof}

\begin{lemma}Let $f\in \Z(\Der_k(A),A)$ be a 2-cocycle, $(L,\alpha)$ an $A/k$-Lie-Rinehart algebra and let $\alpha_f:L(\fal)\rightarrow \D$ be the associated 
$\Dl$-Lie algebra. There is a one-to-one correspondence between
the set of $A\otimes_k A$-linear maps $\rho:L(\fal)\rightarrow \End_k(E)$ and the set of $\psi$-connections $\nabla:L\rightarrow \End_k(E)$ for varying $\psi \in \End_A(E)$.
If $\nabla:L\rightarrow \End_k(E)$ is a $\psi$-connection it follows the corresponding $A\otimes_k A$-linear map map $\rho:L(\fal) \rightarrow \End_k(E)$ has curvature
\[ R_{\rho}(u,v)=R_{\nabla}(x,y)-f(x,y)Id_E +b[\nabla(x),\psi] -a[\nabla(y),\psi] \]
with $u=az+x, v=bz+y\in L(\fal)$. If $\psi=Id_E$ it follows 
\[ R_{\rho}(u,v)=R_{\nabla}(x,y)-f(x,y)Id_E.\]
Conversely, if $\rho:L(\fal)\rightarrow \End_k(E)$ is an $A\otimes_k A$-linear map and $\nabla$ the corresponding $\psi$-connection, it follows the curvature of $\nabla$ satisfies
\[ R_{\nabla}(x,y)=R_{\rho}(i(x),i(y)) +f(x,y)\psi \]
where $i:L\rightarrow L(\fal)$ is the canonical inclusion map.
\end{lemma}
\begin{proof} The first statement is proved earlier in this paper. Let $\nabla:L\rightarrow \End_k(E)$ be a $\psi$-connection and let $\rho:L(\fal)\rightarrow \End_k(E)$ be the corresponding
$A\otimes_k A$-linear map. By definition $\rho(az+x):=a\psi+ \nabla(x)$. We get for any elements $u:=az+x,v:=bz+y\in L(\fal)$ the following calculation:
\[ R_{\rho}(u,v):= [\rho(u),\rho(v)]-\rho([u,v])=\]
\[[a\psi +\nabla(x),b\psi+\nabla(y)]-(\alpha(x)(b)-\alpha(y)(a)+f(x,y))z +\nabla([x,y])=\]
\[R_{\nabla}(x,y)+b[\nabla(x),\psi]-a[\nabla(y),\psi]-f(x,y)\psi\]
and the second claim follows. The third claim follows since $[\nabla(x),Id_E]=[\nabla(y),Id_E]=0$. Let $i:L\rightarrow L(\fal)$ be the left $A$-linear canonical inclusion map. It follows for any elements 
$x,y\in L$ we get
\[ [i(x),i(y)]-i([x,y])=f(x,y)z \in L(\fal).\] We get
\[ R_{\nabla}(x,y):= [\rho(i(x)),\rho(i(y))]-\rho(i([x,y])).\]
We get
\[ i([x,y])=[i(x), i(y)]-f(x,y)z \]
hence
\[ R_{\nabla}(x,y)=[\rho(i(x)),\rho(i(y))]-\rho([i(x),i(y)])+\rho(f(x,y)z)=\]
\[R_{\rho}(i(x),i(y))f(x,y)\rho(z)=R_{\rho}(i(x),i(y)) +f(x,y)\psi.\]
The Lemma is proved.
\end{proof}

\begin{lemma} Let $\nabla:L\rightarrow \End_k(E)$ be a $\psi$-connection where $\psi\in \End_A(E)$. It follows
\[ R_{\nabla}(x,y)(ae)=aR_{\nabla}(x,y)(e) + \alpha(y)(a)[\nabla(x),\psi](e)-\alpha(x)(b)[\nabla(y),\psi](e) \]
for all $a\in A$ and $e\in E$.
\end{lemma}
\begin{proof} The proof is a straightforward calculation.
\end{proof}

\begin{lemma} Let $f\in \Z^2(\Der_k(A),A)$ and let $\nabla:L\rightarrow \End_k(E)$ be a $\psi$-connection with $\psi\in \End_A(E)$. Let $\rho:L(\fal)\rightarrow \End_k(E)$ 
be the associated $A\otimes_k A$-linear map. It follows the curvature $R_{\rho}(u,v)\in \Diff^1(E)$. If $\psi=Id_E$ it follows $R_{\rho}(u,v)\in \Diff^0(E):=\End_A(E)$.
\end{lemma}
\begin{proof} We have seen that if $u:=az+x,v:=bz+y$ it follows
\[ R_{\rho}(u,v)=R_{\nabla}(x,y)-f(x,y)\psi+b[\nabla(x),\psi]-a[\nabla(y),\psi].\] 
The following holds for $R_{\nabla}$.
\[ R_{\nabla}(x,y)(ae)=aR_{\nabla}(x,y)(e)+\alpha(y)(a)[\nabla(x),\psi](e)-\alpha(x)(b)[\nabla(y),\psi](e) \]
for all $a\in A, e\in E$. One checks that
\[ f(x,y)\psi+b[\nabla(x),\psi]-a[\nabla(y),\psi] \in \End_A(E)\]
for all $u,v$.
We get for any $a\in A$ the following:
\[ [R_{\nabla}(x,y),aId_E]= R_{\nabla}(x,y)a - aR-{\nabla}(x,y)=\]
\[\alpha(y)(a)[\nabla(x),\psi]-\alpha(x)(a)[\nabla(y),\psi] \in \End_A(E)\]
hence
\[ R_{\nabla}(x,y)\in \Diff^1(E).\]
It follows
\[ R_{\rho}(u,v)\in \Diff^1(E).\]
One checks that if $\psi=Id_E$ it follows $R_{\rho}(u,v)\in \End_A(E)$ and the Lemma follows.
\end{proof}

\begin{lemma} Let $(\tL,\ta,\tp,[,],D)$ be a $\Dl$-Lie algebra and let $\rho:\tL \rightarrow \End_k(E)$ be an $A\otimes_k A$-linear map. Let $u,v\in \tL$. The following holds:
\[\rho: \tL \rightarrow \Diff^1(E).\]
For any set of elements $u_1,..,u_i\in \tL$ it follows $\rho(u_1)\circ \cdots \circ \rho(u_i)\in \Diff^i(E)$.
Moreover
\[ R_{\rho}(u,v) \in \Diff^2(E) .\]
If $\rho(D)=Id_E$ it follows
\[ R_{\rho}(u,v)\in \Diff^1(E).\]
\end{lemma}
\begin{proof} Let $a\in A$. It follows
\[ [\rho(u),\phi_a]=\tp(u)(a)\rho(D)\in \Diff^0(A):=\End_A(E) \]
hence $\rho(u)\in \Diff^1(A)$. Since $\Diff(E)$ is a filtered associative ring it follows for any set of $i$ elements $u_1,..,u_i \in \tL$ we get
\[ \rho(u_1)\circ \cdots \circ \rho(u_i)\in \Diff^i(E).  \]
We get
\[ [\R_{\rho}(u,v), \phi_a]=\tp(v)(a)[\rho(u),\rho(D)]+\tp(u)(a)[\rho(D),\rho(v)]+ \]
\[ [\tp(u),\tp(v)](a)\rho(D)\circ \rho(D)-\tp([u,v])(a)\rho(D) \in \Diff^1(E) \]
hence $R_{\rho}(u,v)\in \Diff^2(E)$ in general. If $\rho(D)=Id_E$ we get
\[  [\R_{\rho}(u,v), \phi_a]=  [\tp(u),\tp(v)](a)\rho(D)\circ \rho(D)-\tp([u,v])(a)\rho(D) \in \Diff^0(E) \]
hence $R_{\rho}(u,v)\in \Diff^1(E)$.
The Lemma follows.
\end{proof}

\begin{example} The universal ring of an $\tL$-connection.\end{example}

Let $(\tL,\ta,\tp,[,],D)$ be a $\Dl$-Lie algebra and let $\rho:\tL \rightarrow \End_k(E)$ be an object in  $\Mod(\tL, Id)$. This means
$\rho$ is $A\otimes_k A$-linear and $\rho(D)=Id_E$. 

\begin{lemma} The connection $\rho$ induce a connection
\[ \tilde{\rho}:\tL \rightarrow \Der_k(\End_A(E)) \]
defined by
\[ \tilde{\rho}(u)(\phi):=[\rho(u),\phi] \]
with curvature
\[ R_{\tilde{\rho}}(u,v)(\phi)=[R_{\rho}(u,v),\phi] . \]
\end{lemma}
\begin{proof} Let $a\in A, e\in E, u\in \tL$ and $\phi \in \End_A(E)$. We get
\[ \tilde{\rho}(u)(\phi)(ae)=\rho(u)(\phi(ae))-\phi(\rho(u)(ae)=\]
\[a\rho(u)(\phi(e))+\tp(u)(a)\rho(D)(\phi(e))-a\phi(\rho(u)(e))-\tp(u)(a)\phi(\rho(D)(e))=a[\rho(u), \phi](e)=\]
\[ a\tilde{\rho}(u)(\phi)(e) \]
hence $\tilde{\rho}(\phi)\in \End_A(E)$. One checks that for any two $\phi, \psi\in \End_A(E)$ we get
\[ \tilde{\rho}(u)([\phi,\psi])=[\tilde{\rho}(u)(\phi),\psi]+[\phi, \tilde{\rho}(u)(\psi)] .\]
One checks that $R_{\tilde{\rho}}(u,v)(\phi)=[R_{\rho}(u,v),\phi]$ and the Lemma follows.
\end{proof}

We get by \cite{maa141} a non-abelian extension of $A/k$-Lie-Rinehart algebras
\begin{align}
&\label{nonabelian} 0 \rightarrow \End_A(E) \rightarrow \End(\tL, E) \rightarrow \tL \rightarrow 0 
\end{align}
where $\End(\tL,E):=\End_A(E) \oplus \tL $ with the following Lie product. For any $z:=(\phi,u),z':=(\psi,v)\in \End(\tL,E)$ define
\[ [z,z']:=([\phi,\psi]+[\rho(u),\psi]-[\rho(v),\phi]+R_{\rho}(u,v),[u,v]).\]
Define the central element $\tilde{D}:=(0,D)\in \End(\tL,E)$. There is a canonical map 
\[ \alpha_E:\End(\tL, E) \rightarrow \D \]
defined by
\[ \alpha_E(\phi,u):= \ta(u).\]
There is a map
\[ \pi_E: \End(\tL,E)\rightarrow \Der_k(A) \]
defined by
\[ \pi_E(\phi,u):=\tp(u).\]

\begin{proposition} \label{nonabelianconn} Let $(\tL, \ta,\tp,[,],D)$ be a $\Dl$-Lie algebra and let $(\rho,E)$ be an object in $\Mod(\tL, Id)$.
It follows the 5-tuple $(\End(\tL, E),\alpha_E,\pi_E,[,], \tilde{D})$ constructed above is a $\Dl$-Lie algebra.
\end{proposition}
\begin{proof} Let $u\in \tL$ and $\phi\in \End_A(E)$. We get
\[ [\rho(u) ,\phi]a = \rho(u)\circ \phi a-\phi \circ \rho(u)a=\]
\[ a\rho(u) \circ \phi +\tp(u)(a) \phi - a\phi \circ \rho(u)- \phi \tp(u)(a) =\]
\[ a[\rho(u),\phi] \]
hence $\tilde{\rho}(u)(\phi):=[\rho(u),\phi]\in \End_A(E)$. We get a map
\[ \tilde{\rho}:\tL \rightarrow \End_k(\End_A(E)) \]
and one checks that this map induce a map
\[ \tilde{\rho}: \tL \rightarrow \Der_k(\End_A(E)) .\]
One moreover checks that 
\[ R_{\tilde{\rho}}(u,v)(\phi)=[R_{\rho}(u,v),\phi] .\] Hence we get by \cite{maa141}  a non-abelian extension of $A/k$-Lie-Rinehart algebras
\[ 0 \rightarrow \End_A(E) \rightarrow \End(\tL, E) \rightarrow \tL \rightarrow 0 \]
where $\End(\tL, E):=\End_A(E)\oplus \tL$ with the given Lie product. The maps are the canonical maps. Define $\tilde{D}:=(0,D)\in \End(\tL, E)$. Let $z:=(\phi,u)\in \End(\tL,E)$ and let $c\in A$.
Define
\[ zc=(\phi,u)c:=(\phi c, uc)=(c\phi, cu+\tp(u)(c)D)=c(\phi,u)+\tp(u)(0,D).\]
It follows
\[   zc= cz+\pi_E(z)(c)\tilde{D} \]
where we have defined
\[ \pi_E:\End(\tL, E) \rightarrow \D \]
by
\[ \pi_E(\phi,u):=\tp(u)\in \D.\]
It follows $\End(\tL, E)$ is a left and right $A$-module and the map $f: \End(\tL, E)\rightarrow \tL$ is a map of $A\otimes_k A$-modules with $f(\tilde{D})=D$.
Define
\[ \alpha_E:\End(\tL,E)\rightarrow \Der_k(A)\]
by
\[ \alpha_E(\phi,u):= \ta(u)\in \Der_k(A).\]
It follows $\alpha_E$ is a map of $A\otimes_k A$-modules and $k$-Lie algebras. The element $\tilde{D}$ is in the center of $\End(\tL,E)$ hence the 5-tuple
$(\End(\tL,E),\alpha_E,\pi_E,[,],\tilde{D})$ is a $\Dl$-Lie algebra. The sequence \ref{nonabelian} is a non-abelian extension of $\tL$ by $\End_A(E)$ and the Lemma follows.
\end{proof}

\begin{theorem} \label{mainuniversal} Let $(\tL, \ta,\tp,[,],D)$ be a $\Dl$-Lie algebra and let $(\rho, E)$ be an $\tL$-connection with $\rho(D)=Id_E$. There is a canonical map
\[ \rho^!:\End(\tL, E)\rightarrow \Diff^1(E) \]
and $\rho^!$ is a map of $B:=A\otimes_k A$-modules and $k$-Lie algebras. The map $\rho^!$ induce a map $\T(\rho^!):\tUo(\End(\tL,E))\rightarrow \Diff(E)$ of associative rings.
Let $\tUo(\tL,\rho):= Im(\T(\rho^!))$ be the image. We get an exact sequence of rings
\[ 0 \rightarrow ker(\T(\rho^!)) \rightarrow \tUo(\End(\tL, E)) \rightarrow \tUo(\tL, \rho) \rightarrow 0\]
where $\tUo(\End(\tL,E)):= \Uo(\End(\tL, E))/I$ where $I$ is the 2-sided ideal generated by the elements $u \otimes v-v\otimes u - [u,v]$ for $u,v\in \End(\tL,E)$. The rings
$\tUo(\End(\tL,E))$ and $\tUo(\tL,\rho)$ are almost commutative. If $A$ is noetherian and $\tL,E$ finitely generated as left $A$-modules it follows
$\tUo(\End(\tL,E))$ and $\tUo(\tL,\rho)$ are Noetherian rings.
\end{theorem}
\begin{proof} By Proposition \ref{nonabelianconn} we may do the following: Define the map 
\[ \rho^!:\End(\tL,E)\rightarrow \Diff^1(E) \]
by
\[ \rho^!(\phi,u):=\phi+\rho(u).\]
Since $\rho(u)\in \Diff^1(E)$ it follows $\rho^!(\phi,u)\in \Diff^1(E).$ By definition it follows $\rho^!$ is an $A\otimes_k A$-linear map. Let $z:=(\phi,u),z':=(\psi,v)\in \End(\tL,E)$. We get
\[ \rho^!([z,z'])=\rho^!([\phi,\psi]+[\rho(u),\psi]+[\phi,\rho(v)]+R_{\rho}(u,v)+\rho([u,v])=\]
\[ [\phi,\psi]+[\rho(u),\psi]+[\phi,\rho(v)]+[\rho(u),\rho(v)]= \]
\[  [\phi+\rho(u),\psi+\rho(v)]=[\rho^!(z),\rho^!(z')] \]
and the map $\rho^!$ is a map of $k$-Lie algebras. We get an induced map
\[ \T(\rho^!):\Uo(\End(\tL,E))\rightarrow \Diff(E) \]
Since $\T(\rho^!)(u\otimes v-v\otimes u-[u,v]))=0 $ it follows we get an induced exact sequence
\[ 0 \rightarrow ker(\T(\rho^!)) \rightarrow \tUo(\End(\tL,E)) \rightarrow \tUo(\tL, \rho) \rightarrow 0\]
By Lemma \ref{almostcommquot} and Proposition \ref{almostcomm} it follows $\tUo(\End(\tL,E))$ and $\tUo(\tL, \rho)$ are almost commutative associative unital rings. 
If $A$ is Noetherian and $\tL,E$ are finitely generated as left $A$-modules it follows from Theorem \ref{mainnoetherian} 
$Gr(\tUo(\End(\tL,E)))$, $Gr(\tUo(\tL,\rho))$, $\tUo(\End(\tL,E))$ and $\tUo(\tL, \rho)$ are Noetherian rings . The Theorem follows.
\end{proof}

\begin{definition} Let $\tUo(\tL, \rho):=Im(\T(\rho^!))\subseteq \Diff(E)$ be the \emph{universal ring} of the connection $(\rho, E)$.
\end{definition}

Recall the ring $\Uo_E(\tL)$ from Definition \ref{annihilator}.

\begin{corollary} \label{smainalmost}Let $(\tL, \ta,\tp,[,],D)$ be a $\Dl$-Lie algebra and let $\rho: \tL \rightarrow \End_k(E)$ be an object in $\Mod(\tL, Id)$. It follows the ring $\Uo_E(\tL)$ is an almost commutative subring
of $\Diff(E)$. If $A$ is Noetherian and $E$ a finitely generated $A$-module it follows $\Uo_E(\tL)$ is an almost commutative and Noetherian sub ring of $\Diff(E)$.
\end{corollary}
\begin{proof} By definition it follows the ring $\tUo(\tL, \rho)$ is generated as a sub ring of $\Diff(E)$ by elements on the form
\[ (\phi_1+\rho(u_1))\circ \cdots \circ (\phi_i+\rho(u_i)) \]
with $\phi_j\in \End_A(E)$ and $u_j\in \tL$.
There is an exact sequence
\[ 0 \rightarrow ann^{\otimes}(E,\rho) \rightarrow \Uo(\tL) \rightarrow \Uo_E(\tL) \rightarrow 0 \]
and $\Uo_E(\tL) \subseteq \Diff(E) $ is the subring generated by elements on the form
\[ \rho(u_1) \circ \cdots \circ \rho(u_i) \]
for $u_j\in \tL$. Since the element $\rho(u_1) \circ \cdots \circ \rho(u_i) \in \tUo(\tL, \rho)$ 
it follows $\Uo_E(\tL)\subseteq \tUo(\tL, \rho)$ is a sub ring. Hence by Lemma \ref{almostcomm} and Theorem \ref{mainuniversal} 
 it follows $\Uo_E(\tL)$ is an almost commutative sub-ring of $\tUo(\tL, \rho)$ and $\Diff(E)$. If $A$ is Noetherian and $E$ finitely generated as $A$-module, it follows $\Uo_E(\tL)$ is Noetherian and the Theorem is proved.
\end{proof}

Note: The ring $\Diff(E)$ is not almost commutative in general. We get another proof of the almost commutativity of $\Uo(\tL)$:

\begin{corollary} \label{almostcom} For any 2-sided ideal $I\subseteq \Uo(\tL)$, it follows $\Uo(\tL)/I$ is almost commutative. If $A$ is Noetherian and $\tL$ a finitely generated left $A$-module
it follows $\Uo(\tL)$ is Noetherian and almost commutative.
\end{corollary}
\begin{proof} Let $E:=\Uo(\tL)/I$. It follows $E$ is a left $\Uo(\tL)$-module, hence there is a connection $\rho: \tL\rightarrow \End_k(E)$ with $\rho(D)=Id_E$ and $\rho$ an $A\otimes_k A$-linear map.
It follows $ann^{\otimes}(E,\rho)=I$ and $\Uo_E(\tL)\cong \Uo(\tL)/I$ is by Corollary \ref{smainalmost} almost commutative. 
In particular if  $J=(0)$ it follows $\Uo(\tL)$ is almost commutative. There is always a canonical surjective map of graded $A$-algebras
\[ \rho: \sym_A^*(\tL) \rightarrow Gr(\Uo(\tL)) \]
and if $A$ is Noetherian and $\tL$ finitely generated as left $A$-module it follows $\sym^*_A(\tL)$ and $Gr(\Uo(\tL))$ are Noetherian. It follows $\Uo(\tL)$ is Noetherian.
The Corollary follows.
\end{proof}

\begin{example} An application to the study of the Chern classes.\end{example}

In \cite{maa1} the following Theorem is proved:

Let 
\[ k\subseteq A:=U_0 \subseteq U_1\subseteq \cdots \subseteq U_i \subseteq \cdots \subseteq U\]
be an almost commutative PBW-algebra. This means $U$ is an associative unital $k$-algebra with $k$ a commutative unital ring in the center $Z(U)$ of $U$.
The multiplication respects the filtration and the canonical map of graded $A$-algebras
\[ \rho: \sym_A^*(U_1/U_0)\rightarrow Gr(U, U_i) \]
is an isomorphism. The left $A$-module  $L:=U_1/U_0$ is in a canonical way an $A/k$-Lie-Rinehart algebra $(L,\alpha)$ and if $L$ is projective as left $A$-module, there is a 2-cocycle $f\in \Z^2(L,A)$
and an isomorphism
\[  \U(A,L,f)\cong U \]
of filtered associative rings. Hence a class of almost commutative rings may be classified by such pairs $((L,\alpha),f)$. We may ask if it is possible to apply such a classification to the study of the set of 
of 2-sided ideals in $\Uo(\tL)$ and $\Uo_E(\tL)$.

Recall that for a connection $(E,\nabla)$ where $E$ is a rank $r$ projective $A$-module and $\nabla$ has curvature type $f$ for a 2-cocycle $f\in \Z^2(L,A)$ it follows the $k$'th Chern class is determined
by the first Chern class. If $A$ contains a field of characteristic zero we get the formula
\[ c_k(E)=\frac{1}{r^{k-1}}c_1(E)^k.\]
Hence since $\Uo(\tL)$ is almost commutative we may get a condition on the Chern classes of all $\Uo(\tL)$-modules: In some cases all higher Chern classes are determined by $c_1(E)$. If there is a PBW-theorem
for $\Uo(\tL)$ in the case when $\tL$ is a finite rank projective left  $A$-module we may get such results. Such a PBW-theorem is the topic of a forthcoming paper.

\begin{example} An application to the study of the Hodge conjecture.\end{example}

This may have interest for people studying the cycle map, Chern character and the Hodge conjecture. If $X$ is a smooth projective variety there is the Chern character
\[ Ch :\operatorname{K}_0(X)\otimes \mathbb{Q} \rightarrow \H^*(X, \mathbb{C}) \]
where $\H^*(X,\mathbb{C})$ is singular cohomology of $X$ with complex coefficients and $\operatorname{K}_0(X)\otimes \mathbb{Q}$ is the Grothendieck group of $X$ with rational coefficients.
The $p$'th Chern class of a rank $r$ vector bundle $E$ satisfies
\[ c_p(E) \in \H^{p,p}(X)\cap \H^{2p}(X, \mathbb{Z}).\]
If the formula $c_p(E)=\frac{1}{r^{p-1}}c_1(E)^p$ holds for all vector bundles $E$, this may give information on the Hodge conjecture. If we can find elements in  the group 
\[  \H^{p,p}(X)\cap \H^{2p}(X, \mathbb{Z})  \]
that do not lie in the ring generated by $\H^2$ we may get a counterexample to the conjecture. Chose an element $c\in \H^{p,p}(X)\cap \H^{2p}(X, \mathbb{Z})$ that is not in the subring generated by $\H^2$.
Chose any rank $r$ vector bundle $E$ on $X$. It follows $c_p(E)\neq c$. The image of the Chern character equals the image of the cycle map, hence for any cycle  $z\in \operatorname{CH}^*(X)\otimes \mathbb{Q}$
there is a finite rank vector bundle $E$ on $X$ with $Ch(E)=\gamma(z)$. It follows that for any cycle $z\in\operatorname{CH}^*(X)\otimes \mathbb{Q}$ we have $\gamma(z)\neq c$.

Note: When $Y$ is a smooth projective complex variety the Atiyah class
\[ a(E)\in \Ext^1_{\O_Y}(E, \Omega^1_Y \otimes E) \]
is the obstruction for $E$ to have an algebraic (or equivalentely holomorphic) connection and $a(E)$ is "seldom" zero. Hence $E$ seldom has an algebraic connection. The Atiyah class $a(E)$
and its powers $a_p(E)$ can be used to construct to the Chern classes of $E$. Hence the projective situation differs from the affine situation: For a finite rank projective $A$-module $E$ the Atiyah class
\[ a(E)\in \Ext^1_A(E,\Omega^1_A \otimes E) \]
is always zero. Hence there is always an algebraic connection $\nabla: E\rightarrow \Omega^1_A \otimes E$ in the affine situation.

If $E$ is a rank $m$ locally trivial $\O_Y$-module we may always define characteristic classes $Ch_k(E)\in \H^{2k}(Y, \mathbb{C})$  as follows:

\begin{definition}\label{cherncharacter}
Let $Ch_1(E):=c_1(E)$ and let for $k\geq 2$
\begin{align}
&\label{cherncl} Ch_k(E):=\frac{1}{k!}\frac{1}{m^{k-1}}c_1(E)^k\in \H^{2k}(Y,\mathbb{C}) 
\end{align}
be the \emph{$k$'th Chern character} of $E$, where $\H^{2k}(Y,\mathbb{C})$ is singular cohomology of $Y$ with complex coefficients. 
Let
\[ Ch(E):=\sum_{k\geq 0}Ch_k(E)\]
be the \emph{Chern character} of $E$.
\end{definition}

It follows $f^*Ch_k(E)=Ch_k(f^*E)$ for any map $f: X\rightarrow Y$ of smooth projective varieties.
Moreover $Ch_1(E\otimes F)=c_1(E\otimes F)=nc_1(E)+nc_1(F)=nCh_1(E)+mCh_1(F)$.

%If the class $c^*_i(E)$ defines a "Chern class theory" and there is a "unicity theorem"
%for such Chern classes with values in singular cohomology, it would follow that the image of the Chern character 
%(and cycle map) is in the subring generated by $\H^{1,1}(Y)$, since all higher Chern classes of $E$ are determined by $c_1(E)$. 

\begin{lemma} For two locally trivial $\O_Y$-modules $E,F$ of rank $m,n$ the following holds for all $k\geq 1$:
\[ Ch_k(E\otimes F)=\sum_{i+j=k}Ch_i(E)Ch_j(F) \]
\end{lemma}
\begin{proof} We get
\[ Ch_k(E\otimes F):=\frac{1}{k!}\frac{1}{(mn)^{k-1}}c_1(E\otimes F)^k=\]
\[ =\frac{1}{k!}\frac{1}{(mn)^{k-1}}(nc_1(E)+mc_1(F))^k=\]
\[ \frac{1}{k!}\frac{1}{(mn)^{k-1}}\sum_{i=0}^k \binom{k}{i}n^ic_1(E)^im^{k-i}c_1(F)^{k-i}=\]
\[\frac{1}{(mn)^{k-1}}\sum_{i=0}^k \frac{1}{i!}n^ic_1(E)^i\frac{1}{(k-i)!}m^{k-i}c_1(F)^{k-i}=\]
\[\frac{1}{(mn)^{k-1}}\sum_{i=0}^k(mn)^{k-1}Ch_i(E)Ch_{k-i}(F)=\sum_{i+j=k}Ch_i(E)CH_j(F).\]
\end{proof}

We get a map
\[Ch: \oplus_{E} \mathbb{Z}[E] \rightarrow \H^{2*}(Y, \mathbb{Q}) \]
defined by
\[ Ch(\sum_k n_k[E_k]):=\sum_k n_kCh(E_k) .\]
The map $Ch$ is a functorial map and multiplicative with respect to tensor product of locally free sheaves:
\[ Ch(E\otimes F):=\sum_{k\geq 0}Ch_k(E\otimes F)=\]
\[ \sum_{k\geq 0}\sum_{i+j=k}Ch_i(E)Ch_j(F) =\]
\[ (\sum_{i\geq 0}Ch_i(E))(\sum_{j\geq 0} Ch_j(F)) =Ch(E)Ch(F).\]
If the map $Ch$ factored through the Grothendieck group $K_0(Y)$, we would get a theory of Chern classes defined in terms of the first Chern class $c_1$.
The image of the Chern character $Ch$ from Definition \ref{cherncharacter}  is in the sub ring generated by $\H^{1,1}\subseteq \H^{2*}(Y,\mathbb{C})$. 

Let $E$ be any locally free rank $e$ $\O_Y$-module and define for any open subset $U\subseteq Y$ the cohomology class
\begin{align}
&\label{akeu} a_k(E, U):=c_k(E_U)-\frac{1}{e^{k-1}}c_1(E_U)^k \in \H^{2k}(U,Q) 
\end{align}
where $E_U$ is the restriction of $E$ to $U$. Let $U_i$ be an open affine cover of $Y$, with $p_i: U_i \rightarrow Y$ the inclusion map.
If $p_i^*:\H^{2k}(Y,Q)\rightarrow \H^{2k}(U_i,Q)$ is the pull back map, it follows by functoriality of the Chern class construction that
\[ p^*_i(a_k(E,Y))=a_k(E, U_i).\]

%If in general a class $\alpha\in \H^{2k}(Y,Q)$ satisfies $p_i^*(\alpha)=0$ for all $i$ then it followed $\alpha=0$, we could argue as follows:
%If the class $p_i^*(a_k(E,Y))=a_k(E, U_i)=0$ for all $i$ it would follow $a_k(E,Y)=0$. As a consequence we would get the formula
%$c_k(E)=\frac{1}{e^{k-1}}c_1(E)^k$ for any $E$ on $Y$. If this argument holds, it would follow the Chern character $Ch$ from Definition \ref{cherncharacter} would factor 
%through the Grothendieck group.

\begin{proposition} \label{prop} Let $Y$ be a smooth projective variety and let $E$ be a rank $e$ locally trivial $\O_Y$-module. There is an equality
$a_k(E,Y)=0$  if and only if for any open affine set $U\subseteq Y$ there is an equality $a_k(E,U)=0$. 
\end{proposition}
\begin{proof} One implication follows from functoriality. Assume $a_k(E,U)=0$ holds for all open affine sets $U\subseteq Y$ and let $U_i$ be a finite affine open cover of $Y$. It follows 
$a_k(E,U_i)=0$ holds for all $i$. Let $p_i: U_i \rightarrow Y$. It follows $a_k(E,U_i)=p_i^*a_k(E,Y):=a_k(E,Y)_{U_i}=0$ and since the class $a_k(E,Y)$ has the property that  $a_k(E,Y)_{U_i}=0$ for all $i$ 
where $U_i$ is an open cover of $Y$, it follows $a_k(E,Y)=0$ and the Proposition follows.
\end{proof}

\begin{example} A relationship to the universal ring $\Uo(\tL)$.\end{example}

Hence if we can prove that $a_k(E,U_i)=0$ for an affine open cover $U_i$ of $Y$, it follows from Proposition \ref{prop} the image of the Chern character is in the subring generated by $H^{1,1}$.
A PBW-theorem for the universal ring $\Uo(\tL)$ may imply such a result. Recall that there is for any $\Dl$-Lie algerbra $\tL$ a canonical surjective map of graded $A$-algebras
\begin{align}
&\label{rho}  \rho: \sym_A^*(L) \rightarrow Gr(\Uo(\tL)) 
\end{align}
where $L:=\Uo(\tL)^1/\Uo(\tL)^0$. If $\tL$ is projective as left $A$-module we may ask if $\rho$ is an isomorphism. If this is the case, there is an $A/k$-Lie-Rinehart algebra $(L, \alpha)$ and a 2-cocycle $f\in \Z^2(L,A)$
and an isomorphism $\U(A,L,f)\cong \Uo(\tL)$ of filtered associative rings. It follows any left $\Uo(\tL)$-module $E$ corresponds to an $L$-connection $\nabla: L\rightarrow \End_k(E)$
of curvature type $f$. If $k$ contains a field of characteristric zero and $E$ is a finitely generated and projective $A$-module of rank $e$, it follows $c_k(E)=\frac{1}{e^{k-1}}c_1(E)^k$. 
Hence a PBW theorem for $\Uo(\tL)$ may have applications to the study of the Chern character, the cycle map and the Hodge conjecture.

\begin{example} Connections that are not of curvature type $f$.\end{example}

There are examples of vector bundles $\E$ of rank $e$ on smooth complex surfaces $S$ where $c_1(\E)=0$ and where $c_2(\E)\neq 0$. Hence there is an open affine subscheme
$U:=\Spec(A)\subseteq S$ with $c_1(\E_U)=0$ and $c_2(\E_U)\neq 0$. Hence the rank $e$ projective $A$-module $E:=\E_U$ has a connection
\[ \nabla: \Der_{\mathbb{C}}(A)\rightarrow \End_{\mathbb{C}}(E) ,\]
and $\nabla$ does not have curvature type $f$ for a 2-cocycle $f\in \Z^2(\Der_{\mathbb{C}}(A),A)$. Hence the map $\rho$ in \ref{rho} is not an isomorphism in general.
One wants to investigate the relationship between $\rho$ and properties of the Chern classes of vector bundles $\E$ on $S$.

%\begin{example} Almost commutative and Noetherian quotients of the tensor algebra.\end{example}

%Note: From Corollary \ref{almostcom} it follows $\Uo(\tL)$ is an almost commutative ring. Since $\Uo(\tL):=\T^*_k(\tL)/J$ where $J$ is some 2-sided ideal in $\T^*_k(\tL)$, we would like to classify
%2-sided ideals $I$ in the tensor algebra  $\T^*_k(\tL)$ with the property that the quotient $\T^*_k(\tL)/I$ is almost commutative (or Noetherian). The tensor algebra is not Noetherian in general.

%\begin{example} Almost commutative 2-sided ideals in the universal ring $\Uo(\tL)$.\end{example}
%We would like to classify 2-sided ideals in $\Uo(\tL)$ with the property that the quotient $\Uo(\tL)/J$ is almost commutative. Since $\Uo(\tL):=\T^*_k(\tL)/J_1$ we would also like to classify
%2-sided ideals $I$ in the tensor algebra $\T^*_k(\tL)$ with $\T^*_k(\tL)/I$ almost commutative. By Example \ref{idealsinuniversal} it follows such a classification
%will give information on the study of the curvature of a connection.

\begin{example} \label{characteristic} Methods from the theory of rings of differential operators and D-modules. \end{example}

There is an extensive theory of flat connections on complex manifolds, the Riemann-Hilbert correspondence and
modules on rings of differential operators and D-modules. See \cite{borel} and \cite{ginzburg} for an introduction to the subject with references.

The universal ring $\tUo(\tL,\rho)$ defined above is an almost commutative Noetherian ring in many cases. Left and right $\tUo(\tL,\rho)$-modules that are finitely generated over $\tUo(\tL,\rho)$  have many properties
similar to D-modules as studied in \cite{ginzburg}.  We can define the characteristic variety of $(\rho, E)$ using filtrations coming from $\tUo(\tL,\rho)$. 
We may use such filtrations to define the notion of holonomiticy for non-flat connections.  Definition 1.1.11 in \cite{ginzburg} is  algebraic and the only assumption is that $\tUo(\tL,\rho)$ is defined over a field of characteristic zero. Once the characteristic variety $\SS(\rho, E)$ is defined we may define holonomicity as in \cite{ginzburg}. It remains to find out if this leads to a reasonable theory where one can calculate explicit examples. 
Since the universal ring $\tUo(\tL,\rho)$ is Noetherian when $E,\tL$ are finitely generated as $A$-modules it might be worthwile to investigate this.

   In \cite{ginzburg} one uses localization for non-commutative rings to prove that results for modules on
almost commutative Noetherian rings globalize to give constructions valid for modules on sheaves of differential operators on complex manifolds. 
Similar methods can be used to prove that the construction of the universal ring $\tUo(\tL,\rho)$  globalize to give a construction for arbitrary schemes. Since the cohomology and homology of a connection
is defined using the theory of modules over associative rings, the theory uses methods from non-commutative algebra/algebraic geometry and the theory of 
sheaves of rings of differential operators and jet bundles.

Recall the following theorem

\begin{theorem} \label{poisson} Let $\U$ be a Noetherian  almost commutative associative unital ring and $M$ a finitely generated $\U$-module. The characteristic variety
\[ \SS(M) \subseteq \Spec(Gr(U)) \]
is coisotropic with respect to the Poisson structure on $Gr(\U)$.
\end{theorem}
\begin{proof} See \cite{ginzburg}, Theorem 1.2.5.
\end{proof}

In the case above where $A$ is Noetherian, $\tL,E$ finitely generated as left $A$-modules, it follows $\tUo(\End(\tL,E))$ and $\tUo(\tL,\rho)$ are almost commutative Noetherian rings.
Hence Theorem \ref{poisson} applies to any finitely generated $\tUo(\tL,\rho)$ or $\tUo(\End(\tL,E))$-module. The connection $(\rho, E)$ is a finitely generated left $\tUo(\tL,\rho)$-module
and we may construct the characteristic variety $\SS(\rho, E)$. One wants to study the ring $\tUo(\tL,\rho)$, the variety $\SS(\rho, E)$ and its relationship
to the connection $\rho$ and the Chern classes $c_i(E)$ for non-flat connections $(\rho, E)$.

%\begin{corollary} Assume the ring $k$ contains a field of characteristic zero and assume $(E,\rho)$ is a connection $\rho: \Der_k(A)\rightarrow \End_k(E)$. Let $L:=\Der_k(A)$ 
%and let $f\in \Z^2(L,A)$ be a 2-cocycle. If there is an inclusion of  ideals $I(f) \subseteq ker(\T(\rho))$ where
%\[  \T(\rho): \T^*_B(L(0)) \rightarrow \Diff(E) \]
%is the map induced by $\rho$ it follows there is an equality of Chern classes 
%\[  c_k(E)=\frac{1}{(k-1)!}c_1(E)^k \]
%in $\H^{2k}(L,A)$.
%\end{corollary}
%\begin{proof} The claim follows from Proposition \ref{fixed} and Corollary \ref{det}
%\end{proof}

%Hence the universal enveloping algebra  $\U(A,L,f)$ defined and studied in \cite{maa1} is Morita equivalent to the 
%associative ring $\T^*_{I(\fal)}(L(0))$ for all 2-cocycles $\fal$. Hence we may for all $\fal$ construct the category
%$\Mod(\U(A,L,f)) $ using quotients of the tensor algebra $\T^*_B(L(0))$ of the trivial abelian extension $L(0)$ of $\operatorname{D}^1(A,0)$. Hence the module categories
%for the family of algebras $\U(A,L,f)$ may be constructed using quotients of one fixed tensor algebra $\T^*_B(L(0))$ using the family of 2-sided ideals $I(f)$ constructed in Proposition \ref{fixed}. 

%When $A$ is Noetherian and $L$ a finitely generated and
%projective $A$-module, it follows $\U(A,L,f)$ is Noetherian. The proof uses the PBW-theorem for $\U(A,L,f)$ (see \cite{maa1} for a proof). 

Sridharan studied the enveloping algebra $\U(k,\lg,f)$ in \cite{sridharan} for $k$ a fixed commutative ring, $\lg$ a Lie algebra over $k$ and $f$ a 2-cocycle for $\lg$.  
He gave a complete description of the deformation groupoid of $\lg$ in the case when $\lg$ is a $k$-Lie algebra with a basis  as $k$-module.
Rinehart studied in \cite{rinehart} the enveloping algebra $\U(A,L)$ for an arbitrary  Lie-Rinehart algebra and proved the PBW-teorem for $\U(A,L)$ in the case when $L$ is a projective $A$-module.
He used this theorem to study cohomology and homology of $L$-connections. 
Tortella gave in \cite{tortella} a simultaneous generalization of the construction of Sridharan and Rinehart for holomorphic Lie-algebroids on complex manifolds and proved a PBW-theorem for the sheaf of enveloping algebras of such holomorphic Lie algebroids. In \cite{maa1} I gave an algebraic construction of the enveloping algebra $\U(A,L,f)$ for any Lie-Rinehart algebra $L$ and any 2-cocycle $f$. I also gave some algebraic proofs of results in Tortellas paper and a proof of the PBW-theorem for $\U(A,L,f)$ in the case when $L$ is a projective $A$-module. I used the PBW-theorem to give a complete determination of the deformation groupoid of 
$(L,\alpha)$ in the case when $L$ is a projective $A$-module. Note that the paper \cite{maa1} contains some errors in the section on families of connections (Definition 5.1).
In \cite{maa141} I give a classification of all non-abelian extensions of a $\Dl$-Lie algebra $(\tL, \ta, \tp, [,], D)$ by an $A$-Lie algebra $(W,[,])$
with $aw=wa$ for all $a\in A$ and $w\in W$ such that $\tL$  is projective as left $A$-module. This classification generalize the classification given in \cite{tortella1} to $\Dl$-Lie algebras. In \cite{tortella1}
the classification is done for Lie-Rinehart algebras.

\end{document}